\documentclass[10pt]{article}

\usepackage{latexsym}
\usepackage{amsfonts}
\usepackage{amsmath}
\usepackage{amsthm}
\usepackage{amssymb}
\usepackage{mathrsfs}
\usepackage{dsfont}    
\usepackage{bbold}     
\usepackage[english]{babel}
\usepackage{caption}
\usepackage{epsfig}
\usepackage{float}
\usepackage{psfrag}
\usepackage{graphicx}
\usepackage{epsfig}
\usepackage{hyperref}
\newtheorem{theorem}{Theorem}
\newtheorem{defi}[theorem]{Definition}
\newtheorem{lemma}[theorem]{Lemma}
\newtheorem{proposition}[theorem]{Proposition}
\newtheorem{rmk}[theorem]{Remark}
\newtheorem{hypo}{Hypothesis}

\newcommand{\HH}{{\mathcal H}}
\newcommand{\ZZZ}{\mathds{Z}}
\newcommand{\CCC}{\mathds{C}}
\newcommand{\NNN}{\mathds{N}}

\newcommand{\RRR}{\mathds{R}}
\newcommand{\TTT}{\mathds{T}}
\newcommand{\uno}{\mathds{1}}

\newcommand{\calA}{{\mathcal A}}

\newcommand{\CCCC}{{\mathcal C}}
\newcommand{\DD}{{\mathcal D}}

\newcommand{\calG}{{\mathcal G}}
\newcommand{\calH}{{\mathcal H}}
\newcommand{\calI}{{\mathcal I}}

\newcommand{\calL}{{\mathcal L}}

\newcommand{\calQ}{{\mathcal Q}}
\newcommand{\RR}{{\mathcal R}}

\newcommand{\TT}{{\mathcal T}}
\newcommand{\calU}{{\mathcal U}}
\newcommand{\VV}{{\mathcal V}}

\newcommand{\gots}{{\mathfrak s}}

\newcommand{\ol}{\overline}

\newcommand{\Fullbox}{{\rule{2.0mm}{2.0mm}}}

\newcommand{\EP}{\hfill\Fullbox\vspace{0.2cm}}
\newcommand{\prova}{\noindent{\it Proof. }}

\newcommand{\e}{\varepsilon}
\newcommand{\al}{\alpha}
\newcommand{\de}{\delta}
\newcommand{\be}{\beta}
\newcommand{\z}{\zeta}
\newcommand{\n}{\nu}
\newcommand{\m}{\mu}
\newcommand{\x}{\xi}

\newcommand{\ka}{\kappa}
\newcommand{\g}{\gamma}
\newcommand{\om}{\omega}
\newcommand{\h}{\eta}
\newcommand{\la}{\lambda}

\newcommand{\f}{\varphi}
\newcommand{\s}{\sigma}
\newcommand{\id}{\mathbb{1}}
\newcommand{\del}{\partial}

\newcommand{\oo}{\omega}

\newcommand{\ff}{\mathtt{f}}

\newcommand{\CC}{\boldsymbol{C}}

\newcommand{{\resonance}}{relevant self-energy cluster }

\def\ins#1#2#3{\vbox to0pt{\kern-#2 \hbox{\kern#1 #3}\vss}\nointerlineskip}

\newcommand{\te}{\theta}

\makeindex
\begin{document}



\title{\bf Quasi-periodic solutions for fully nonlinear forced reversible Schr\"odinger equations}

\author
{\bf Roberto Feola, Michela Procesi 
\vspace{2mm}
\\ \small 
Address: Universit\`a di Roma La Sapienza, P.le A. Moro 5, 00185 Roma
\\ \small 
E-mail: feola@mat.uniroma1.it, mprocesi@mat.uniroma1.it}


\date{}

\maketitle

\begin{abstract}
In this paper we consider a class of fully nonlinear forced and reversible Schr\"odinger equations and prove existence and stability of quasi-periodic solutions. We use a  Nash-Moser algorithm together with a reducibility theorem on the linearized operator in a neighborhood of zero. Due to the presence of the highest order derivatives in the non-linearity the classic KAM-reducibility argument fails and one needs to use a wider class of changes of variables such has diffeomorphisms of the torus and pseudo-differential operators.  This procedure automtically produces a change of variables, well defined on the phase space of the equation, which diagonalizes the operator linearized at the solution. This gives the linear stability. 

\end{abstract}

\tableofcontents

\setcounter{equation}{0}
\section{Itroduction}
\label{sec:1}
In this paper we study a class of reversible forced fully non linear Schr\"odinger equations of the form
\begin{equation}\label{mega}
iu_{t}=u_{xx}+\e 	\ff(\oo t,x,u,u_{x},u_{xx}),
\quad x\in\TTT:=\RRR/2\pi\ZZZ, 
\end{equation}
where $\e>0$ is a small parameter, the nonlinearity is quasi-periodic in time with
diophantine frequency vector $\oo\in \RRR^{d}$
and $\ff(\f,x,z)$, with $\f\in\TTT^{d}$, $z=(z_{0},z_{1},z_{2})\in\CCC^{3}$
is in $C^{q}(\TTT^{d+1}\times\CCC^{3};\CCC)$ in the real sense (i.e. as function of ${\rm Re}(z)$ and ${\rm Im}(z)$).
For this class we prove
existence and stability of quasi-periodic solutions with Sobolev regularity for
all $\la$ in an appropriate positive measure Cantor-like set. 
The study of this kind of solutions for the ``classic'' autonomous 
semi-linear NLS (where the nonlinearity $\ff$
does not contain derivatives) was one of the first successes of KAM theory for Pde's. See for instance
\cite{W1, K1, KP}. More recently a series of papers have appeared concerning dispersive
semi-linear Pde's where the nonlinearity contains derivatives of order $\de\leq n-1$,  here $n$ is
the order of the highest derivative appearing in the linear constant coefficients term.
We mention in particular  \cite{ZGY} for the reversible NLS
 and \cite{LY} for the Hamiltonian case.  The key point of 
 the aforementioned papers is to apply KAM theory
 by using an appropriate generalization of the so called
 "\emph{Kuksin Lemma}''. This idea has been introduced by Kuksin in \cite{Ku2} 
 to deal with non-critical unbounded perturbations, i.e. $\de<n-1$, with the purpose of studying KdV type equations, see also \cite{KaP}. The previously mentioned results require that the equation is semi-linear and dispersive; in the "weakly dispersive" case of the derivative Klein-Gordon equation we mention the results \cite{BBiP1}-\cite{BBiP2}, also based on KAM theory.
Note that our equation  is   fully nonlinear, namely
 the second spatial derivative appears also in the nonlinearity. Hence the KAM approach seems to fail and one has to develop different strategies. The first breakthrough result for fully nonlinear Pde's is due to Iooss-Plotnikov-Toland
  who studied in  \cite{IPT} the existence of periodic solutions for water-waves; we mention also the papers
  by Baldi \cite{Ba1}, \cite{Ba2} on periodic solutions for the Kirchoff and Benjamin-Ono equations. 
  These papers are based on Nash-Moser methods and the key point is to apply appropriate diffeomorphism of the torus  and pseudo-differential operators in order to invert the operator linearized at an approximate solution. Note that these results do not cover the linear stability of the solutions and they do not work in the quasi-periodic case. Quite recently this problem has been overcome by Berti, Baldi, Montalto who studied 
 fully nonlinear perturbations of the KdV equation first in \cite{BBM1}, for the forced case, then in \cite{BBM2} for the autonomous. This was the first result for quasi-periodic solutions for quasi linear Pde's and the main purpose of the present paper is to generalize their strategy to cover the NLS. 
 
 \smallskip
 
 We now give an overview of the main problems and strategies which appear in the study of quasi-periodic solutions in Pde's referring always, for simplicity, to the forced case. 
  A quasi-periodic solution, with frequency $\oo\in\RRR^{d}$, for an equation such as \eqref{mega} is a function 
 of the form $\mathtt{u}(t,x)=u(\oo t,x)$ where 
 $$
u(\f,x) : \TTT^{d}\times\TTT\to\CCC.
$$
 In other words we look for non-trivial 
$(2\pi)^{d+1}-$periodic solutions $u(\f,x)$ of
\begin{equation}\label{mega2}
i\oo\cdot\del_{\f}u=u_{xx}+\e  \ff(\oo t,x,u,u_{x},u_{xx})
\end{equation}
in the Sobolev space 
\begin{equation}\label{1.1}
H^{s}:=H^{s}(\TTT^{d}\times\TTT; \CCC):=
\{u(\f,x)=\!\!\!\sum_{(\ell,k)\in\ZZZ^{d}\times\ZZZ}\!\!\!u_{\ell,k}e^{i(\ell\cdot\f+k\cdot x)} :
||u||^{2}_{s}:=\sum_{i\in\ZZZ^{d+1}}|u_{i}|^{2}\langle i\rangle^{2s}<+\infty
\}.
\end{equation}
where $s>\gots_{0}:=(d+2)/2>(d+1)/2$, 
$i=(\ell,k)$ and $\langle i\rangle:=\max(|\ell|, |k|,1)$, $|\ell|:=\max\{|\ell_{1}|,\ldots,
|\ell_{n}|\}$. For $s\geq\gots_{0}$ $H^{s}$ is a Banach Algebra and 
$H^{s}(\TTT^{d+1})\hookrightarrow C(\TTT^{d+1})$ continuously. As in \cite{BBM1} we consider the frequency vector 
\begin{equation}\label{dio}
\oo=\la \bar{\oo}\in \RRR^{d}, \quad \la\in\Lambda:=\left[\frac{1}{2},\frac{3}{2}\right], \quad
|\bar{\oo}\cdot\ell|\geq \frac{\g_{0}}{|\ell|^{d}}, \quad \forall\; \ell\in\ZZZ^{d}\backslash\{0\}.
\end{equation}
We impose the {\em reversibility condition}
\begin{hypo} \label{hyp2}
Assume that $\ff$ is such that

\begin{itemize}

\item[(i)] $\ff(\f,-x,-z_{0},{z}_{1},-z_{2})=
-\ff(\f,x,z_{0},{z}_{1},z_{2})$.

\item[(ii)] $\ff(-\f,x,z_{0},z_{1},z_{2})=
\ol{{\ff}(\f,x,\bar{z}_{0},\bar{z}_{1},\bar{z}_{2})}$,
\item[(iii)] $
\ff(\f,x,0)\neq 0\, ,\quad  \del_{z_{2}}\ff\in\RRR\backslash\{0\} $,

\noindent where $\del_{z}= \del_{\rm{Re}(z)}-i\, \del_{\rm{Im}(z)}$.
\end{itemize}
\end{hypo}
Our main result is stated in the following:

\begin{theorem}\label{teo1}
There exists $s:=s(d)>0$, $q=q(d)\in\NNN$ such that
for every nonlinearity $\ff\in C^{q}(\TTT^{d+1}\times\CCC^{3};\CCC)$ that satisfies
Hypothesis \ref{hyp2} and  for all $\e\in(0,\e_{0})$, with $\e_{0}=\e_{0}(\ff,d)$
small enough, there exists a Cantor set $\CCCC_{\e}\subset\Lambda$ of asymptotically full Lebesgue measure, i.e.
\begin{equation}\label{asy}
|\CCCC_{\e}|\to 1 \quad  as \quad \e\to0,
\end{equation}
such that for all $\lambda\in\CCCC_{\e}$ the perturbed NLS equation (\ref{mega2})
has solution $u(\e,\lambda)\in H^{s}$ such that $u(t,x)=-\bar{u}(-t,-x)$, 
with $||u(\e,\lambda)||_{s}\to0$ as $\e\to0$. In addition, $u(\e,\lambda)$ is {\rm linearly stable}.
\end{theorem}

 Finding such a solution is equivalent to finding zeros of a nonlinear functional on the prescribed Sobolev space.
 In forced cases, the starting point is to consider functionals $F(\la,\e,u)$  that for $\e=0$
 are linear with constant coefficients and have purely imaginary spectrum which accumulate to zero. See \eqref{4'} for the NLS case. 
 Note that this is a perturbative problem since $F(\la,0,0)=0$. However the linearized operator $d_{u}F(\la,0,0)$ is not invertible and one needs to use a generalized Implicit Function Theorem.
 
\begin{figure}[ht] 
\centering 
\ins{260pt}{-20pt}{$F(u)$}
\ins{285pt}{-143pt}{$u_0$}
\ins{245pt}{-143pt}{$u_1$}
\ins{222pt}{-143pt}{$u_2$}
\includegraphics[width=3in]{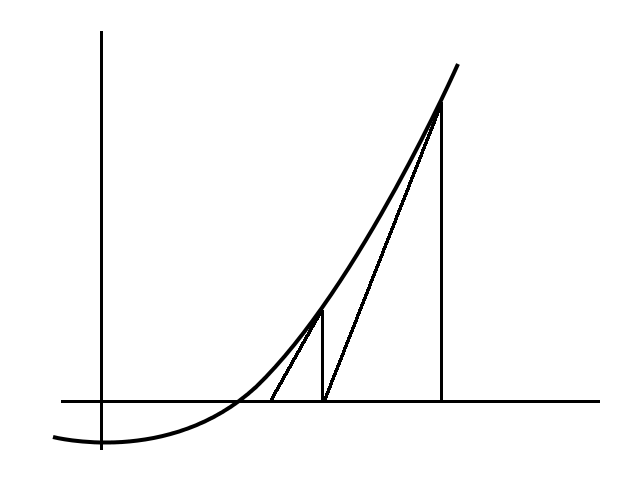} 
\vskip.2truecm 
\caption{Three steps of the Newton algorithm
$u_{n+1} := u_n - (d_{u}F(\la,\e,u_n))^{-1} [F(\la,\e, u_n)] $ }
\label{fig.newton} 
\end{figure}

Typically this method is based on a Newton-like scheme, which  relies on the
invertibility of the linearized equation in a whole neighborhood of the unperturbed solution,
in our case $u=0$; see Figure \ref{fig.newton}.


On a purely formal level one can state an
 abstract ``Nash-Moser'' scheme (see for instance \cite{BBP},\cite{BCP} and our Section \ref{secNM})
which says that if $\la$ is such that for all $n$ the operator $(d_{u}F(\la,\e, u_n))^{-1}$ is well-defined
and bounded from $H^{s+\mu}$ to $H^{s}$ for some $\mu$, then a solution of \eqref{mega2} exists.
Then the problem reduces to proving that such set of parameters $\la$ is non-empty, or even better
that it has asymptotically full measure.

If we impose some symmetry such as a Hamiltonian or a reversible structure  the linearized operator
$d_{u}F(\la,\e, u)$ is self-adjoint and it is easy to obtain lower bounds 
on its eigenvalues,  
implying its invertibility with bounds
on the $L^2$-norm of the inverse for ``most" parameters $ \la $, this is the so called \emph{first Mel'nikov condition}.  However  this information
is not enough to prove the convergence of the algorithm: one needs
estimates on the high Sobolev norm of the inverse, which  do not follow only
from bounds on the eigenvalues. 

Naturally, if $d_{u}F(\la,\e,u)$ were diagonal, passing from $L^2$ to $H^s$ norm would be trivial,
but the problem is that the operator which diagonalizes $d_{u}F(\la,\e,u)$ may not be bounded in
$H^s$. The property of an operator to be diagonalizable via a ``smooth'' change of variables
is known as \emph{reducibility} and in general is connected to the fact that the matrix is regular
semi-simple, namely its eigenvalues are distinct (see \cite{El1}  for the finite dimensional case).
When dealing with infinite dimensional matrices, one also has to give quantitative estimates
on the difference between two eigenvalues: this is usually referred to as the \emph{second order
Mel'nikov condition} (note that this can be seen as a condition on $\la$).  
Naturally one does not need to diagonalize a matrix in order to invert it,
indeed in the case of  Pde's on tori, where the eigenvalues are multiple,  
 the first results have been proved without any reducibility. See for instance Bourgain in \cite{B1,B3,B5}, 
 Berti-Bolle in \cite{BB1,BB2},  Wang \cite{W2}.
 These papers rely on the so called ``multi-scale'' analysis based on first Mel'nikov condition and 
 geometric properties of ``separation of singular sites''.
 Note that this method does not imply reducibility and linear stability of the solutions. Indeed there are very few results on reducibility on tori. We mention Geng-You in \cite{GY} for the smoothing NLS, Eliasson-Kuksin in \cite{EK} 
 for the NLS and Procesi-Procesi \cite{PP1,PP2} for the resonant NLS. All the aforementioned papers, both using KAM or multiscale,  
 are naturally on semi linear Pde's with no derivatives in the non linearity. This problem is at this moment 
 completely open and all the results are in the one dimensional case. Here as we said the first results were obtain by KAM methods using the Kuksin lemma.   Roughly speaking the aim of a reducibility scheme is to iteratively conjugate 
 an operator $D+ \e M$, where $D$ is diagonal, to  $D_++\e^2 M_+$ where $D_+$ is again diagonal. Clearly the conjugating transformation must be bounded. The equation which defines the change of variables is called the {\emph {homological equation} } while the operators $D,D_+$ are called the \emph{normal form}.  When $M$ contains derivatives it turns out that $D_+$ can only be diagonal in the space variable (with coefficients depending on time). The  purpose of the Kuksin Lemma is to show that such an algorithm can  be run,  namely that one can solve the homological equation also when the normal form is diagonal only in the space variable (as is  $D_+$).  Unfortunately, if $M$ has the same order (of derivatives) as $D$ this scheme seems to fail and one is not able to find a bounded solution of the Homological equation. The breakthrough idea taken from pseudo-differential calculus is to conjugate   $D+ \e M$ to an operator $D_++\e M_+$ where $D_+$ is again diagonal while $M_+$ is of lower order w.r.t. $M$. After a finite number of such steps one obtains an operator of the form $D_F+ \e M_F$ where $D_F$ is diagonal and $M_F$ is bounded. At this point one can apply a  KAM  reducibility scheme in order to diagonalize. Note that in principle one needs only to invert $D_F+ \e M_F$ which could be done by a multiscale argument, however since we are working in one space dimension one can show that the second Mel'nikov condition can be imposed. This gives the stronger stability result. This scheme, i.e.  Nash Moser plus reducibility of the linearized operator, is very reminiscent of the classical KAM scheme. The main difference is that  we do not apply the changes of variables that diagonalize the linearized operator. The KAM idea instead, is to change variables at each step in order to ensure that the linearized operator is approximately diagonal (and the solution is approximately at the origin).  Unfortunately this changes of variables destroy the special structure of the linearized operator of a Pde, but this property is strongly needed in the first part of our strategy namely in order to conjugate    $D+ \e M$ to $D_F+ \e M_F$.
 
 Regarding our reversibility condition (actually a very natural condition appearing in various works, starting from Moser \cite{Moser-Pisa-66}) some comments are in order. First of all some  symmetry conditions are needed in order to have existence, in order to exclude the presence of dissipative terms. Also such conditions guarantee that the eigenvalues of the linearized operator are all imaginary. All this properties could be imposed by using a Hamiltonian structure, however preserving the symplectic structure during our  Nash-Moser iteration is not straightforward. 
 Another property which follows by the reversibility is that the spectrum of the operator linearized at zero is simple, this is not true in the Hamiltonian case, see \cite{F}.
 A further step is to consider autonomous equations as done in \cite{BBM2}. In this paper we decided to restrict our attention to the forced case where one does not have to handle the bifurcation equation. In the paper \cite{BB3} (see also \cite{BBM2}) the authors show that one can reduce the autonomous case to a forced one, by choosing appropriate coordinates at each Nash-Moser step.   Since the forced case contains all the difficulties related to the presence of derivatives, we are fairly confident that this set of ideas can be used to cover the case of the autonomous NLS.


\subsection{Notations and scheme of the proof}

{\bf Vector NLS.}
We want to ``double'' the variables and study a ``vector'' NLS.  Let us define
\begin{equation}
{\bf u}:=(u^{+},u^{-})\in
H^{s}\times H^{s}.
\end{equation}
On the space $H^{s}\times H^{s}$ we consider the natural norm $||{\bf u}||_{s}:=\max\{||u^{+}||_{s},||u^{-}||_{s}\}$
(we denote by $||\cdot||_{s}$ the usual Sobolev norm on $H^{s}(\TTT^{d+1};\CCC)$).
 We consider also the {\emph{real}} subspace
\begin{equation}\label{5}
\calU:=\left\{{\bf u}=(u^{+},u^{-})\; : \; \ol{u^{+}}=u^{-}\right\}
\end{equation}
in which we look for the solution.

\begin{defi}\label{maria}
Given $\ff\in C^{q}$, we define the ''vector'' NLS as
\begin{equation}\label{4'}
F({\bf u}):=\oo\cdot\del_{\f}{\bf u}+i(\del_{xx}{\bf u}+\e f(\f,x,{\bf u}))=0, \quad
f(\f,x,{\bf u}):=\left(\begin{matrix} f_{1}(\f,x,u^{+},u^{-},u_{x}^{+},u_{x}^{-},u_{xx}^{+},u_{xx}^{-}) \\ f_{2}(\f,x,u^{+},u^{-},u_{x}^{+},u_{x}^{-},u_{xx}^{+},u_{xx}^{-}) \end{matrix}\right)
\end{equation}
where the functions $f=(f_{1},f_{2})$ extend  $(\ff,\bar{\ff})$ in the following sense.
The $f_{j}$ are in $C^{q}(\TTT^{d+1}\times\RRR^{6}\times\RRR^{6}; \RRR^{2})$, and moreover 
on the subspace $\calU$ they satisfy 
\begin{equation}\label{4bis}
\begin{aligned}
f&=(\ff,\bar{\ff})\\
\del_{z^{+}_{2}}f_{1}=\del_{z^{-}_{2}}f_2, \quad  &\del_{z^{+}_{i}}f_1=\ol{\del_{z^{-}_{i}}f_{2}}, \;\; i=0,1,\;\;\;
\del_{z^{-}_{i}}f_1=\ol{\del_{z^{+}_{i}}f_2}, \quad i=0,1,2,\\ \del_{\overline {z^{+}_{i}}}f_{1}& =\del_{\overline {z^{+}_{i}}}f_{2}=\del_{\overline {z^{-}_{i}}}f_{1} =\del_{\overline {z^{-}_{i}}}f_{2}=0 \;\;
{\rm where}\;\; \del_{\overline {z^{\s}_{j}}}= \del_{{\rm Re}\,z^{\s}_j}+ i \del_{{\rm Im}\,z^{\s}_j}, \;\; \s=\pm,
\end{aligned}
\end{equation} 
Note that this extension is trivial in the analytic case.
\end{defi}

By Definition \ref{maria} the (\ref{4'}) reduces to (\ref{mega2}) on the subspace $\calU$ (see the first line in (\ref{4bis})). The advantage of working on (\ref{4'}) is that the linearized operator 
$dF(\bf u):=\calL({\bf u})$ for  ${\bf u}\in\calU$ is self-adjoint.
Note that the linearized operator of (\ref{mega2}) is actually self-adjoint, but even at $\e=0$ is not diagonal. To diagonalize one needs to complexify and then to give meaning to $f\in C^{q}$, thus 
 we introduce the extension.  
 
 \noindent By Hypothesis \ref{hyp2} one has that \eqref{4'}, restricted to $\calU$,
 is reversible with respect to the involution
\begin{equation}\label{b}
S : u(t,x) \to -\bar{u}(t,-x), \quad S^{2}=\id,
\end{equation}
namely, setting
$V(t,u):=-i(u_{xx}+\e \ff(\oo t,x,u,u_{x},u_{xx}))
$
we have
\begin{equation*}
-SV(-t,u)=V(t,Su).
\end{equation*}
Hence the subspace of ``reversible'' solutions
\begin{equation}\label{reversiblesol}
u(t,x)=-\bar{u}(-t,-x).
\end{equation}
is invariant.
 It is then natural to  look for ``reversible''  solutions, i.e. $u$ which satisfy (\ref{reversiblesol}). To formalize this condition we introduce  spaces of odd or even functions in $x\in \TTT$. For all $s\geq0$, we set
\begin{equation}\label{SPACES}
\begin{aligned}
X^{s}&:=\left\{u\in H^{s}(\TTT^{d}\times\TTT) : \;\;\; u(\f,-x)=-u(\f,x), \;\; u(-\f,x)=\bar{u}(\f,x)
 \right\},\\ 
Y^{s}&:=\left\{u\in H^{s}(\TTT^{d}\times\TTT) : \;\;\; u(\f,-x)=u(\f,x), \;\; u(-\f,x)=\bar{u}(\f,x)
 \right\},\\ 
Z^{s}&:=\left\{u\in H^{s}(\TTT^{d}\times\TTT) : \;\;\; u(\f,-x)=-u(\f,x), \;\; u(-\f,x)=-\bar{u}(\f,x)
 \right\},
\end{aligned}
\end{equation}
Note that reversible solutions means $u\in X^{s}$, moreover an operator reversible w.r.t. the involution $S$ maps $X^s$ to $Z^s$. 
\begin{defi}\label{stizzi}
 We denote with bold symbols the spaces ${\bf G^{s}}:=G^{s}\times G^{s}\cap \mathcal U$ where
$G^{s}$ is $H^s$, $X^{s},Y^{s}$ or $Z^{s}$.

\noindent We denote by  $H_{x}^{s}:=H^{s}(\TTT)$ the Sobolev spaces of functions of $x\in\TTT$ only, same for all the subspaces $G^s_x$ and ${\bf G}^s_x$.
\end{defi}
 
 \begin{rmk}\label{rmkphase}
Given a family of linear operators $A(\f) : H_{x}^{s}\to H_{x}^{s}$ for $\f\in\TTT^{d}$, we can associate it to an operator $A : H^{s}(\TTT^{d+1})\to H^{s}(\TTT^{d+1})$ by considering each matrix elements of $A(\f)$ as a multiplication operator. This identifies a subalgebra of linear operators on $H^{s}(\TTT^{d+1})$. An operator 
$A$ in the subalgebra identifies uniquely its corresponding ``phase space'' operator $A(\f)$. 
With reference to the Fourier basis this sub algebra is called ``T\"opliz-in-time'' matrices (see formul{\ae}  \eqref{eq:2.16}, \eqref{pozzo10}).

\end{rmk}

\begin{rmk}\label{rmkspace}
Part of the proof is to control that, along the algorithm, the operator $d_u F(\lambda,\e,u)$ maps
the subspace ${\bf X}^{0}$  into ${\bf Z}^{0}$. In order to do this, we will introduce the notions of 
''reversible'' and ``reversibility-preserving'' operator in the next Section.
\end{rmk}

The proof is based on four main technical propositions.
First we apply an (essentially standard) Nash-Moser iteration scheme which produces a Cauchy sequence of functions converging to a  {\em solution} on a possibly empty Cantor like set.
\begin{proposition}\label{teo4}
Fix $\g\leq\g_{0}, \mu>\tau> d$. There exist $q\in\NNN$, depending only on $\tau,d,\mu$, such that for any nonlinearity
$\ff\in C^{q}$ satisfying Hypothesis \ref{hyp2} the following holds. Let $F({\bf u})$ be defined in Definition
\ref{maria}, then there exists a small constant $\epsilon_0>0$ such that
for any $\e$ with
 $0<\e\g^{-1}<\epsilon_0$, there exist constants $C_{\star}, N_0\in \NNN$, a sequence of functions ${\bf u}_{n}$ and a sequence of sets $\calG_{n}(\g,\tau,\mu)\equiv\calG_{n}\subseteq\Lambda$ such that
\begin{equation}\label{teo41}
{\bf u}_{n} : \calG_{n}\to{\bf X}^{0}, \quad ||{\bf u}_{n}||_{\gots_{0}+\mu,\g}\leq1, \quad
||{\bf u}_{n}-{\bf u}_{n-1}||_{\gots_{0}+\mu,\g}\leq C_{\star}\e \g^{-1} (N_{0})^{-\left(\frac{3}{2}\right)^{n}(18+2\mu)}.
\end{equation}
Here $||\cdot||_{s,\g}$ is an appropriate weighted Lipschitz norm, see \eqref{lnorm}.
Moreover the sequence converges in $||\cdot||_{\gots_{0}+\mu,\g}$ to a function ${\bf u}_{\infty}$
such that
\begin{equation}\label{teo42}
F({\bf u}_{\infty})=0, \quad \forall\; \lambda\in\calG_{\infty}:=\cap_{n\geq0}\calG_{n}.
\end{equation}
\end{proposition}
In the Nash-Moser scheme the main point is to invert, with appropriate bounds, $F$ linearized at any ${\bf u}_n$. 
Following the classical Newton scheme we define 
\begin{equation}\label{pozzoprofondo}
{\bf u}_{n+1}= {\bf u}_n -\Pi_{N_{n+1}}\calL^{-1}({\bf u}_n)\Pi_{N_{n+1}}F({\bf u}_{n}), \quad \calL({\bf u}):=d_{{\bf u}}F({\bf u})
\end{equation}
where $\Pi_{N}$ is the projection on trigonometric polynomials of degree $N$ and
$N_{n}:=(N_{0})^{\left(\frac{3}{2}\right)^{n}}$.
In principle we do not know whether this definition is well posed since $\calL({\bf u})$ 
may not be invertible.
Thus one introduce $\calG_{n}$ as the set where such inversion is possible and bounded in high Sobolev norms. Unfortunately this sets are often difficult to study.  In order to simplify this problem we prove that $\calL({\bf u})$ can be diagonalized in a whole neighborhood of zero.
A major point is to prove that the diagonalizing changes of variables are bounded in high Sobolev norms.
This reduction procedure is quite standard when the non linearity $\ff$ does not contain derivatives. 
In this simpler case $\calL({\bf u})$ is a diagonal matrix plus a small bounded perturbation.
In our case this is not true,  indeed
\begin{equation}\label{1.2.1}
\begin{aligned}
\calL({\bf u}) &= \oo\cdot\del_{\f}\id+i(\id+  A_{2}(\f,x))\del_{xx}+
i A_{1}(\f,x)\del_{x}+i  A_{0}(\f,x)\\
\end{aligned}
\end{equation}
where $A_{i} : {\bf H}^{s}\to {\bf H}^{s}$ are defined in (\ref{eq:3.2}) and $\id$ is the $2\times2$ identity.
Hence the reduction requires a careful analysis which we perform in Sections \ref{sec:3} and \ref{sec:4}.  More precisely in Section \ref{sec:3} we perform a series of changes of variables which conjugate $\calL$ to an operator $\calL_4$ which is the sum of an unbounded {\em  diagonal} operator plus a small \emph{bounded} remainder. Then in section \ref{sec:4} we perform a KAM reduction algorithm. Putting this two steps together, in Section \ref{sec:5} we obtain:
\begin{proposition}\label{teo2}
Fix $\g\leq\g_{0}, \tau> d$. There exist ${\h},q\in\NNN$, depending only on $\tau,d$, such that for any nonlinearity
$\ff\in C^{q}$ satisfying the Hypotheses \ref{hyp2}, 
there exists $\epsilon_0>0$ such that for any $\e$ with $0<\e\g^{-1}<\epsilon_0$, for any set $\Lambda_{o}\subseteq\Lambda$ and for any Lipschitz family ${\bf u}(\lambda)\in{\bf X}^{0}$ defined on  $\Lambda_{o}$  with $||{\bf u}||_{\gots_{0}+{\h},\g}\leq1$ the following holds.
There exist Lipschitz functions $\mu^{\infty}_{h}:\Lambda \to i \RRR$  of the form
\begin{equation}\label{formamu}
\mu^{\infty}_{h}:=\mu^{\infty}_{\s,j}=-\s i m j^{2}+r^{\infty}_{\s,j}, \quad m\in\RRR, \quad h=(\s,j)\in\CC\times\NNN,
\quad \sup_{h}|r^{\infty}_{h}|_{\g}\leq C\e, 
\end{equation}
with $\CC:=\{+1,-1\}$,
such that  $\mu^{\infty}_{\s,j}=-\mu^{\infty}_{-\s,j}$  and setting
\begin{equation}\label{martina10}
\Lambda_{\infty}^{2\g}({\bf u}):=\left\{
\la\in\Lambda_{o} : |\la\bar\oo\cdot\ell\!+\!\mu^{\infty}_{\s,j}(\la)-\!\mu^{\infty}_{\s',j'}(\la)|\geq
\frac{2\g|\s j^{2}-\s' j'^{2}|}{\langle\ell\rangle^{\tau}}, \;\forall\ell\in\ZZZ^{d}, \forall (\s,j),(\s',j')\in\CC\times\NNN
\right\},
\end{equation}
we have:

\noindent
(i) For $\la\in \Lambda_\infty^{2\g}$ there exist linear
bounded operators $W_{1},W_{2} : {\bf X}^{\gots_{0}}\to{\bf X}^{\gots_0}$ with
bounded inverse, such that $\calL({\bf u})$ defined in \eqref{pozzoprofondo} satisfies
\begin{equation}\label{1.2.2}
\begin{aligned}
\calL({\bf u})=W_{1}&\calL_{\infty}W_{2}^{-1}, \quad \calL_{\infty}=\oo\cdot\del_{\f}\id+\DD_{\infty}
\quad
{\rm with}\quad
\DD_{\infty}=diag_{k\in\CC\times\NNN}\{\mu^{\infty}_{h}\}, 
\end{aligned}
\end{equation}
for any $k=(\s,j)\in\CC\times\NNN$. Moreover, for any $s\in(\gots_{0},q-{\h})$, if $||{\bf u}||_{s+{\h},\g}<+\infty$, then
$W_{i}^{\pm1}$ are bounded operators ${\bf X}^{s}\to{\bf X}^{s}$. 

\noindent
(ii) under the same assumption of $(i)$, for any $\f\in\TTT^{d}$ the $W_{i}$ define changes of variables on the phase space
\begin{equation}\label{1.2.4}
W_{i}(\f), W_{i}^{-1}(\f) : {\bf X}^{s}_{x}\to {\bf X}^{s}_{x}, \quad i=1,2,
\end{equation}
see Remark \ref{rmkphase}. Such operators satisfy the bounds
\begin{equation}\label{eq:4.9madonna}
||(W_{i}^{\pm1}(\f)-\id){\bf h}||_{{\bf H}^{s}_{x}}\leq
\e\g^{-1} C(s)(||{\bf h}||_{{\bf H}^{s}_{x}}+||{\bf u}||_{s+\h+\gots_{0}}
||{\bf h}||_{{\bf H}^{1}_{x}}).
\end{equation}
\end{proposition}

\begin{rmk}\label{ciccio}
The purpose of item $(ii)$ is to prove that a function ${\bf h}(t)\in {\bf X}_{x}^{s}$  is a solution of the linearized NLS \eqref{1.2.1}
if and only if the function ${\bf v}(t):=W_{2}^{-1}(\oo t)[{\bf h}(t)]\in {\bf H}^{s}_{x}$ 
solves the constant
coefficients dynamical system
\begin{equation}\label{1.2.6}
\left(\begin{matrix} \del_{t}v \\ \del_{t}\bar{v}\end{matrix}\right) 
+\left(\begin{matrix} \DD_{\infty} &0 \\ 0 &-\DD_{\infty}\end{matrix}\right)
\left(\begin{matrix} v \\ \bar{v}\end{matrix}\right)=\left(\begin{matrix} 0\\ 0 \end{matrix}\right)
\;\Leftrightarrow\quad \dot{v}_{j}+ \mu^{\infty}_{+,j} v_j=0
\quad j\in\NNN,
\end{equation}
  
\noindent
Since the eigenvalues are all
imaginary we have that 
\begin{equation}\label{1.2.7}
||v(t)||^{2}_{H_{x}^{s}}=\sum_{j\in\NNN}|v_{j}(t)|^{2}\langle j\rangle^{2s}=\sum_{j\in\NNN}
|v_{j}(0)|^{2}\langle j\rangle^{2s}=||v(0)||^{2}_{H_{x}^{s}},
\end{equation}
that means that the Sobolev norm in the space of functions depending on $x$, is constant in time.
\end{rmk}

\noindent
Once $\calL({\bf u})$ si diagonal it is trivial to invert it in an explicit Cantor like set. In section \ref{sec:5} we prove
\begin{lemma}\label{inverseofl}({\bf Right inverse of $\calL$})
Under the hypotheses of Proposition \ref{teo2}, set
\begin{equation}\label{eq:4.4.18}
\z:=4\tau+\h+8.
\end{equation}
where $\h$ is fixed in Proposition \ref{teo2}. 
Consider a Lipschitz family ${\bf u}(\la)$ with $\la\in\Lambda_{o}\subseteq\Lambda$ such that
\begin{equation}\label{eq:4.4.19}
||{\bf u}||_{\gots_{0}+\z,\g} \leq1.
\end{equation}
Define the set
\begin{equation}\label{eq:primedimerda}
P_{\infty}^{2\g}({\bf u}):=
\left\{
\la\in\Lambda_{o} : |\la\bar\oo\cdot\ell+\mu^{\infty}_{\s,j}(\la)|\geq
\frac{2\g j^{2}}{\langle\ell\rangle^{\tau}}, \;\;\forall\ell\in\ZZZ^{d}, \;\;\forall\; (\s,j)\in\CC\times\NNN
\right\}.
\end{equation}
There exists $\epsilon_{0}$, depending only on the data of the problem, such that
if $\e\g^{-1}<\epsilon_{0}$
then, for any $\la\in\Lambda^{2\g}_{\infty}({\bf u})\cap P^{2\g}_{\infty}({\bf u})$ 
(see (\ref{martina10})),  and for any Lipschitz family
${\bf g}(\la)\in {\bf Z}^{s}$, 
the equation $\calL {\bf h}:=\calL(\la,{\bf u}(\la)){\bf h}= {\bf g}$,
where $\calL$ is the linearized operator in (\ref{pozzoprofondo}), admits a solution
\begin{equation}\label{eq:4.4.22}
{\bf h}:=\calL^{-1}{\bf g}:=W_{2}\calL_{\infty}^{-1}W_{1}^{-1}{\bf g}\in{\bf X}^{s},
\end{equation}
such that 
\begin{equation}\label{eq:4.4.23}
||{\bf h}||_{s,\g}\leq C(s)\g^{-1}\left(||{\bf g}||_{s+2\tau+5,\g}+
||{\bf u}||_{s+\z,\g}||{\bf g}||_{\gots_{0},\g}
\right), \qquad \gots_{0}\leq s \leq q-\z.
\end{equation}
\end{lemma}
By formula (\ref{eq:4.4.23}) we have good bounds on the inverse of $\calL({\bf u}_{n})$ in the set
$\Lambda^{2\g}_{\infty}({\bf u}_{n})\cap P^{2\g}_{\infty}({\bf u}_{n})$. It is easy to see that this sets have positive measure for all $n\geq0$. Now in the Nash-Moser proposition \ref{teo4} we defined the sets $\calG_n$ in order to ensure bounds on the inverse of  $\calL({\bf u}_{n})$, thus we have the following
\begin{proposition}[{\bf Measure estimates}]\label{measure}
Set $\g_n:=(1+2^{-n})\g$ and consider the set $\calG_{\infty}$ of Proposition \ref{teo4}
with $\mu=\zeta$ defined in Lemma \ref{inverseofl} . We have
\begin{subequations}\label{eq136tot}
\begin{align}
&\cap_{n\geq0}\Lambda^{2\g_n}_{\infty}({\bf u}_n)\cap P^{2\g_n}_{\infty}({\bf u}_n)\subseteq \calG_{\infty}, \label{eq136b}\\
&|\Lambda\backslash\calG_{\infty}|\to 0, \;\; {\rm as} \;\; \g\to0. \label{eq136}
\end{align}
\end{subequations}
\end{proposition}
 Formula (\ref{eq136b}) is essentially trivial. One just need to look at Definition \ref{invertibility}
 and item $(N1)_{n}$ of Theorem \ref{NM}, which fix the sets $\calG_{n}$. The (\ref{eq136}) is more delicate. The first point is that we reduce to computing the measure of the left hand side of (\ref{eq136b}).
 It is simple to show that each $\Lambda^{2\g_n}_{\infty}({\bf u}_n)\cap P^{2\g_n}_{\infty}({\bf u}_n)$ has measure $1-O(\g)$, however in principle as $n$ varies this sets are unrelated and then the intersection might be empty.
 We need to study the dependence of the Cantor sets on the function ${\bf u}_{n}$.
Indeed $\Lambda_{\infty}^{2\g}({\bf u})$ is constructed by imposing infinitely many {\em second Melnikov conditions}. We show that this conditions imply a {\bf finitely} many second Melnikov conditions on a whole neighbourhood of ${\bf u}$. 
\begin{lemma}\label{teo3}
 Under the hypotheses of Proposition \ref{teo2}, for $N$ sufficiently large, for any $0<\rho<\g/2$ and for any Lipschitz family ${\bf v}(\la)\in {\bf X}^{0}$ 
with $\la \in\Lambda_{o}$ such that
\begin{equation}\label{martina}
\sup_{\la\in\Lambda_{o}}||{\bf u}-{\bf v}||_{\gots_{0}+\h}\leq \e C \rho N^{-\tau},
\end{equation}
we have the following.  
For all $\la\in\Lambda_{\infty}^{2\g}({\bf u})$ there exist  invertible
and reversibility-preserving (see  Section \ref{accidenti} for a precise definition) 
transformations $V_{i}$ for $i=1,2$ such that
\begin{equation}\label{martina2}
V_{1}^{-1}\calL({\bf v})V_{2}=\oo\cdot\del_{\f}\id+diag_{h\in\CC\times\NNN}\{\mu^{(N)}_{h}\}+E_{1}\del_{x}+E_{0} : {\bf X}^{0}\to{\bf Z}^{0},
\end{equation}
where $\mu^{(N)}_{h}$ have the same form of $\mu^{\infty}_{k}$ in (\ref{1.2.2}) with bounds
\begin{equation}\label{matrina5}
|r^{\infty}_{h}-r^{(N)}_{h}|_{\g}\leq \e C ||{\bf u}-{\bf v}||_{\gots_{0}+\h,\g}+C \e N^{-\kappa},
\end{equation}
for an appropriate $\kappa$ depending only on $\tau$.
More precisely $\Lambda_{\infty}^{2\g}({\bf u})\subset \Lambda_{N}^{\g-\rho}({\bf v})$  with
$$
\Lambda_{N}^{\g-\rho}({\bf v}):=\left\{
\la\in\Lambda_{o} : |\la\bar\oo\cdot\ell\!+\!\mu^{(N)}_{\s,j}(\la)-\!\mu^{(N)}_{\s',j'}(\la)|\geq
\frac{(\g-\rho)|\s j^{2}-\s' j'^{2}|}{\langle\ell\rangle^{\tau}}, \;\forall \,|\ell|<N,\; \forall (\s,j),(\s',j')\in\CC\times\NNN
\right\}.
$$
Finally the $V_{i}$ satisfy bounds like (\ref{eq:4.9madonna}) and  the remainders satisfy 
\begin{equation}\label{martina3}
||E_{0}{\bf h}||_{s}+||E_{1}{\bf h}||_{s}\leq \e C N^{-\kappa}(||{\bf h}||_{s}+||{\bf v}||_{s+\h}||{\bf h}||_{\gots_{0}}).
\end{equation}
\end{lemma}
\noindent
Since the ${\bf u}_n$ are a rapidly converging Cauchy sequence this proposition allows us to prove that $\calG_{\infty}$ has asymptotically full measure.

\setcounter{equation}{0}
\section{ An Abstract  Existence Theorem }
\label{secNM}
In this Section we prove an Abstract Nash-Moser theorem in Banach spaces.
This abstract formulation essentially shows a method to find solutions
of implicit function problems.
The aim is to apply the scheme to prove Proposition \ref{teo4} to the functional $F$
defined in (\ref{4'}).

\subsection{Nash-Moser scheme}\label{anm}

Let us consider a scale of Banach spaces $(\HH_s, \| \ \|_s)_{s \geq 0}$, such that
$$
\forall  s \leq s', \  \ \HH_{s'} \subseteq \HH_s \ \ {\rm and} \  \ \|u\|_s \leq \|u \|_{s'} \, , \ \forall u \in \HH_{s'} \, , 
$$
and define
$ \HH := \cap_{s\geq 0} \HH_s $.

We assume that there is a non-decreasing family $(E^{(N)})_{N \geq 0}$ of subspaces of $\HH$ such that
$\cup_{N \geq 0} E^{(N)}$ is dense in $\HH_s$ for any $s\geq 0$, and that there are projectors
$$
\Pi^{(N)}: \HH_0 \to E^{(N)}
$$
satisfying:  for any $s\ge0$ and any $\nu\ge0$ there is a positive
constant $ C :=C(s, \nu)$ such that

\begin{itemize}

\item[(P1)] $\|\Pi^{(N)}u\|_{s+\nu}\le CN^{\nu}\|u\|_{s}$
for all $u\in \HH_{s}$, 

\item[(P2)] $\|(\uno - \Pi^{(N)})u\|_{s}\le C
N^{-\nu}\|u\|_{s+\nu}$ for all $u\in \HH_{s+\nu}$.
\end{itemize}
In the following we will work with parameter families of functions in $\HH_s$, more precisely we consider
$
u=u(\lambda)\in {\rm Lip}(\Lambda,\HH_s)
$
where $\Lambda\subset \RRR$.
We define:
\begin{itemize}
\item {\it sup norm:}  $\displaystyle{||f||_{s}^{sup}:=||f||^{sup}_{s,\Lambda}:=\sup_{\la\in\Lambda}||f(\la)||_{s}}$,
\item {\it Lipschitz semi-norm:}  $\displaystyle{||f||_{s}^{lip}:=||f||_{s,\Lambda}^{lip}:=
\sup_{\substack{\la_{1},\la_{2}\in\Lambda \\ \la_{1}\neq\la_{2}}}
\frac{||f(\la_{1})-f(\la_{2})||_{s}}{|\la_{1}-\la_{2}|}}$,
\end{itemize}
and for $\g>0$ the weighted  Lipschitz norm
\begin{equation}\label{lnorm}
||f||_{s,\g}:=||f||_{s,\Lambda,\g}:=||f||^{sup}_{s}+\g||f||_{s}^{lip}.
\end{equation}

Let us consider a $ C^2 $ map $F: [0,\e_0]\times\Lambda \times{\mathcal H}_{\gots_0 + \nu} \to {\mathcal H}_{\gots_0} $ for some $ \nu > 0 $ and assume the following
\begin{itemize}
\item[(F0)] $F$ is of the form
$$ F(\e,\la,u)= L_\la u +\e f(\la,u)$$
where, for all $\la\in \Lambda$, $L_\la$ is a linear operator which preserves all  the subspaces $E^{(N)}$. 

\item[(F1)]  {\em reversibility} property:
\begin{equation}\label{Fsimme}
 \exists \, A_s, B_s \subseteq {  \calH}_s \
\mbox{closed subspaces  of } {  \calH}_s, \, s \geq 0,  \, \mbox{such that} \ 
F: A_{s+\nu}\to  B_s  \, . 
\end{equation} 

\end{itemize}
We assume also  the following tame properties: given $ S' > \gots_0 $, $\forall s \in[\gots_0, S')$, for all Lipschitz map
$u(\la)$ such that
$\|u\|_{\gots_0,\g} \le1$, $(\e,\la)\in[0,\e_{0})\times\Lambda $, 
\begin{itemize}
\item[(F2)] $\|f(\la,u)\|_{s,\g} ,\|L_\la u\|_{s,\g}\le
C(s)(1+\|u\|_{s+\nu,\g})$,

\item[(F3)] $\|d_{u}f(\la,u)[h]\|_{s,\g}\le C(s)\big(\|u\|_{s+\nu,\g}\|h\|_{\gots_0+\nu,\g}+\|h\|_{s+\nu,\g}\big)$,

\item[(F4)] $\|d^{2}_{u}f(\la,u)[h,v]\|_{s,\g} \le C(s)\big(
\|u\|_{s+\nu,\g}\|h\|_{\gots_0+\nu,\g}\|v\|_{\gots_0+\nu,\g}+\|h\|_{s+\nu,\g}\|v\|_{\gots_0+\nu,\g}
+\|h\|_{\gots_0+\nu,\g}\|v\|_{s+\nu,\g}\big)$,

for any two Lipschitz maps $h(\la)$, $v(\la)$.
\end{itemize}
\begin{rmk}\label{porcarever}
Note that $(F1)$  implies $d_uF(\e,\la,v): A_{s+\nu}\to B_s$ for all $v\in A_s$.
\end{rmk}
We denote
\begin{equation}\label{lilli}
\calL(u)\equiv \calL(\la,u):= L_\lambda + \e d_u f(\la,u)\,,
\end{equation}
we have the following definition.
\begin{defi}[{\bf Good Parameters}]\label{invertibility}
Given  $\mu>0,$ $N>1$
 let 
\begin{equation}\label{eq106}
\ka_{1}=6\mu+12\nu, \quad  \ka_{2}=11\mu+25\nu\,,
\end{equation}
for any Lipschitz familty ${ u}(\la)\in E^{(N)}$  with 
$||{ u}||_{\gots_{0}+\mu,\g}\leq1$,
we define the set of good parameters 
 $\la\in\Lambda$  as:
\begin{subequations}\label{eq104}
\begin{align}
\calG_N({u }):= &\left\{\la\in \Lambda \; : \;\; ||\calL^{-1}({u }){h }||_{\gots_{0},\g}\leq C(\gots_{0})\g^{-1}||{h }||_{\gots_{0}+\mu,\g},\right. \label{eq104b}\\
& \;\;\;\quad ||\calL^{-1}({u }){h }||_{s,\g}\leq C(s)\g^{-1}\left(||{h }||_{s+\mu,\g}+
||{u }||_{s+\mu,\g}||{h }||_{\gots_{0},\g}\right), \forall \gots_{0}\leq s\leq \gots_{0}+\ka_{2}-\mu,\label{eq104a}\\
 &\;\;\;\left.\text{ for all Lipschitz maps ${h }(\la)$}\right\}.\nonumber
 \end{align}
\end{subequations} 
\end{defi}

Clearly, Definition \ref{invertibility} depends on  $\mu$ and $ N$.

Given $N_0>1$ we set
$$ N_n= (N_0)^{(\frac32)^n}\,,\quad \HH_n:= E^{(N_n)}\,, \quad A_n:= A_s\cap \HH_n$$
same for the subspace $B$.

In the following, we shall write $a\leq_{s}b$ to denote $a\leq C(s) b$, for some contant $C(s)$ depending on $s$.
In general, we shall write $a\lessdot b$ if there exists a constant $C$, depending only on the data of the problem, such that $a\leq C b$.

 \begin{theorem}\label{NM}{  (Nash-Moser algorithm)}
 Assume $ F$ satisfies $(F0)-(F4)$ and fix $\g_{0}>0$, $\tau>d+1$.
Then, there exist constants $\epsilon_{0}>0, C_{\star}>0,  N_{0}\in \NNN$,  such that 
for  all $\gamma \le \g_0$ and $\e\g^{-1}<\epsilon_{0}$ the following properties hold for any $n\geq0$:

\noindent
$(N1)_{n}$ there exists a function 
\begin{equation}\label{eq107}
{u}_{n} : \calG_{n} \subseteq \Lambda \to {A}_{n},  \qquad 
||{u}_n||_{\gots_{0}+\mu,\g}\leq1,
\end{equation}
where the sets $\calG_{n}$ are defined inductively by $\calG_{0}:=\Lambda$ and
$\calG_{n+1}:= \calG_n\cap \calG_{N_{n}}({u}_{n})$, such that
\begin{equation}\label{eq110}
||F({u}_{n})||_{\gots_{0},\g}\leq C_{\star}\e N_{n}^{-\ka_{1}}.
\end{equation}
Moreover one has that ${h}_{n}:={u}_{n}-{u}_{n-1}$ (with ${h}_{0}=0$)
satisfies
\begin{equation}\label{eq109}
||{h}_{n}||_{\gots_{0}+\mu,\g}\leq C_{\star}\e\g^{-1}N_{n}^{-\kappa_{3}}, \;\;\;
\kappa_{3}:=9\nu+2\mu.
\end{equation} 
The Lipschitz norms are defined on the sets $\calG_{n}$.
\noindent
$(N2)_{n}$ the following estimates in high norms hold:
\begin{equation}\label{eq111}
||{u}_{n}||_{\gots_{0}+\ka_{2},\g}+\g^{-1}||F({u}_{n})||_{\gots_{0}+\ka_{2},\g}
\leq C_{\star}\e\g^{-1} N_{n}^{\ka_{1}}.
\end{equation}

Finally, setting $\calG_{\infty}:=\cap_{n\geq0}\calG_{n}$,
the sequence $({u}_{n})_{n\geq0}$ converges in norm 
$||\cdot||_{\gots_{0}+\mu, \calG_{\infty},\g}$ to a function ${u}_{\infty}$
such that
\begin{equation}\label{eq112}
F(\la,{u}_{\infty}(\la))\equiv0, \quad \;\; 
\sup_{\la\in \calG_{\infty}}||{u}_{\infty}(\la)||_{\gots_{0}+\mu}
\leq C \e\g^{-1}.
\end{equation}
 \end{theorem}

\proof We proceed by induction.

\noindent
 We set ${  u}_{0}={  h}_{0}=0$, we get  $(N1)_{0}$ and $(N2)_{0}$ by fixing
$C_{\star}\geq\max\left\{||f(0)||_{\gots_{0}}N_{0}^{\ka_{1}}, ||f(0)||_{\gots_{0}+\ka_{2}}N_{0}^{-\ka_1}
\right\}$ .

\noindent
We assume inductively $(Ni)_{n}$ for $i=1,2,3$ for some $n\geq0$ and prove
 $(Ni)_{n+1}$ for $i=1,2,3$.
 
\noindent
By $(N1)_{n}$, $u_n\in A_n$ satisfies the conditions in Definition \ref{invertibility}. 
Then, by definition, $\la\in \calG_{n+1}$ implies that  $\calL_{n}:=\calL({  u}_{n})$ is invertible 
with estimates (\ref{eq104}), (used with ${  u}={  u}_{n}$ and $N=N_{n}$).

\noindent
Set
\begin{equation}\label{eq113}
{  u}_{n+1}:={  u}_{n}+h_{n+1}\in {   A}_{n+1}, \quad
{  h}_{n+1}:=-\Pi_{n+1}\calL_{n}^{-1}\Pi_{n+1} F({  u}_{n}),
\end{equation}
which is well-defined. Indeed, $F({  u}_{n})\in { B }_{s}$ implies, since
$\calL_{n}$ maps $A_{s+\nu}\to B_s$ that ${  h}_{n+1}\in {   A}_{n+1}$.
By definition
\begin{equation}\label{eq115}
F({  u}_{n+1})=F({  u}_{n})+\calL_{n}{  h}_{n}+\e \calQ({  u}_{n},{  h}_{n+1}),
\end{equation}
where, by condition $(F0)$ we have
\begin{equation}\label{eq116}
\calQ({  u}_{n},{  h}_{n+1}):=f({  u}_{n}+{  h}_{n+1})-f({  u}_{n})-
d_uf({  u}_{n}){  h}_{n+1}\,,
\end{equation}
which is at least quadratic in ${  h}_{n+1}$.
Then, using the definition of ${  h}_{n+1}$ in (\ref{eq113}) we obtain
\begin{equation}\label{eq117}
\begin{aligned}
F({  u}_{n+1})&=F({  u}_{n})-\calL_{n}\Pi_{n+1}\calL_{n}^{-1}\Pi_{n+1} F({  u}_{n})+
\e\calQ({  u}_{n},{  h}_{n+1})\\
&=\Pi^{\perp}_{n+1}F({  u}_{n})+
\calL_{n}\Pi^{\perp}_{n+1}\calL_{n}^{-1}\Pi_{n+1} F({  u}_{n})+\e\calQ({  u}_{n},{  h}_{n+1})\\
&=\Pi_{n+1}^{\perp}F({  u}_{n})+\Pi^{\perp}_{n+1}\calL_{n}\calL_{n}^{-1}\Pi_{n+1}F({  u}_{n})
+[\calL_{n},\Pi^{\perp}_{n+1}]\calL_{n}^{-1}F({  u}_{n})+\e\calQ({  u}_{n},{  h}_{n+1}),
\end{aligned}
\end{equation}
hence, by using the fact that  by $(F0)$
$[\calL_{n},\Pi^{\perp}_{n+1}]=\e [d_uf(\la,u_n),\Pi^{\perp}_{n+1}]$, one has
\begin{equation}\label{eq118}
F({  u}_{n+1})=\Pi^{\perp}_{n+1}F({  u}_{n})+\e[d_u f({  u}_{n}),\Pi^{\perp}_{n+1}]
\calL_{n}^{-1}\Pi_{n+1}F({  u}_{n})+\e\calQ({  u}_{n},{  h}_{n+1}).
\end{equation}

Now we need a technical Lemma to deduce the estimates (\ref{eq110}) and (\ref{eq111}) at the step
$n+1$. This Lemma guarantees that the scheme is quadratic, and the high norms of the
approximate solutions and of the vector fields do not go to fast to infinity.
\begin{lemma}\label{lemma3.7} 
Set for simplicity
\begin{equation}\label{eq119}
K_{n}:=||{  u}_{n}||_{\gots_{0}+\ka_{2},\g}+\g^{-1}||F({  u}_{n})||_{\gots_{0}+\ka_{2},\g}, \quad
k_{n}:=\g^{-1}||F({  u}_{n})||_{\gots_{0},\g}.
\end{equation}
Then, there exists a constant $C_{0}:=C_{0}( \mu,d,\ka_{2})$ such that
\begin{equation}\label{eq120}
K_{n+1}\leq C_{0}N_{n+1}^{2\mu+4\nu}(1+k_{n})^{2}K_{n}, \qquad
k_{n+1}\leq C_{0}N_{n+1}^{-\ka_{2}+\mu+2\nu}K_{n}(1+k_{n})+C_{0}N_{n+1}^{2\nu+2\mu}k^{2}_{n}
\end{equation}
\end{lemma}

\prova
First of all, we note that, by conditions $(F2)-(F4)$, 
 $\calQ({  u}_{n},\cdot)$ satisfies
\begin{subequations}\label{eq121}
\begin{align}
||\calQ({  u}_{n},{  h})||_{s,\g}&\leq||{  h}||_{\gots_{0}+\nu,\g}\left(||{  h}||_{s+\nu,\g}+
||{  u}_{n}||_{s+\nu,\g}||{  h}||_{\gots_{0}+\nu,\g}
\right),\;\;\; \forall {  h}(\la)\label{eq121a}\\
||\calQ({  u}_{n},{  h})||_{\gots_{0}+\nu,\g}&\leq_{s}N_{n+1}^{2\nu}
||{  h}||_{\gots_{0}+\nu,\g}^{2}, \;\;\; \forall {  h}(\la)\in {  H}_{n+1}\label{eq121b}
\end{align}
\end{subequations}
where ${  h}(\la)\in A_{n+1}$ is a Lipschitz family of functions depending on a parameter.
 The bound (\ref{eq121b}) is nothing but the
(\ref{eq121a}) with $s=\gots_{0}+\nu$,  where we used the fact that $||{  u}_{n}||_{\gots_{0}+\nu}\leq1$
and the smoothing properties $(P1)$, that hold because ${  u}_{n}\in {   A}_{n}$ by definition and $ {  h} \in A_{n+1}$ by hypothesis.

Consider ${  h}_{n+1}$ defined in (\ref{eq113}). then we have
\begin{subequations}\label{eq123}
\begin{align}
||{  h}_{n+1}||_{\gots_{0}+\ka_{2},\g}&\stackrel{(\ref{eq104a})}{\leq}_{\!\!\!\!\gots_{0}+\ka_{2}}
\g^{-1}N_{n+1}^{\mu}\left(||F({  u}_{n})||_{\gots_{0}+\ka_{2},\g}+||{  u}_{n}||_{\gots_{0}+\ka_{2},\g}||F({  u}_{n})||_{\gots_{0},\g}
\right),\label{eq123a}\\
||{  h}_{n+1}||_{\gots_{0},\g}&\stackrel{(\ref{eq104b})}{\leq}_{\!\!\!\!\gots_{0}}\g^{-1}
N_{n+1}^{\mu}||F({  u}_{n})||_{\gots_{0},\g},\label{eq123b}
\end{align}
\end{subequations}
Moreover, recalling that by (\ref{eq113}) one has 
${  u}_{n+1}={  u}_{n}+{  h}_{n+1}$,
we get, by (\ref{eq123}),
\begin{equation}\label{eq124}
||{  u}_{n+1}||_{\gots_{0}+\ka_{2},\g}\leq
||{  u}_{n}||_{\gots_{0}+\ka_{2},\g}
\left(1+\g^{-1}N_{n+1}^{\mu}||F({  u}_{n})||_{\gots_{0},\g}
\right)+
\g^{-1}N_{n+1}^{\mu}||F({  u}_{n})||_{\gots_{0}+\ka_{2},\g}.
\end{equation}

Now, we would like to estimate the norms of $F({  u}_{n+1})$. First of all, we can estimate 
the term
$R_{n}:=[d_u f({  u}_{n}),\Pi^{\perp}_{n+1}]
\calL_{n}^{-1}\Pi_{n+1}F({  u}_{n})$ in (\ref{eq118}), without using the commutator structure, 
\begin{subequations}\label{eq122}
\begin{align}
||R_{n}||_{s,\g}&\leq_{s}\g^{-1}N_{n+1}^{\mu+2\nu}\left(||F({  u}_{n})||_{s,\g}+||{  u}_{n}||_{s,\g}||F({  u}_{n})||_{\gots_{0},\g}
\right),\label{eq122a}\\
||R_{n}||_{\gots_{0},\g}&\leq_{\gots_{0}+\ka_{2},\g}\g^{-1}N_{n+1}^{-\ka_{2}+\mu+2\nu}\left(
||F({  u}_{n})||_{\gots_{0}+\ka_{2},\g}+||{  u}_{n}||_{\gots_{0}+\ka_{2},\g}||F({  u}_{n})||_{\gots_{0},\g}
\right),\label{eq122b}
\end{align}
\end{subequations}
where we used the (\ref{eq104}) to estimate $\calL_{n}^{-1}$, the $(F3)$ for $d_u f$ and 
the smoothing estimates $(P1)-(P2)$.
By (\ref{eq118}), (\ref{eq122b}), (\ref{eq121b}) and using $\e\g^{-1}\leq1$ we obtain, 
\begin{equation}\label{eq125}
\begin{aligned}
||F({  u}_{n+1})||_{\gots_{0},\g}\leq_{\gots_{0}}&
||\Pi^{\perp}_{N_{n+1}}F({  u}_{n})||_{\gots_{0},\g}+
+\e N_{n+1}^{2\n}||{  h}||_{\gots_{0},\g}^{2}\\
&+\e\g^{-1}N_{n+1}^{-\ka_{2}+\mu+2\nu}\left(
||F({  u}_{n})||_{\gots_{0}+\ka_{2},\g}+||{  u}_{n}||_{\gots_{0}+\ka_{2},\g}||F({  u}_{n})||_{\gots_{0},\g}
\right)\\
\stackrel{(P2)}{\leq}_{\!\!\!\gots_{0}+\ka_{2}}&
N_{n+1}^{-\ka_{2}+\mu+2\nu}\left(
||F({  u}_{n})||_{\gots_{0}+\ka_{2},\g}+||{  u}_{n}||_{\gots_{0}+\ka_{2},\g}||F({  u}_{n})||_{\gots_{0},\g}\right)\\
&+
\e\g^{-2}N_{n+1}^{2\nu+2\mu}||F({  u}_{n})||_{\gots_{0},\g}^{2}.
\end{aligned}
\end{equation}
Following the same reasoning as in (\ref{eq125}), by using the estimates
(\ref{eq122a}), (\ref{eq121a}), (\ref{eq123}) and (P2), we get the estimate
in high norm
\begin{equation}\label{eq126}
\!\!||F({  u}_{n+1})||_{\gots_{0}+\ka_{2},\g}\!\leq\!\!
\left(||F({  u}_{n})||_{\gots_{0}+\ka_{2},\g}\!+\!||{  u}_{n}||_{\gots_{0}+\ka_{2},\g}||F({  u}_{n})||_{\gots_{0},\g}
\right)\left(\!1\!+\! N_{n+1}^{\mu+2\nu}\!+\!N_{n+1}^{2\mu+4\mu}\!\g^{-1}\!||F({  u}_{n})||_{\gots_{0},\g}\!\right).
\end{equation}
From the (\ref{eq125}) follows directly the  second of the (\ref{eq120}),
while collecting togheter  (\ref{eq124}) and (\ref{eq126}) one obtain the first of (\ref{eq120}).
\EP 

\noindent
By \eqref{eq110} 
we have that
\begin{equation}\label{eq127}
k_{n}\leq\e\g^{-1}C_{\star}N_{n}^{-\ka_{1}}\leq1,
\end{equation}
if $\e\g^{-1}$ is small enough. Then one has, for $N_{0}$ large enough,
\begin{equation}\label{eq128}
\begin{aligned}
K_{n+1}&\stackrel{(\ref{eq127}), (\ref{eq120})}{\leq}N_{n+1}^{4\nu+2\mu}2K_{n}\leq
C_{\star}\e\g^{-1} N_{n}^{\frac{3}{2}(4\nu+2\mu)}N_{n}^{\ka_{1}}\leq
C_{\star}\e\g^{-1}N_{n+1}^{\ka_{1}}
\end{aligned}
\end{equation}
where we used the fact that, by formula (\ref{eq106}), 
one has $3(2\nu+\mu)+\ka_{1}=\frac{3}{2} \ka_{1}$.
This proves the $(N3)_{n+1}$. In the same way,
\begin{equation}\label{eq129}
\begin{aligned}
k_{n+1}&\stackrel{(N2)_{n}, (\ref{eq120})}{\leq}
2N_{n+1}^{-\ka_{2}+\mu+2\nu}\e\g^{-1}N_{n}^{\ka_{1}}C_{0}+
\e^{2}\g^{-2}C_{0}N_{n}^{\frac{3}{2}(2\nu+2\mu)}N_{n}^{-2\ka_{1}}\leq
\e\g^{-1} C_{\star} N_{n+1}^{-k_{1}},
\end{aligned}
\end{equation}
where we used again the formula \eqref{eq106}.
This proves the $(N2)_{n+1}$.
The bound (\ref{eq109}) follows by $(N2)_{n}$ and by using Lemma \ref{lemma3.7} to
 estimate the norm of ${  h}_{n}$.
 Then we get
\begin{equation}\label{eq131}
||{  u}_{n+1}||_{\gots_{0}+\mu,\g}\leq ||{  u}_{0}||_{\gots_{0}+\mu,\g}+\sum_{k=1}^{n+1}||{  h}_{i}||_{\gots_{0}+\mu,\g}\leq \sum_{k=1}^{\infty}C_{\star}\e\g^{-1}N_{k}^{-\kappa_{3}}\leq 1,
\end{equation}
if $\e\g^{-1}$ is small enough. This means that $(Ni)_{n+1}$, $i=1,2$, hold.

Now, if $\e\g^{-1}$ is small enough, we have by $(N1)_{n}$ that the sequence
$({  u}_{n})_{\geq0}$ is a Cauchy sequence in norm $||\cdot||_{\gots_{0}+\mu,\g}$,
on the set $\calG_{\infty}=\cap_{n\geq0}\calG_{n}$. Hence, we have that
${  u}_{\infty}:=\lim_{n\to\infty}{  u}_{n}$
 solves the equation
since
\begin{equation}\label{eq132}
||F({  u}_{\infty})||_{\gots_{0}+\mu,\gamma}\leq\lim_{n\to\infty}||F({  u}_{n})||_{\gots_{0}+\mu,\gamma}
\leq\lim_{n\to\infty}N_{n}^{\mu}C_{\star}\e N_{n}^{-\ka_{1}}=0.
\end{equation}
This concludes the proof of Theorem \ref{NM}.
\EP

\subsection{Reversible Operators}\label{accidenti}

We now specify to $\HH_s={\bf H}^{s}:= H^s(\TTT^{d+1},\CCC)\times H^s(\TTT^{d+1},\CCC)\cup\calU$, with the notations of Definition \ref{stizzi};
recall that
\begin{equation}\label{1.2.21tris}
||{\bf h}||_{s,\g}:=||{\bf h}||_{H^s\times H^s,\g}=\max\{||h^{+}||_{s,\g}
,||h^{-}||_{s,\g}\}.
\end{equation}

Since we are working one the space of 
functions which are odd in space,  it is more convenient
to use the sine basis in space  instead of the exponential one. Namely  for $u$ odd in space
we have the two equivalent representations:

$$
u(\f,x)=\sum_{\ell\in\ZZZ^{d},j\in\ZZZ}u_{j}(\ell)e^{i(\ell\cdot\f+jx)} =\sum_{\ell\in\TTT^{d},j\in\NNN}\tilde{u}_{j}(\ell)e^{i\ell\cdot\f}\sin jx, $$
setting
$\tilde{u}_{j}(\ell)=2i u_{j}(\ell)$, since $ u_{j}=-u_{-j}
$.
Then we have also two equivalent $H^s$ norms differing by a factor $2$.
In the following we will use the second one which we denote by $\|\cdot \|_s$, because it is more suitable to deal with
odd functions and odd operators.
The same remark holds also for even functions, in that case we will
use the cosine basis of $L_{x}^{2}$.

We will also use this notation.
From a dynamical point of view our solution ${\bf u}(\f,x)\in {\bf H}^{s}(\TTT^{d}\times\TTT)$ can be
seen as a map
\begin{equation}\label{1.2.20}
\TTT^{d}\ni\f \rightarrow h(\f):={\bf u}(\f,x)\in {\bf H}_{x}^{s}:=H_{x}^{s}(\TTT)\times H_{x}^{s}(\TTT)\cap\calU.
\end{equation}
In other words words we look for a curve in the phase space ${\bf H}_{x}^{s}$ that solves (\ref{4'}).
We will denote the norm of $h(\f):=(u(\f,x),\bar{u}(\f,x))$
\begin{equation}\label{1.2.21}
||h(\f)||^{2}_{{\bf H}^{s}_{x}}:=\sum_{j\in\ZZZ}|u_{j}(\f)|^{2}\langle j\rangle^{2s}.
\end{equation}
It can be interpreted as the norm of the function at time a certain time $t$, with $\oo t \leftrightarrow \f$.
The same notation is used also if the function $u$ belongs to some subspaces of even or odd
functions in $H_{x}^{s}$.

Let  $a_{i,j}\in H^{s}(\TTT^{d}\times\TTT)$, on the multiplication operator 
$A=(a_{i,j})_{i,j=\pm1} : {\HH}_{s}\to{\HH}_{s}$, we define the norm
\begin{equation}\label{1.2.21bis}
||A||_{s}:=\max_{i,j=\pm1}\{||a_{i,j}||_{s}\}, \quad ||A||_{s,\g}:=\max_{i,j=\pm1}\{||a_{i,j}||_{s,\g}\}
\end{equation}

Recalling the definitions (\ref{SPACES}), we set,

\begin{defi}\label{reversible}
An operator $R: H^s\to H^s$ is ``{\rm reversible}'' with respect to the reversibility (\ref{b}) if
\begin{equation}\label{rever}
R : X^{s}\to Z^{s}, \quad s\geq0
\end{equation}
We say that $R$ is ``{\rm reversibility-preserving}'' if
\begin{equation}\label{rever2}
R : G^{s}\to G^{s}, \quad {\rm for} \quad G^{s}=X^{s}, Y^{s}, Z^{s}, \quad s\geq0.
\end{equation}
In the same way, we say that $A:{\bf X}^{s}\to {\bf Z}^{s}$, for $s\geq0$ is ``reversible'', while
$A: {\bf G}^{s}\to{\bf G}^{s}$, for ${\bf G}^{s}={\bf X}^{s},{\bf Y}^{s},{\bf Z}^{s}$, $s\geq0$ is ``{reversibility-preserving}''.
\end{defi}

\begin{rmk}
Note that, since  ${\bf X}^{s}= X^s\times X^s \cap\calU$, 
 Definition \ref{reversible} guarantees that
a reversible operator  preserves also the subspace $\calU$, namely $(u,\bar{u})\stackrel{R}{\to}(z,\bar{z})\in{ H}^{s}\times H^s\cap\calU$.
\end{rmk}

\begin{lemma}\label{rmk:rev}
Consider operators $A,B,C$  of the form 
$$
A:=\begin{pmatrix}a_{1}^1(\f,x) & a_{1}^{-1}(\f,x)\\  a_{-1}^1(\f,x) & a_{-1}^{-1}(\f,x) \end{pmatrix}, \quad B:=i \begin{pmatrix}a_{1}^1(\f,x) & a_{1}^{-1}(\f,x)\\  -a_{-1}^1(\f,x) & -a_{-1}^{-1}(\f,x) \end{pmatrix}\,,\quad C:= B\del_x.
$$
One has that $A$ is reversibility-preserving if and only if $a_{\s}^{\s'}\in Y^{s}$ for $\s,\s'=\pm1$.
Moreover $B $ is reversible if and only if $A$ is reversibility-preserving.
Finally $C$ is reversible if and only if $a_{\s}^{\s'}\in X^{s}$.
\end{lemma}

\prova
The Lemma is proved by simply noting that for $u\in X^{s}$
\begin{equation}\label{rmkrev1}
\begin{aligned}
&a_{\s}^{\s'} u\in X^{s}, \;\; i\s a_{\s}^{\s'}\cdot u\in Z^{s}, \quad \forall \;a_{\s}^{\s'}\in Y^{s}, \qquad
ia_{\s}^{\s'}\cdot u_{x}\in Z^{s}, \quad \forall \; a_{\s,\s'}\in X^{s},
\end{aligned}
\end{equation}
using  that $u_{x}\in Y^{s}$ if $u\in X^{s}$.
The fact that the subspace $\calU$ is preserved, follows by the hypothesis
that  $a_{\s}^{\s'}=\ol{a_{\s'}^{\s}}$, that guarantees, for instance $R{\bf u}=(z_{1},z_{2})$ with $z_{1}=\ol{z_{2}}$.
\EP

\subsection{Proof of Proposition \ref{teo4}}
We now prove that our equation \eqref{4'} satisifies the hypotheses of the abstract Nash-Moser theorem.
We fix $\nu=2$ and consider the operator $F : {\bf H}^{s}\to {\bf H}^{s-2}$,
\begin{equation}\label{eq103}
F({\bf u}):=F(\la,{\bf u})=
\left(\begin{matrix}\la\bar{\oo}\cdot\del_{\f}u+i\del_{xx}u \\ \la\bar{\oo}\cdot\del_{\f}\bar{u}-i\del_{xx}\bar{u} \end{matrix}\right)+\e\left(\begin{matrix} i f_{1}(\f,x,u,\bar{u},u_{x},\bar{u}_{x},u_{xx}, \bar{u}_{xx}) \\
- i f_{2}(\f,x,\bar{u},u,\bar{u}_{x},u_{x}, \bar{u}_{xx},u_{xx})
\end{matrix}
\right)
\end{equation}
For simplicity we  write
\begin{equation}\label{eq114}
F({\bf u}):=F(\la,{\bf u})=L_{\oo}{\bf u}+\e f({\bf u})
\end{equation}
where (recall $\oo=\la\bar{\om}$)
\begin{equation}\label{eq114bis}
L_{\la}\equiv L_{\oo}:=
\left(\begin{matrix} \oo\cdot\del_{\f}+i\del_{xx} & 0 \\ 0 & \oo\cdot\del_{\f}-i\del_{xx}\end{matrix}\right),
\qquad f({\bf u}):=\left(\begin{matrix} i f_1(\f,x,u,\bar{u},u_{x},\bar{u}_{x},u_{xx}, \bar{u}_{xx}) \\
- i f_2(\f,x,\bar{u},u,\bar{u}_{x},u_{x}, \bar{u}_{xx},u_{xx})\end{matrix}
\right)
\end{equation}
\noindent
Hypothesis $(F0)$ is trivial.  Hypothesis $(F1)$  holds true with $A_s= {\bf X}^s$, $B_s={\bf Z}^s$ by Hypothesis \ref{hyp2}.

\bigskip

\noindent
Hypotheses $(F2)-(F4)$ follow from the fact that $\mathtt f$ is a $C^q$ composition operator, see Lemmata \ref{lemA2}, \ref{lemA3}. 
Let us discuss in detail the property $(F3)$, which we will use in the next section.

Take ${\bf u}\in{\bf X}^{s}$, then by our extension rules we have
\begin{equation}\label{eq114tris}
\e d_{\bf u} f({\bf u}):= i \sum_{j=0}^2A_j(\f,x,{\bf u}) \del_x^j
\end{equation}
 where, by \eqref{4bis},  the coefficients of the linear operators $A_{j}=A_j(\f,x,{\bf u})$ have the form 
 \begin{equation}\label{eq:3.2}
A_{2}:= \left(\begin{matrix} a_{2} & b_{2} \\ -\bar{b}_{2} & -a_{2}
\end{matrix}\right), \quad A_{1}:=\left(\begin{matrix} a_{1} & b_{1} \\ -\bar{b}_{1} & -\bar{a}_{1}
\end{matrix}\right), \quad A_{0}:=\left(\begin{matrix} a_{0} & b_{0} \\ -\bar{b}_{0} & -\bar{a}_{0}
\end{matrix}\right).
 \end{equation}
with
 \begin{equation}\label{eq:3.222}
 \begin{aligned}
 a_{i}(\f,x)&:=\e (\del_{z^+_{i}}f_{1})
 (\f,x,u,\bar{u},u_{x},\bar{u}_{x},u_{xx},\bar{u}_{xx}),\\
  b_{i}(\f,x)&:=\e (\del_{z^-_{i}}f_{1})
 (\f,x,u,\bar{u},u_{x},\bar{u}_{x},u_{xx},\bar{u}_{xx}).
 \end{aligned}
 \end{equation}

Thanks to Hypothesis \ref{hyp2}, and Remark \ref{porcarever} 
one has that $d_{{\bf u}}f({\bf u}) : {\bf X}^{0}\to {\bf Z}^{0} $ and hence
\begin{equation}\label{eq:3.223}  
a_{i}, b_{i}\in Y^{s},\;\; i=0,2, \qquad a_{1},b_{1}\in X^{s}.
\end{equation}
By (\ref{eq:3.2}) and Lemma \ref{rmk:rev}, the (\ref{eq:3.223}) implies
\begin{equation}\label{eq:3.3}
iA_{2}, iA_{0} : {\bf X}^{0} \to {\bf Z}^{0}, \quad iA_{1}\del_{x} : {\bf X}^{0} \to {\bf Z}^{0}.
\end{equation}
then the operator $\calL= d_{\bf u} F$ maps ${\bf X}^{0}$ to ${\bf Z}^{0}$, i.e. it is \emph{reversible} according to Definition \ref{reversible}.

The coefficients $a_{i}$ and $b_{i}$ and their derivative $d_{u^{\s}}a_{i}({\bf u})[h]$
with respect to $u^{\s}$ 
in the direction $h$, for $h\in H^{s}$, satisfy the following tame estimates.

\begin{lemma}\label{lem:3.1}
 For all $\gots_{0}\leq s\leq q-2$, 
$||u||_{\gots_{0}+2}\leq1$ we have, for any $i=0,1,2$, $\s=\pm1$
\begin{subequations}
\begin{align}
||b_{i}({\bf u})||_{s}, ||a_{i}({\bf u})||_{s}&\leq \e C(s)(1+||{\bf u}||_{s+2}),\label{eq:3.4a}\\
||d_{u^{\s}}b_{i}({\bf u})[h]||,||d_{u^{\s}}a_{i}({\bf u})[h]||_{s}&\leq\e C(s)(|| h||_{s+2}+
||{\bf u}||_{s+2}|| h||_{\gots_{0}+2})\label{eq:3.4b}. 
\end{align}
\end{subequations}
If moreover $\la\to {\bf u}(\la)\in{\bf H}^{s}$ is a Lipschitz family such that 
$||{\bf u}||_{s,\g}\leq 1$, then
\begin{equation}\label{eq:3.5}
||b_{i}({\bf u})||_{s,\g},||a_{i}({\bf u})||_{s,\g}\leq \e C(s)(1+||{\bf u}||_{s+2,\g}).
\end{equation}
\end{lemma}

\prova To prove the (\ref{eq:3.4a}) it is enough to apply Lemma \ref{lemA2}(i) to the 
function $\del_{z_{i}^{\s}}f_{1}$, for any $i=0,1,2$ and $\s=\pm1$ which holds for $s+1\leq q$.
Now, let us write, for any  $i=0,1,2$  and $\s,\s'=\pm$,
\begin{equation}\label{eq:3.6}
\begin{aligned}
d_{u^{\s}}a_{i}({\bf u})[ h]&\stackrel{(\ref{eq:3.2})}{=}\e
\sum_{k=0}^{2}(\del^{2}_{z_{k}^{\s}z^{+}_{i}}f_{1})(\f,x,u,\bar{u},u_{x},\bar{u}_{x},
u_{xx},\bar{u}_{xx})\del^{k}_{x}h, \\
d_{u^{\s}}b_{i}({\bf u})[ h]&\stackrel{(\ref{eq:3.2})}{=}\e
\sum_{k=0}^{2}(\del^{2}_{z_{k}^{\s}z^{-}_{i}}f_{1})(\f,x,u,\bar{u},u_{x},\bar{u}_{x},
u_{xx},\bar{u}_{xx})\del^{k}_{x}h, 
\end{aligned}
\end{equation}

Then, by Lemma \ref{lemA2}(i) applied on $\del^{2}_{z_{k}^{\s'}z_{i}^{\s}}f_{1}$ we obtain
\begin{equation}\label{eq:3.7}
||(\del^{2}_{z_{k}^{\s'}z_{i}^{\s}}f_{1})(\f,x,u,\bar{u},u_{x},\bar{u}_{x},u_{xx},\bar{u}_{xx})||_{s}
\leq C(s)||f||_{C^{s+2}}(1+||u^{\s}||_{s+2}),
\end{equation}
for $s+2\leq q$. The bound (\ref{eq:3.4b}) follows by (\ref{A5}) using the (\ref{eq:3.7}).
To prove the (\ref{eq:3.5}) one can reason similarly.
\EP

This Lemma ensures property $(F3)$. Properties $(F2)$ and $(F4)$ are proved in exactly in the same way, 
for property $(F4)$ just consider derivatives of $f$ of order $3$.

\smallskip
We have verified all the Hypotheses of Theorem \ref{NM}, which ensures the existence of a solution defined on some {\em possibly empty} set of parameters $\calG_\infty$.
This concludes the proof of Proposition \ref{teo4}.

\setcounter{equation}{0}
\section{The diagonalization algorithm: regularization}
\label{sec:3}

For ${\bf u}\in{\bf X}^{0}$ we consider the linearized operator
\begin{equation}\label{merdaccia}
\calL({\bf u}):=L_{\oo}+\e d_{\bf u} f({\bf u})=\omega\cdot \partial_\f \id +i (E+A_2(\f,x,{\bf u}) \del_{xx} + i A_1\del_x +i A_0\,,\quad E=\begin{pmatrix}1&0\\0&-1\end{pmatrix}
\end{equation}
with $d_{\bf u} f({\bf u})$ defined in formula \eqref{eq114tris} and $\|{\bf u}\|_{\gots_0+2}$ small. 
In this Section we prove 
\begin{lemma}\label{lem:3.8}
Let $\ff\in C^{q}$ satisfy  the Hypotheses of Proposition \ref{teo4} and assume $q>\h_{1}+\gots_{0}$  
where 
\begin{equation}\label{eq:3.2.0}
\h_{1}:=d+2\gots_{0}+10.
\end{equation}
There exists $\epsilon_0>0$ such that, if $\e\g_{0}^{-1}\leq \epsilon_0$ (see (\ref{dio}) for the definition of $\g_0$) then, for any $\g\le \g_0$ and for all 
 ${\bf u} \in{\bf X}^{0}$ depending in a Lipschitz way on  $\la\in \Lambda$,
 if
\begin{equation}\label{eq:3.2.1b}
||{\bf u}||_{\gots_{0}+\h_{1},\g}\leq1,
\end{equation}
 then, for $\gots_{0}\leq s\leq q-\h_{1}$, 
the following holds.

\noindent (i) 
 There exist invertible maps $\VV_1, \VV_2: {\bf H}^{0}\to{\bf H}^{0}$ such that
 $\calL_{4}:=\VV_{1}^{-1}\calL\VV_{2}$
 with
 \begin{equation}\label{eq:3.5.9}
 \calL_{4}:=
 \oo\cdot\del_{\f}\id+i\left(\begin{matrix}m\!\!\! & 0 \\ 0 & \!\!\! -m\end{matrix}\right)\del_{xx}+
 i\left(\begin{matrix}0\!\!\! & q_{1}(\f,x) \\ -\bar{q}_{1}(\f,x) & \!\!\! 0\end{matrix}\right)\del_{x}
 +i\left(\begin{matrix}q_{2}(\f,x) \!\!\! & q_{3}(\f,x) \\ -\bar{q}_{3}(\f,x) &\!\!\! -\bar{q}_{2}(\f,x) \end{matrix}\right).
 \end{equation}
The  $\VV_{i}$ are reversibility-preserving and moreover for all ${\bf h}\in {\bf X}^{0}$
 \begin{equation}\label{eq:3.2.3}
 ||\VV_{i}{\bf h}||_{s,\g}+||\VV_{i}^{-1}{\bf h}||_{s,\g}\leq 
 C(s)(|| {\bf h}||_{s+2,\g}+
 ||{\bf u}||_{s+\h_{1},\g}||{\bf h}||_{\gots_{0}+2,\g}), \quad i=1,2.
\end{equation}
\noindent
(ii) The coefficient $m:=m({\bf u})$ of $\calL_{4}$ satisfies
\begin{subequations}\label{eq:3.2.4}
\begin{align}
|m({\bf u})-1|_\g&\leq \e C,\label{eq:3.2.4a}\\
|d_{{\bf u} }m({\bf u})[{\bf h}]|&\leq \e C||{\bf h}||_{\h_{1}}.\label{eq:3.2.4b}
\end{align}
\end{subequations}
\noindent
(iii) The operators  $q_{i}:=q_{i}({\bf u})$, are such that
\begin{subequations}\label{eq:3.2.7}
\begin{align}
||q_{i}||_{s,\g}&\leq \e C(s)(1+||{\bf u}||_{s+\h_{1},\g}),\label{eq:3.2.6a}\\
||d_{{\bf u} }(q_{i})({\bf u})[{\bf h}]||_{s}&\leq \e C(s)(||{\bf h}||_{s+\h_{1}}+||{\bf u}||_{s+\h_{1}}+
||{\bf h}||_{\gots_{0}+\h_{1}}),\label{eq:3.2.6b}
\end{align}
\end{subequations}
Finally 
$\calL_{4}$ is reversible.
\end{lemma}

The rest of the Section is devoted to the proof of this Lemma. We divide it in four steps. 
at each step we construct a \emph{reversibility-preserving} change of variable $\TT_{i}$ that conjugates\footnote{Actually in the third step we only are able to conjugate $\calL_{2}$ to $\rho\calL_{3}$, where $\rho$ is a suitable function. This is the reason why $\calL$ is semi-conjugated to $\calL_{4}$. } $\calL_{i}$
to $\calL_{i+1}$ where $\calL_{0}:=\calL$ and 
\begin{equation}\label{elle}
\begin{aligned}
\calL_{i}:=\oo\cdot\del_{\f}\id+i\left(\begin{matrix}1+a_{2}^{(i)}(\f,x) \!\!\!\!&\!\!\! b_{2}^{(i)}(\f,x) \\ \!\!\! -\bar{b}_{2}^{(i)}(\f,x) &\!\!\!\! -1-a_{2}^{(i)}(\f,x)\end{matrix}\right)\del_{xx}\!+\! i
\left(\begin{matrix}a^{(i)}_{1}(\f,x)\!\! &\!\! b^{(i)}_{1}(\f,x)\\ \!\! -\bar{b}^{(i)}_{1}(\f,x) &\!\! -\bar{a}^{(i)}_{i}(\f,x)\end{matrix}\right)\del_{x}\!
+\!
\left(\begin{matrix}a^{(i)}_{0}(\f,x) \!\! & \!\! b^{(i)}_{0}(\f,x) \\ -\bar{b}^{(i)}_{0}(\f,x) \!\! &\!\! -\bar{a}^{(i)}_{0}(\f,x)\end{matrix}\right).
\end{aligned}
\end{equation}
On the transformation we need to prove bounds like
\begin{subequations}\label{eq:B15}
\begin{align}
||\TT_{i}({\bf u}){\bf h}||_{s,\g}+||\TT^{-1}_{i}({\bf u}){\bf h}||_{s,\g}&\leq C(s) (||{\bf h}||_{s,\g}+||{\bf u}||_{s+\ka_{i},\g}||{\bf h}||_{\gots_{0}}), \label{eq:B15a}\\
||d_{{\bf u}}\TT_{i}({\bf u})[{\bf h}]{\bf g}||_{s}, ||d_{{\bf u}}\TT^{-1}_{i}({\bf u})[{\bf h}]{\bf g}||_{s}&\leq \e C(s)\left(||{\bf g}||_{s+1}||{\bf h}||_{\gots_{0}+\ka_{i}}+
||{\bf g}||_{2}||{\bf h}||_{s+\ka_{i}}+||{\bf u}||_{s+\ka_{i}}||{\bf g}||_{2}||{\bf h}||_{\gots_{0}}
\right),\label{eq:B15b}
\end{align}
\end{subequations}
for suitable $\ka_{i}$.
Moreover the coefficients in (\ref{elle}) satisfy
\begin{subequations}\label{eq:3.10}
\begin{align}
||a_{j}^{(i)}({\bf u})||_{s,\g}, ||b_{j}^{(i)}({\bf u})||_{s,\g}&\leq \e C(s)(1+||{\bf u}||_{s+\ka_{i},\g}),\label{eq:3.10a}\\
||d_{u^{\s}}a_{j}^{(i)}({\bf u})[h]||_{s}, ||d_{u^{\s}}b_{j}^{(i)}({\bf u})[h]||_{s}&\leq \e C(s)(||{\bf h}||_{s+\ka_{i}}+||{\bf u}||_{s+\ka_{i}}+
||{\bf h}||_{\gots_{0}+\ka_{i}}),\label{eq:3.10b}
\end{align}
\end{subequations}
for $j=0,1,2$ and $i=1,\ldots,4$.

\bigskip
\noindent
{\bf Step 1. Diagonalization of the second order coefficient}
\bigskip

We first diagonalize the term $E +A_{2}$ in \eqref{merdaccia}.
By a direct calculation, one can see that the matrix $(E +A_{2})$ has eigenvalues
$\la_{1,2}=\pm\sqrt{(1+a_{2})^{2}-|b_{2}|^2}$. Hence we set
$a_{2}^{(1)}(\f,x)=\la_{1}-1$. We have that $a_{2}^{(1)}\in\RRR$ because $a_{2}\in\RRR$ and $a_i,b_i$ are small. The diagonalizing matrix is
\begin{equation}\label{fiorellino}
\TT^{-1}_{1}:=\frac{1}{2}\left(\begin{matrix}2+a_{2}+a_{2}^{(1)} & b_{2} \\
-\bar{b}_{2}& -(2+a_{2}+a_{2}^{(1)})\end{matrix}
\right)\;\Rightarrow\; \TT^{-1}_{1}(E +A_{2})\TT_{1}=\left(\begin{matrix}1+a^{(1)}_{2}(\f,x) & 0\\ 0 & -1-a_{2}^{(1)}(\f,x)\end{matrix}\right)
\end{equation}

The tame estimates (\ref{eq:3.10}) for $a_{2}^{(1)}$ and the (\ref{eq:B15}) on $\TT^{-1}_{1}$ follow with $\ka_{1}=2$ by 
(\ref{eq:3.4a}), (\ref{eq:3.2.1b}) and  (\ref{A5}). 
The bound on $\TT_{1}$ 
 follows since
$
{\rm det} \TT^{-1}_{1}=(|b_{2}|^{2}-(2+a_{2}+a_{2}^{(1)})^{2})/4, 
$
and by using the same strategy  as for $a_{2}^{(1)}$.

One has
\begin{equation}\label{eq:3.12}
\begin{aligned}
\calL_{1}&:=\TT^{-1}_{1}\calL  \TT_{1}=\oo\cdot\del_{\f}\id+i  \TT^{-1}_{1}(E+A_{2}) \TT_{1} \del_{xx}+i\left[2 \TT^{-1}_{1}(E+A_{2})\del_{x} \TT_{1}
+ \TT^{-1}_{1}A_{1} \TT_{1}\right]\del_{x}\\
&+i\left[-i \TT^{-1}_{1}(\oo\cdot\del_{\f} \TT_{1})+ \TT^{-1}_{1}(E+A_{2})\del_{xx} \TT_{1}+ \TT^{-1}_{1}A_{1}\del_{x}\TT_{1}+\TT^{-1}_{1}A_{0}\TT_{1}\right];
\end{aligned}
\end{equation}
the (\ref{eq:3.12}) has the form (\ref{elle}) and this identifies the coefficients $a_{j}^{(1)},b_{j}^{(1)}$.
Note that  the matrix of the second order operator is now diagonal.
Moreover, by (\ref{A5}), (\ref{eq:3.10}) on $a_{2}^{(1)}$, (\ref{eq:3.4a}) and (\ref{eq:3.4b})
one obtains the bounds (\ref{eq:3.10}) for the remaining coefficients $a_{j}^{(1)},b_{j}^{(1)}$ with $\ka_{1}:=5$. Then we can fix $\ka_{1}=5$ in all the bounds (\ref{eq:B15a})-(\ref{eq:3.10b}) even if for some of the coefficients there are better bounds. 

Finally, since the matrix $\TT^{-1}_{1}$ is $E+A_{2}$ plus a diagonal matrix
with even components, it has the same parity properties of $A_{2}$, then maps
${\bf Y}^{s}$ to ${\bf Y}^{s}$ and ${\bf X}^{s}$ to 
${\bf X}^{s}$, this means that it is reversibility-preserving and hence $\calL_{1}$ is reversible.
In particular one has that $a_{2}^{(1)}, a^{(1)}_{0},b^{(1)}_{0}\in Y^{0}$ and $a^{(1)}_{1},b^{(1)}_{1}\in X^{0}$
then by Lemma \ref{rmk:rev}.

\begin{rmk}\label{rmk:3.4}
We can note that in the \emph{quasi-linear} case this first step can be avoided. Indeed
in that case one has
$
\del_{\bar{z}_{2}}f\equiv0, 
$
so that the matrix $A_{2}$ is already diagonal, with real coefficients.
\end{rmk}


\bigskip
\noindent
{\bf Step 2. Change of the space variable}
\bigskip

We consider a $\f-$dependent family of diffeomorphisms of the $1-$dimensional torus $\TTT$
of the form
\begin{equation}\label{20}
y=x+\x(\f,x),
\end{equation}
where $\x$ is as small real-valued funtion, $2\pi$ periodic in all its arguments. The change of variables
$(\ref{20})$ induces on the space of functions the invertible linear operator
\begin{equation}\label{21}
(\TT_{2} h)(\f,x):=h(\f,x+\x(\f,x)), \; {\rm with\;\; inverse} \;\;\;\;
(\TT_{2}^{-1}v)(\f,y)=v(\f,y+\widehat{\x}(\f,y))\,,
\end{equation}
where $y\to y+\widehat{\x}(\f,y)$ is the inverse  diffeomorphism of $(\ref{20})$.
With a slight abuse of notation we  extend the operator to ${\bf H}^s$:
\begin{equation}\label{eq:3.14}
{{\TT_{2}}} : {\bf H}^{s}\to {\bf H}^{s}, \quad {\TT_{2}}{\bf h}=\left(
\begin{matrix} (\TT_{2} h)(\f,x) \\ (\TT_{2} \bar{h})(\f,x)\end{matrix}\right).
\end{equation}

Now we have to calculate the conjugate ${\TT_{2}}^{-1}\calL_{1}{\TT_{2}}$ 
of the operator $\calL_{1}$ in 
$(\ref{eq:3.12})$.

\noindent
The conjugate $\TT_{1}^{-1} a\TT_{2}$ of any multiplication operator
$a : h(\f,x) \to a(\f,x)h(\f,x)$ is the multiplication operator $(\TT_{2}^{-1}a) : v(\f,y)\to(\TT_{2}^{-1}a)(\f,y)v(\f,y)$.
The conjugate of the differential operators will be
\begin{equation}\label{22}
\begin{aligned}
\TT_{2}^{-1}\oo\cdot\del_{\f}\TT_{2}&=\oo\cdot\del_{\f}+[\TT_{2}^{-1}(\oo\cdot\del_{\f}\x)]\del_{y},\qquad
\calA^{-1}\del_{x}\TT_{2}&=[\TT_{2}^{-1}(1+\x_{x})]\del_{y},\\
\TT_{2}^{-1}\del_{xx}\calA&=[\TT_{2}^{-1}(1+\x_{x})^{2}]\del_{yy}+[\TT_{2}^{-1}(\x_{xx})]\del_{y},
\end{aligned}
\end{equation}
where all the coefficients
are periodic functions of $(\f,x)$.
Thus
we have obtained $\calL_{2}={\TT_{2}}^{-1}\calL_{1}{\TT_{2}} $ where $\calL_{2}$ has the form
(\ref{elle}).
Note that  the second rows are the complex conjugates of the first, this is due to the fact that $\TT_{2}$ trivially preserves the subspace $\calU$.
We have
\begin{equation}\label{24}
\begin{aligned}
&1+a_{2}^{(2)}(\f,y)=\TT_{2}^{-1}[(1+a_{2}^{(1)})(1+\x_{x})^{2}], \;\; \quad b^{(2)}_{1}(\f,y)=\TT_{2}^{-1}[b^{(1)}_{1}(1+\x_{x})],\\
& a^{(2)}_{1}(\f,y)=
\TT_{2}^{-1}((1+a_{2}^{(1)})\x_{xx})-i\TT_{2}^{-1}(\oo\cdot\del_{\f}\x)+\TT_{2}^{-1}[a^{(1)}_{1}(1+\x_{x})],\\
& a^{(2)}_{0}(\f,y)=\TT_{2}^{-1}[a^{(1)}_{0}], \;\;\; \quad
 b^{(2)}_{0}(\f,y)=\TT_{2}^{-1}[b^{(1)}_{0}].
\end{aligned}
\end{equation}

We are looking for $\x(\f,x)$ such that the coefficient of the second order differential operator does not depend on
$y$, namely
\begin{equation}\label{25}
\TT_{2}^{-1}[(1+a_{2}^{(1)})(1+\x_{x})^{2}]=1+a_{2}^{(2)}(\f),
\end{equation}
for some function $a_{2}^{(2)}(\f)$. Since $\TT_{2}$ operates only on the space variables, the $(\ref{25})$
is equivalent to 
\begin{equation}\label{26}
(1+a_{2}^{(1)}(\f,x))(1+\x_{x}(\f,x))^{2}=1+a_{2}^{(2)}(\f).
\end{equation}
Hence 
we have to set
\begin{equation}\label{27}
\x_{x}(\f,x)=\rho_{0}, \qquad \rho_{0}(\f,x):=(1+a_{2}^{(2)})^{\frac{1}{2}}(\f)(1+a_{2}^{(1)}(\f,x))^{-\frac{1}{2}}-1,
\end{equation}
that has solution $\x$ periodic in $x$ if and only if $\int_{\TTT}\rho_{0} dy=0$. This condition
implies
\begin{equation}\label{28}
a_{2}^{(2)}(\f)=\left(\frac{1}{2\pi}\int_{\TTT}(1+a_{2}^{(1)}(\f,x))^{-\frac{1}{2}}\right)^{-2}-1
\end{equation}
Then we have the solution (with zero average) of $(\ref{27})$
\begin{equation}\label{29bis}
\x(\f,x):=(\del_{x}^{-1}\rho_{0})(\f,x),
\end{equation}
where $\del_{x}^{-1}$ is defined by linearity as
\begin{equation}\label{30}
\del_{x}^{-1}e^{ikx}:=\frac{e^{ikx}}{ik}, \quad \forall \; k\in\ZZZ\backslash\{0\}, \quad \del_{x}^{-1}=0.
\end{equation}
In other word $\del_{x}^{-1}h$ is the primitive of $h$ with zero average in $x$.
  Thus, conjugating $\calL_{1}$ through the operator ${\TT_{2}}$ in (\ref{eq:3.14}) , we obtain the  operator 
 $\calL_{2}$ in (\ref{elle}).
 
Now we start by proving that the coefficient $a_{2}^{(2)}$ satisfies tame estimates like (\ref{eq:3.10})
with $\ka_{2}=2$.
Let us write
\begin{equation}\label{eq:3.46}
a_{2}^{(2)}(\f)=\psi\left(G[g(a_{2}^{(1)})-g(0)]\right)-\psi(0), \quad \psi(t):=(1+t)^{-2}, \quad 
Gh:=\frac{1}{2\pi}\int_{\TTT}h d x, \quad g(t):=(1+t)^{-\frac{1}{2}}.
\end{equation}
Then  one has, for $\e$ small,
\begin{equation}\label{eq:3.47}
||a_{2}^{(2)}||_{s}\stackrel{(\ref{A13})}{\leq}C(s)||G[g(a_{2}^{(1)})-g(0)]||_{s}\leq C(s)||g(a_{2}^{(1)})-g(0)||_{s}
\stackrel{(\ref{A13})}{\leq} C(s)||a_{2}^{(1)}||_{s}.
\end{equation}
In the first case we used (\ref{A13}) on the function $\psi$ 
with $u=0, p=0, h=G[g(a_{2}^{(1)})-g(0)]$, while in the second case
we have set $u=0,p=0, h=a_{2}^{(1)}$ and used the estimate on $g$. Then we used 
the (\ref{eq:3.10}) and the bound (\ref{eq:3.4a}), with $s_{0}=\gots_{0}$ which holds for $s+2\leq q$.
 By (\ref{eq:3.46}), we get for $\s=\pm1$
\begin{equation}\label{eq:3.48}
d_{u^{\s}}a_{2}^{(2)}({\bf u})[h]=\psi'\left(G[g(a_{2}^{(1)})-g(0)]\right)
G\left[g'(a_{2}^{(1)})d_{u^{\s}}a_{2}^{(1)}[h]\right]
\end{equation}
Using (\ref{A5}) with $s_{0}=\gots_{0}$,
Lemma \ref{lemA2}(i) to estimate the functions $\psi'$ and $g'$, as done in (\ref{eq:3.47})),
and by the (\ref{eq:3.4b}) we get (\ref{eq:3.10b}).
The (\ref{eq:3.10a}) follows by (\ref{eq:3.47}), (\ref{eq:3.10b}) and Lemma \ref{lemA3}.
\noindent
The second step is to give tame estimates on the function $\x=\del_{x}^{-1}\rho_{0}$ defined in 
(\ref{27}) and (\ref{29bis}). It is easy to check that, estimates 
(\ref{eq:3.10})  are satisfied also by $\rho_{0}$. They follow by using
the estimates on $a_{2}^{(2)}$ and the estimates (\ref{eq:3.10}), (\ref{eq:3.4a}), (\ref{eq:3.4b}), (\ref{eq:3.5}) 
for $a_{2}^{(1)}$.
By defining $|u|_{s}^{\infty}:=||u||_{W^{s,\infty}}$ and using Lemma \ref{A}(i) we get
\begin{subequations}
\begin{align}
|\x|^{\infty}_{s}&\leq C(s)||\x||_{s+\gots_{0}}\leq C(s)||\rho_{0}||_{s+\gots_{0}}\leq
\e C(s)(1+||{\bf u}||_{s+\gots_{0}+2}),\label{eq:3.50a}\\
|d_{u^{\s}}\x({\bf u})[h]|^{\infty}_{s }&\leq
\e C(s)(||h||_{s+\gots_{0}+2}+||{\bf u}||_{s+\gots_{0}+2}||h||_{\gots_{0}+2}),\label{eq:3.50b}
\end{align}
\end{subequations}
and hence, by Lemma \ref{lemA3} one has
\begin{equation}\label{eq:3.51}
|\x|^{\infty}_{s,\g}\leq \e C(s)(1+||{\bf u}||_{s+\gots_{0}+2,\g}),
\end{equation}
for any $s+\gots_{0}+2\leq q$. 
The diffomorphism $x\mapsto x+\x(\f,x)$ is well-defined if $|\x|_{1,\infty}\leq1/2$, but it is easy to note that
this condition is implied requiring $\e C(s)(1+||{\bf u}||_{\gots_{0}+3})\leq1/2$.
\noindent
Let us study the inverse diffeomorphism $(\f,y)\mapsto(\f, y+\widehat{\x}(\f,y))$ of $(\f,x)\mapsto(\f, x+\g(\f,x))$.
Using Lemma \ref{change}(i) on the torus $\TTT^{d+1}$, one has
\begin{equation}\label{eq:3.52}
|\widehat{\x}|^{\infty}_{s }\leq C|\x|^{\infty}_{s }\leq \e C(s)(1+||{\bf u}||_{s+\gots_{0}+2}).
\end{equation}
By definition we have that $\widehat{\x}(\f,y)+\x(\f,y+\widehat{\x}(\f,y))=0$, which implies, for $\s=\pm1$,
\begin{equation}\label{eq:B10}
\begin{aligned}
|d_{u^{\s}}\widehat{\x}({\bf u})[h]|^{\infty}_{s }
&\leq \e C(||h||_{\gots_{0}+2}+||{\bf u}||_{s+\gots_{0}+3}||h||_{\gots_{0}+2}).
\end{aligned}
\end{equation}
Now, thanks to bounds (\ref{eq:3.52}) and (\ref{eq:B10}), using again Lemma \ref{lemA3} with $p=\gots_{0}+3$, we obtain
\begin{equation}\label{eq:B11}
|\widehat{\x}|^{\infty}_{s,\g} \leq \e C(s)(1+||{\bf u}||_{s+\gots_{0}+3,\g}).
\end{equation}

We have to estimate $\TT_{2}({\bf u})$ and $\TT_{2}^{-1}({\bf u})$. 
By using (\ref{A20c}), (\ref{eq:3.51}) and (\ref{eq:B11}), we get 
the (\ref{eq:B15a}) with $\ka_{2}=\gots_{0}+3$,
Now, since $$d_{{\bf u}}(\TT_{2}({\bf u})g)[{\bf h}]:=d_{{\bf u}}g(\f,x+\x(\f,x;{\bf u}))=(\TT_{2}({\bf u}) g_{x})d_{{\bf u}}\x({\bf u})[{\bf h}],$$ 
we get the (\ref{eq:B15b}) using the (\ref{A7}), (\ref{eq:3.50b}) and (\ref{eq:B15a}). 
The (\ref{eq:B15b}) on $\TT^{-1}_{2}$ follows by the same reasoning.
Finally, using the bounds 
(\ref{A7}), (\ref{eq:B15}), 
(\ref{eq:B11}), (\ref{eq:3.5}), Lemma \ref{lem:3.1} and 
$||{\bf u}||_{\gots_{0}+\h_{1},\g} \leq1$, one has the (\ref{eq:3.10a}) on the coefficients
$a_{j}^{(2)},b_{j}^{(2)}$ for $j=0,1$ in (\ref{24}). 
Now, by definition (\ref{24}), we can write
\begin{equation}\label{eq:B20}
a_{1}^{(2)}=\TT_{2}^{-1}({\bf u})\rho_{1}, \quad
\rho_{1}:=(1+a_{2}^{(1)})\x_{xx}-i\oo\cdot\del_{\f}\x+a_{1}^{(1)}(1+\x_{x}),
\end{equation}
so that, thanks to bounds in Lemma \ref{lem:3.1}, and
(\ref{eq:3.50a}), (\ref{eq:3.50b}), (\ref{A7}) and recalling that
$||{\bf u}||_{\gots_{0}+\h_{1}}\leq1$, we get the (\ref{eq:3.10a}) on $\rho_{1}$.
Now, the (\ref{eq:3.10b}) on $a_{1}^{(2)}$ follows by using the chain rule, setting
$\ka_{2}=\gots_{0}+5$
and for $s+\gots_{0}+5\leq q$. 
The same bounds on the coefficients 
$a_{0}^{(2)}, b_{0}^{(2)}$ are obtained in the same way.

\begin{rmk}\label{rmk2}
Note that $\x$ is a real function and  $\x(\f,x)\in X^{0}$ since $a\in Y^{0}$. This implies that
the operators $\TT_2$  and $\TT_2^{-1}$ map $X^{0}\to X^{0}$ and 
$Y^{0}\to Y^{0}$, namely preserves the parity 
properties of the functions. Moreover we have that $a_{2}^{(2)}, a^{(2)}_{0}, b^{(2)}_{0}\in Y^{0}$,
 while
$a^{(2)}_{1},b^{(2)}_{1}\in X^{0}$. Then 
then by Lemma \ref{rmk:rev},   one has that
the operator $\calL_{2}$ is reversible. 
\end{rmk}

\bigskip
\noindent
{\bf Step 3. Time reparametrization}
\bigskip

In this section we want to make constant the coefficient of the highest order spatial derivative
operator $\del_{yy}$ of $\calL_{2}$, by a quasi-periodic reparametrization of time.
We consider a diffeomorphism of the torus $\TTT^{d}$ of the form
\begin{equation}\label{32}
\theta=\f+\oo\al(\f), \quad \f\in\TTT^{d}, \quad \al(\f)\in\RRR,
\end{equation}
where $\al$ is a small real valued function, $2\pi-$periodic in all its arguments. The induced
linear operator on the space of functions is
\begin{equation}\label{332}
(\TT_{3} h)(\f,y):=h(\f+\oo\al(\f),y),\quad
{\rm with \;\; inverse} \quad
(\TT_{3}^{-1}v)(\theta,y)=v(\theta+\oo\widehat{\al}(\theta),y),
\end{equation}
where $\f=\theta+\oo\widehat{\al}(\theta)$ is the inverse diffeomorphism of
$\theta=\f+\oo\al(\f)$. We extend the operator
\begin{equation}\label{eq:3.15}
{\TT_{3}} : {\bf H}^{s}\to {\bf H}^{s}, \quad ({\TT_{3}}{\bf h})(\f,x)=\left(
\begin{matrix} (\TT_{3} h)(\f,x) \\ (\TT_{3} \bar{h})(\f,x)\end{matrix}
\right).
\end{equation}

 By conjugation, we have that the differential operator becomes
\begin{equation}\label{33}
\TT_{3}^{-1}\oo\cdot\del_{\f}\TT_{3}=\rho(\theta)\oo\cdot\del_{\theta},\qquad
\TT_{3}^{-1}\del_{y}\TT_{3}=\del_{y}, \qquad \rho(\theta):=\TT_{3}^{-1}(1+\oo\cdot\del_{\f}\al).
\end{equation}
We have obtained ${\TT_{3}}^{-1}\calL_{2}{\TT_{3}}=\rho\calL_{3}$ with $\calL_{3}$ as in (\ref{elle})
where 
\begin{equation}\label{eq:3.16}
\begin{aligned}
&1+a_{2}^{(3)}(\theta):=(\TT_{3}^{-1}(1+a_{2}^{(2)}))(\te),\\
&\rho(\theta) {a}^{(3)}_{j}(\theta,y):=(\TT_{3}^{-1}a^{(2)}_{j})(\te,y), \quad \;\;\;
\rho(\theta) {b}^{(3)}_{j}(\theta,y):=(\TT_{3}^{-1}b^{(2)}_{j})(\te,y), \;\; \quad i=0,1.
\end{aligned}
\end{equation}
We look for solutions $\al$ such that the coefficients of the highest order derivatives
($i\oo\cdot\del_{\theta}$ and $\del_{yy}$) are proportional, namely
\begin{equation}\label{35}
(\TT_{3}^{-1}(1+a_{2}^{(2)}))(\theta)=m\rho(\theta)=m\TT_{3}^{-1}(1+\oo\cdot\del_{\f}\al)
\end{equation}
for some constant $m$, that is equivalent to require that
\begin{equation}\label{36}
1+a_{2}^{(2)}(\f)=m(1+\oo\cdot\del_{\f}\al(\f)),
\end{equation}
By setting 
\begin{equation}\label{37}
m=\frac{1}{(2\pi)^{d}}\int_{\TTT^{d}}(1+a_{2}^{(2)}(\f))d\f,
\end{equation}
we can find the (unique) solution of $(\ref{36})$ with zero average
\begin{equation}\label{38}
\al(\f):=\frac{1}{m }(\oo\cdot\del_{\f})^{-1}(1+a_{2}^{(2)}-m)(\f),
\end{equation}
where $(\oo\cdot\del_{\f})^{-1}$ is defined by linearity
$$
(\oo\cdot\del_{\f})^{-1}e^{i\ell\cdot\f}:=\frac{e^{i\ell\cdot\f}}{i\oo\cdot\ell}, \; \ell\neq0, \quad
(\oo\cdot\del_{\f})^{-1}1=0.
$$
thanks to this choice of $\al$  
we have ${\TT}_{3}^{-1}\calL_{2}{\TT}_{3}=\rho \calL_{3}$ with $1+a_{2}^{(3)}(\theta)=m$.

\noindent
First of all, note that the bounds (\ref{eq:3.2.4})  on the coefficient $m$ in (\ref{37}) follow
by the (\ref{eq:3.10}) for $a_{2}^{(2)}$.
Moreover the function $\al({\f})$ defined in (\ref{38}) satisfies the tame estimates:
\begin{subequations}\label{eq:B24}
\begin{align}
|\al|^{\infty}_{s }&\leq\e \g_{0}^{-1}C(s)(1+||{\bf u}||_{s+d+\gots_{0}+2}),\label{eq:B24a}\\
|d_{{\bf u} }\al({\bf u})[{\bf h}]|^{\infty}_{s }&\leq \e \g_{0}^{-1}C(s)
(||{\bf h}||_{s+d+\gots_{0}+2}+
||{\bf u}||_{s+d+\gots_{0}+2}||{\bf h}||_{d+\gots_{0}+2}),\label{eq:B24b}\\
|\al|^{\infty}_{s,\g} &\leq \e\g_{0}^{-1}C(s)
(1+||{\bf u}||_{s+d+\gots_{0}+2,\g}).\label{eq:B24c}
\end{align}
\end{subequations}
Since $\oo=\la\bar{\oo}$ and by (\ref{dio}) one has
$|\bar\oo\cdot\ell|\geq3\g_{0}|\ell|^{-d}$, $\forall\;\ell\neq0$, then one has the (\ref{eq:B24a}).
One can prove similarly the (\ref{eq:B24c}) by using 
(\ref{eq:3.10a}), (\ref{eq:3.2.4}) and the fact $(\oo\cdot\ell)^{-1}=\la^{-1}(\bar\oo\cdot\ell)^{-1}$.
To prove (\ref{eq:B24b}) we compute
\begin{equation}\label{eq:B26}
d_{{\bf u} }\al({\f ;{\bf u}})[{\bf h}]=(\la\bar\oo\cdot\del_{\f})^{-1}
\left(\frac{ d_{\bf u} (1+a_{2}^{(2)}({\bf u}))[{\bf h}]m-(1+a_{2}^{(2)})d_{{\bf u} }m({\bf u})[{\bf h}]}{m^{2}}
\right)
\end{equation}
and  use the estimates (\ref{eq:3.10a}), (\ref{eq:3.10b}) and (\ref{eq:3.2.4}). Finally, the diffeomorphism (\ref{32}) is well-defined if
$|\al|^{\infty}_{1 }\leq1/2$. This is implied by (\ref{eq:B24a}) and (\ref{eq:3.2.1b}) for $\e$ small enough.

\noindent
The inverse diffeomorphism $\theta\to\theta+\oo\widehat{\al}(\theta)$ of (\ref{32}) satisfies 
the same estimates in (\ref{eq:B24}) with $d+\gots_{0}+3$.
The (\ref{eq:B24a}), (\ref{eq:B24c}) on $\widehat{\al}$ follow by the bounds (\ref{A18}), (\ref{A19}) in Lemma \ref{change}
and (\ref{eq:B24a}), (\ref{eq:B24c}). 
As in the second step the estimate on $d_{{\bf u}}\widehat{\al}({\bf u})[{\bf h}]$
follows
by the chain rule using Lemma \ref{change}(iii), (\ref{A6}),  (\ref{eq:B24a}), (\ref{eq:B24b})
on $\al$ and   (\ref{A2}) with $a=d+\gots_{0}+3$, $b=d+\gots_{0}+1$
and $p=s-1$, $q=2$
one has the (\ref{eq:B24b}) for $\widehat{\al}$.

\noindent
We claim that the operators $\TT_{3}({\bf u})$ and $\TT_{3}^{-1}({\bf u})$ defined in (\ref{332}), satisfy for any 
${\bf g},{\bf h}\in {\bf H}^{s}$ the (\ref{eq:B15}) with $\ka_{3}:=d+\gots_{0}+3$.
Indeed to prove estimates (\ref{eq:B15a}),  we apply Lemma \ref{change}(ii) and the estimates
(\ref{eq:B24a}), (\ref{eq:B24c}) on $\al$ and $\widehat{\al}$  obtained above.
Now, since
\begin{equation}\label{equ:B33}
\begin{aligned}
d_{{\bf u} }(\TT_{3}&({\bf u}){\bf g})[{\bf h}]=\TT_{3}({\bf u})(\oo\cdot\del_{\f}{\bf g})d_{{\bf u} }\al({\bf u})[{\bf h}]
\end{aligned}
\end{equation}
then (\ref{A7}), (\ref{eq:B24b}) and (\ref{eq:B15a}),  imply (\ref{eq:B15b}). Reasoning in the same way
one has that (\ref{eq:B24a}), (\ref{eq:B15b}) imply (\ref{eq:B15b}) on $\TT_{3}^{-1}$.

\noindent
By the (\ref{33}) one has $\rho=1+\TT_{3}^{-1}(\oo\cdot\del_{\f}\al)$. By using 
the (\ref{A21a}), (\ref{A21b}), the bounds (\ref{eq:B24}) on $\al$ and (\ref{eq:3.2.1b}) one can prove
\begin{subequations}\label{eq:B34}
\begin{align}
|\rho-1|^{\infty}_{s,\g}&\leq
 \e\g_{0}^{-1}C(s)(1+||{\bf u}||_{s+d+\gots_{0}+4,\g}) \label{eq:B34b}\\
|d_{{\bf u} }\rho({\bf u})[{\bf h}]|^{\infty}_{s}&\leq 
\e\g_{0}^{-1}C(s)(||{\bf h}||_{s+d+\gots_{0}+3}+
||{\bf u}||_{s+d+\gots_{0}+4}||{\bf h}||_{d+\gots_{0}+3}).\label{eq:B34c}
\end{align}
\end{subequations}
Bounds (\ref{eq:3.10}) on the coefficients $a_{j}^{(3)}, b_{j}^{(3)}$ follows, with $\ka_{3}:=d+\gots_{0}+5$, by using the
(\ref{eq:B34}) on $\rho$, the (\ref{eq:B15}) on $\TT_{3}$ and $\TT^{-1}_{3}$,
the (\ref{A5})-(\ref{A7}) and the condition (\ref{eq:3.2.1b}).
\begin{rmk}\label{rmk3}
 Note that $\al$ is a real function and $\al\in X^{0}$, then the operators $\TT_{3}$ and $\TT_{3}^{-1}$ map $X^{0}\to X^{0}$ and $Y^{0}\to Y^{0}$.
 Moreover we have that $m\in \RRR$ , $a^{(3)}_{0}, b^{(3)}_{0}\in Y^{0}$, 
 while
$a^{(3)}_{1}, b^{(3)}_{1}\in X^{0}$. Then 
then by Lemma \ref{rmk:rev},  one has that
the operator $\calL_{3}$ is reversible. 
\end{rmk}

In the following we rename $y=x$ and $\theta=\f$

\bigskip
\noindent
{\bf Step 4. Descent Method: conjugation by multiplication operator}
\bigskip

The aim of this section is to conjugate the operator $\calL_{3}$ to an operator
$\calL_{4}$ which has zero on the diagonal of the first order spatial differential operator.

 We consider an operator
of the form
\begin{equation}\label{eq:3.5.1}
\TT_{4}:=\left(\begin{matrix}1+z(\f,x) & 0 \\ 0 & 1+\bar{z}(\f,x) \end{matrix}\right),
\end{equation}
where $z : \TTT^{d+1}\to \CCC$ is small enough so that $\TT_{4}$ is invertible.  By a direct calculation we have that
$\calL_{4}$ has the form (\ref{elle}) where the second order coefficients are those of $\calL_{3}$ while\footnote{We use $\TT_{4}$ to cancel $a_{1}^{(4)}$, then to avoid apices we rename the remaining coefficients coherently with the definition of $\calL_{4}$.}
\begin{equation}\label{eq:3.5.10}
\begin{aligned}
a_{1}^{(4)}(\f,x)&:=2m\frac{z_{x}(\f,x)}{1+z(\f,x)}+a^{(3)}_{1}(\f,x),\quad
q_{2}(\f,x)\equiv a^{(4)}_{0}(\f,x):=\frac{-i(\oo\cdot\del_{\f}z)(\f,x)+mz_{xx}}{1+z(\f,x)}+a^{(3)}_{0}(\f,x),\\
q_{1}(\f,x)&\equiv b_{1}^{(4)}(\f,x):=b^{(3)}_{1}(\f,x)\frac{1+\bar{z}(\f,x)}{1+z(\f,x)}, \quad
 q_{3}(\f,x)\equiv b^{(4)}_{0}(\f,x):=b^{(3)}_{0}(\f,x)\frac{1+\bar{z}(\f,x)}{1+z(\f,x)}.
\end{aligned}
\end{equation}
We look for $z(\f,x)$ such that $a_{1}^{(4)}\equiv0$. If we look for solutions of the form
$1+z(\f,x)=\exp(s(\f,x))$ we have that $a_{1}^{(4)}=0$ becomes
\begin{equation}\label{eq:3.25}
\begin{aligned}
({\rm Re}(s))_{x}(\f,x)&=-\frac{1}{2m}{\rm Re}(a^{(3)}_{1})(\f,x), \qquad \;\;
({\rm Im}(s))_{x}(\f,x)=-\frac{1}{2m}{\rm Im}(a^{(3)}_{1})(\f,x),
\end{aligned}
\end{equation}
that have unique (with zero average in $x$) solution
\begin{equation}\label{eq:3.26}
\begin{aligned}
({\rm Re}s)(\f,x)&=-\frac{1}{2m}\del_{x}^{-1}{\rm Re}(a^{(3)}_{1})(\f,x), \qquad \;\;
({\rm Im}s)(\f,x)=-\frac{1}{2m}\del_{x}^{-1}{\rm Im}(a^{3}_{1})(\f,x)
\end{aligned}
\end{equation}
where $\del_{x}^{-1}$ is defined in (\ref{30}).  

\noindent
The function $s$ defined in (\ref{eq:3.26}) satisfies the following tame estimates:
\begin{subequations}\label{eq:B40}
\begin{align}
||s||_{s,\g}&\leq \e C(s)(1+||{\bf u}||_{s+d+\gots_{0}+5,\g}),\label{eq:B40b}\\
||d_{{\bf u} }s({\bf u})[{\bf h}]||_{s}&\leq\e C(s)(||{\bf h}||_{s+d+\gots_{0}+4}+
||{\bf u}||_{s+d+\gots_{0}+5}||{\bf h}||_{d+\gots_{0}+4}).\label{eq:B40c}
\end{align}
\end{subequations}
The (\ref{eq:B40}) follow by (\ref{eq:3.2.4}),  used to estimate $m$,
the estimates (\ref{eq:3.10}), on the coefficient of $a_{1}^{(3)}$, and (\ref{eq:3.2.1b}). Since by definition one has
$$
z(\f,x)=\exp(s(\f,x))-1,
$$
clearly the function $z$ satisfies the same estimates (\ref{eq:B40b})-(\ref{eq:B40c}).

\noindent
The estimates (\ref{eq:B40b})-(\ref{eq:B40c}) on the function $z(\f,x)$ imply directly
the tame estimates in (\ref{eq:B15}) on the  operator $\TT_{4}$  defined in (\ref{eq:3.5.1}).
The bound (\ref{eq:B15a}) on the operator $\TT_{4}^{-1}$
follows in the same way.
In order to prove the (\ref{eq:B15b}) we
 note that
\begin{equation*}
d_{{\bf u} }\TT_{4}^{-1}({\bf u})[{\bf h}]=-\TT_{4}^{-1}({\bf u})d_{{\bf u} }\TT_{4}({\bf u})[{\bf h}]\TT_{4}^{-1}{(\bf u)},
\end{equation*}
then, using the (\ref{eq:3.2.1b}) and the (\ref{eq:B15b}) on $\TT_{4}$
we get the (\ref{eq:B15b}) on $\TT_{4}^{-1}$.
\noindent
We show that the coefficients in (\ref{eq:3.5.10}), for $i=1,2,3$ satisfy
the tame estimates in (\ref{eq:3.10}) with $\ka_{4}=d+\gots_{0}+7$
%
that simply are the (\ref{eq:3.2.6a}), (\ref{eq:3.2.6b}).
The strategy to prove the tame bounds on $q_{i}$ is the same used in (\ref{eq:B20}) on $a_{1}^{(2)}$.
Collecting together the loss of regularity at each step one gets $\h_{1}$ as in (\ref{eq:3.2.0}).

\begin{rmk}\label{rmk:3.7}
Since $a^{(3)}_{1}\in X^{0}$, then
$s(\f,x)\in Y^{0}$, 
so that the operator $\TT_{4}$ does not change the parity properties of functions. This implies
that the operator $\calL_{4}$, defined in (\ref{eq:3.5.9}), is reversible.
\end{rmk}

The several steps performed in the previous sections  (semi)-conjugate the linearized operator
$\calL$ to the operator $\calL_{4}$ defined in (\ref{eq:3.5.9}), namely
\begin{equation}\label{tra}
\calL=\VV_{1}\calL_{4}\VV_{2}^{-1}, \quad \VV_{1}:=\TT_{1}{\TT_{2}} {\TT_{3}}\rho\TT_{4}, \quad
\VV_{2}=\TT_{1} {\TT_{2}}{\TT_{3}}\TT_{4}.
\end{equation}
where $\rho$ is the multiplication operator by the function $\rho$ defined in (\ref{33}).
Now by Lemma \ref{lem5}, the operators
$\VV_{1}$ and $\VV_{2}$ defined in (\ref{tra}) satisfy, using (\ref{eq:3.2.1b}), the (\ref{eq:3.2.3}). Note that
we used that $\h_{1}>d+2\gots_{0}+7$.
The estimates in (ii) and (iii) heve been already proved, hence
 the proof of  Lemma \ref{lem:3.8} has been completed.
\EP

The following Lemma is a consequence of the discussion above.
\begin{lemma}\label{lem:3.9}
Under the Hypotheses of Lemma \ref{lem:3.8} possibly with smaller $\epsilon_{0}$, 
if 
(\ref{eq:3.2.1b}) holds, one has that the$\TT_i$ $i\neq3 $ identify operators $\TT_{i}(\f)$, 
 of the phase space ${\bf H}^{s}_{x}:={\bf H}^{s}(\TTT)$. Moreover they are invertible
and the following estimates hold for $\gots_{0}\leq s\leq q-\h_{1}$:
\begin{subequations}
\begin{align}
||(\TT_{i}^{\pm1}(\f)-\id){\bf h}||_{{\bf H}^{s}_{x}}&\leq \e 
C(s)(||{\bf h}||_{{\bf H}^{s}_{x}}+||{\bf u}||_{s+d+2\gots_{0}+4}
||{\bf h}||_{{\bf H}^{1}_{x}}), \quad i=1,2,4,\label{eq:3.93f}
\end{align}
\end{subequations}
\end{lemma}
\prova
$\TT_{1}$ and $\TT_{4}$ are multiplication operators then,
it is enough to perform the proof  on any component
$(\TT_{i})_{\s}^{\s'}$, for $\s,\s'=\pm1$ and $i=1,4$, that are simply multiplication
operators from $H_{x}^{s} \to H_{x}^{s}$. One has
\begin{equation}\label{eq:3.91}
\begin{aligned}
||(\TT_{i})_{\s}^{\s'}(\f)h||_{H^{s}_{x}}&\stackrel{(\ref{A5})}{\leq} C(s)(
||(\TT_{i})_{\s}^{\s'}(\f)||_{H^{s}_{x}}||h||_{H^{1}_{x}}+
||(\TT_{i})_{\s}^{\s'}(\f)||_{H^{1}_{x}}||h||_{H^{s}_{x}})\\
&{\leq}
||(\TT_{i})_{\s}^{\s'}||_{{s+\gots_{0}}}||h||_{H^{s}_{x}}+
||(\TT_{i})_{\s}^{\s'}||_{1+\gots_{0}}||h||_{H^{s}_{x}})\\
&\stackrel{(\ref{eq:B15a})}{\leq} C(s)(||h||_{H^{s}_{x}}+||{\bf u}||_{s+\gots_{0}+2}||h||_{H^{1}_{x}}),
\end{aligned}
\end{equation}
where we used also (\ref{eq:3.2.1b}). 
 In the same way one can show that
\begin{equation}\label{eq:3.92}
||((\TT_{i})_{\s}^{\s'}(\f,\cdot)-\id)h||_{s}\leq
\e C(s)(||h||_{H^{s}_{x}}+||{\bf u}||_{s+\gots_{0}+2}||h||_{H^{1}_{x}}).
\end{equation}
and hence
the bound (\ref{eq:3.93f}) follow. Note that we used the simple fact that
given a function $v \in H^{s}(\TTT^{d+1};\CCC)$ then $||v(\f)||_{H^{s}_{x}}\leq C ||v||_{s+\gots_{0}}$.
Now, for fixed $\f\in\TTT^{d}$ one has $\TT_{2}(\f)h(x):=h(x+\x(\f,x))$.
We can bound, by using the (\ref{A20a}) on the change of  variable $\TTT\to\TTT$, $x\to x+\x(\f,x)$, 
\begin{equation}\label{eq:3.94}
\begin{aligned}
||\TT_{2}(\f)h||_{H^{s}_{x}}&\leq C(s)(||h||_{H^{s}_{x}}+|\x(\f)|_{W^{s,\infty}(\TTT)}||h||_{H^{1}_{x}})\\
&\stackrel{(\ref{eq:3.50a})}{\leq}
C(s)(||h||_{H^{s}_{x}}+||{\bf u}||_{s+\gots_{0}+2}
||h||_{H^{1}_{x}})
\end{aligned}
\end{equation}
where we have used also the fact $|\x(\f)|_{W^{s,\infty}(\TTT)}\leq |\x|^{\infty}_{s+\gots_{0}}$.
One can prove (\ref{eq:3.93f}) by using (\ref{A20b}), (\ref{eq:3.2.1b}) and (\ref{eq:3.50a}).
The estimates 
(\ref{eq:3.93f}) hold for $\TT_{2}^{-1}(\f) : h(y) \to h(y+\widehat{\x}(\f,y))$
thanks to the (\ref{eq:3.52}).
\EP

Note that the fact that $\TT_{3}$ maps $H^{s}_{x}\to H^{s}_{x}$ is trivial.

\setcounter{equation}{0}
\section{The diagonalization algorithm: KAM reduction}
\label{sec:4}

In this section we diagonalize the operator 
$\calL_{4}$ in (\ref{eq:3.5.9}) in Section \ref{sec:3}. 
In order to implement our procedure we pass to Fourier coefficients and  introduce an "off diagonal decay norm" which is stronger that the standard operatorial one. We  also 
define the reversibility properties of the operators, in terms of the Fourier coefficients. 

Consider the bases $\{e_{k}=e^{i\ell\cdot\f}\sin jx : k=(\ell,j)\in\ZZZ^{d}\times\NNN\}$ and
 $\{e_{k}=e^{i\ell\cdot\f}\cos jx : k=(\ell,j)\in\ZZZ^{d}\times\ZZZ_+\}$  for functions  which are odd (resp. even)  in $x$.
Then any linear operator $A : {\bf G}^{0}_1\to {\bf G}_2^{0}$, 
where ${\bf G}^{0}_{1,2}= {\bf X}^{0}, {\bf Y}^{0}$, ${\bf Z}^{0}$, can be represented by 
an infinite dimensional matrix
$$
A:=(A_{i}^{i'})_{i,i'\in\CC\times\ZZZ_+\times\ZZZ^{d}}, \quad (A_{\s}^{\s'})_{k}^{k'}=(Ae_{k'},e_{k})_{L^{2}(\TTT^{d+1})},\quad 
(A_{\s}^{\s'})u=\sum_{k,k'}(A_{\s}^{\s'})_{k}^{k'}u_{k'}e_{k},
$$
where $(\cdot,\cdot)_{L^{2}(\TTT^{d+1})}$ is the usual scalar product on $L^{2}$,
we are denoting $i=(\s,k)=(\s,j,p)\in\CC\times\ZZZ_+\times\ZZZ^{d}$ and $\CC:=\left\{+1,-1\right\}$.

In the case of functions which are odd in $x$ we set the {\em extra } matrix coefficients  (corresponding to $j=0$) to zero.

\begin{defi}
({\bf s-decay norm}). Given an infinite dimensional matrix 
$A:=(A_{i}^{i'})_{i,i'\in\CC\times\ZZZ_+\times\ZZZ^{d}}$ we define the norm of off-diagonal decay
\begin{equation}\label{decay}
\begin{aligned}
|A|^{2}_{s}:=\sup_{\s,\s'\in \CC}|A_{\s}^{\s'}|^{2}_{s}&:=
\sup_{\s,\s'\in\CC}\sum_{h\in\ZZZ_+\times\ZZZ^{d}}\langle h\rangle^{2s}
\sup_{k-k'=h}|A_{\s,k}^{\s',k'}|^{2}
\end{aligned}
\end{equation}
If one has that $A:=A(\la) $ for $\la\in\Lambda\subset\RRR$, we define
\begin{equation}\label{2.1}
|A|_{s}^{sup}:=\sup_{\la\in\Lambda}|A(\la)|_{s}, \quad 
|A|^{lip}_{s}:=\sup_{\la_{1}\neq\la_{2}}\frac{|A(\la_{1})-A(\la_{2})|_{s}}{|\la_{1}-\la_{2}|}, \quad
|A|_{s,\g}:=|A|_{s}^{sup}+\g|A|^{lip}_{s}.
\end{equation}
\end{defi}

The decay norm we have introduced in (\ref{decay}) is suitable for the problem we are studying.
 Note that
\begin{equation*}
\;\;\;
\forall\; s\leq s' \;\; \Rightarrow \;\; |A_{\s}^{\s'}|_{s}\leq|A_{\s}^{\s'}|_{s'}.
\end{equation*}
Moreover norm (\ref{decay}) gives information on the polynomial off-diagonal decay of the matrices,
indeed
\begin{equation}\label{decay2}
\begin{aligned}
&|A_{\s,k}^{\s,k'}|\leq \frac{|A_{\s}^{\s'}|_{s}}{\langle k-k' \rangle^{s}}, \quad\forall\; k,k'\in\ZZZ_+\times\ZZZ^{d},\quad {\rm and} \quad 
|A_{i}^{i}|\leq |A|_{0}, \quad |A_{i}^{i}|^{lip}\leq |A|_{0}^{lip}.
\end{aligned}
\end{equation}

We have the following important result:

\begin{theorem}\label{KAMalgorithm}
Let $f\in C^{q}$ satisfy the Hypotheses of Proposition \ref{teo4} with $q>\h_{1}+\be+\gots_{0}$
where $\h_{1}$ defined in (\ref{eq:3.2.0}) 
and $\be=7\tau+5$ for some $\tau>d$.
Let $\g\in(0,\g_{0})$,
$\gots_{0}\leq s\leq q-\h_{1}-\be$ and ${\bf u}(\la)\in{\bf X}^{0}$ be a family of functions depending 
on a Lipschitz way on a parameter
$\la\in\Lambda_{o}\subseteq\Lambda:[1/2,3/2]$. 
Assume that
\begin{equation}\label{eq:4.2}
||{\bf u}||_{\gots_{0}+\h_{1}+\be, \Lambda_{o},\g}\leq1.
\end{equation}
Then there exist constants $\epsilon_{0}$, $C$, depending only on the data of the problem,
such that,
if $\e\g^{-1}\leq\epsilon_{0}$,
then there exists a sequence of purely imaginary numbers as in Proposition \ref{teo2},
namely
\begin{equation}\label{eq:4.4}
\mu^{\infty}_{h}:= 
\mu_{\s,j}^{\infty}(\la):=\mu_{\s,j}^{\infty}(\la,{\bf u})=-\s i mj^{2}+ r_{\s,j}^{\infty}, \quad \forall\;h=(\s,j)\in \CC\times\NNN, \quad 
\forall\;\la\in\Lambda,
\end{equation}
where $m$ is defined in (\ref{37}) with 
\begin{equation}\label{eq:4.5}
|r_{\s,j}^{\infty}|_{\g}\leq \e C, \quad \forall\; \s\in\CC,\; j\in\NNN.
\end{equation}
and
such that, for any $\la\in \Lambda_{\infty}^{2\g}({\bf u})$, defined in (\ref{martina10}),
there exists a bounded, invertible linear operator
$\Phi_{\infty}(\la) : {\bf H}^{s}\to {\bf H}^{s}$, with bounded inverse $\Phi_{\infty}^{-1}(\la)$,
such that
\begin{equation}\label{eq:4.6}
\begin{aligned}
\calL_{\infty}(\la):=\Phi_{\infty}^{-1}(\la)\circ\calL_{4}\circ\Phi_{\infty}(\la)&=
\la\bar{\oo}\cdot\del_{\f}\id+i \DD_{\infty},\\
{\rm where}\;\;\;\;\;\;\;\;\;
 \DD_{\infty}:={\rm diag_{h\in\CC\times\NNN}}\{\mu_{h}(\la)\},& \quad 
\end{aligned}
\end{equation}
with $\m_{h}$ definend in (\ref{eq:4.4}) and $\calL_4$ in \eqref{eq:3.5.9}.
Moreover, the transformations $\Phi_{\infty}(\la)$, $\Phi_{\infty}^{-1}$ satisfy
\begin{equation}\label{eq:4.8}
|\Phi_{\infty}(\la)-\id|_{s,\Lambda_{\infty}^{2\g},\g}+
|\Phi_{\infty}^{-1}(\la)-\id|_{s,\Lambda_{\infty}^{2\g},\g}\leq
\e \g^{-1} C(s)(1+||{\bf u}||_{s+\h_{1}+\be,\Lambda_{o},\g}).
\end{equation}
In addition to this, for any $\f\in\TTT^{d}$, for any $\gots_{0}\leq s\leq q-\h_{1}-\be$
the operator
$\Phi_{\infty}(\f) : {\bf X}^{s}_{x}\to{\bf X}^{s}_{x}$ is an invertible operator
of the phase space ${\bf X}_{x}^{s}:={\bf X}^{s}(\TTT)$ with inverse
$(\Phi_{\infty}(\f))^{-1}:=\Phi_{\infty}^{-1}(\f)$ and
\begin{equation}\label{eq:4.9}
||(\Phi_{\infty}^{\pm1}(\f)-\id){\bf h}||_{{\bf H}^{s}_{x}}\leq
\e\g^{-1} C(s)(||{\bf h}||_{{\bf H}^{s}_{x}}+||{\bf u}||_{s+\h_{1}+\be+\gots_{0}}
||{\bf h}||_{{\bf H}^{1}_{x}}).
\end{equation}
\end{theorem}

\begin{rmk}
It is important to note that thanks to Reversibility Hypothesis \ref{hyp2}, the operator
$\calL_{\infty} : {\bf X}^{0} \to {\bf Z}^{0}$ i.e. it is reversible. 
\end{rmk}
 
The main point of the Theorem \ref{KAMalgorithm} is that 
the bound on the low norm of $u$ in (\ref{eq:4.2}) guarantees the bound 
on \emph{higher} norms (\ref{eq:4.8}) for the transformations $\Phi_{\infty}^{\pm1}$. This is 
fundamental in order  to get the estimates on the inverse of $\calL$ in high norms.

Moreover, the definition (\ref{martina10}) of the set where the second Melnikov conditions
hold, depends only on the final eigenvalues. Usually in KAM theorems, the non-resonance conditions
have to be checked, inductively, at each step of the algorithm. This 
formulation, on the contrary, allow us to discuss the measure estimates only once. 
Indeed, the functions $\m_{h}(\la)$ are well-defined even if
$\Lambda_{\infty}=\emptyset$, so that, we will perform the measure estimates as the last step of the proof of Theorem \ref{teo1}.

\subsection{Functional setting and notations}

\subsubsection{The off-diagonal decay norm }

Here we want to show some important properties of the norm $|\cdot|_{s}$. Clearly the same results hold for the 
norm $|\cdot|_{{\bf H}^{s}}:=|\cdot|_{H^{s}\times H^{s}}$. Moreover we will introduce 
some charatterization  of the operators we have to deal with during the diagonalization procedure.
 
First of all we have following classical results.

\begin{lemma}\label{ bubbole}{\bf Interpolation.} For all $s\geq s_{0}>(d+1)/2$ 
there are $C(s)\geq C(s_{0})\geq1$ 
such that 
if $A=A(\la)$ and $B=B(\la)$ depend on the parameter
 $\la\in\Lambda\subset\RRR$
in a Lipschitz way, then
\begin{subequations}
\begin{align}
|AB|_{s,\g}&\leq C(s)|A|_{s_0,\g}|B|_{s,\g}
+C(s_{0})|A|_{s,\g}|B|_{s_0,\g},\label{eq:2.11a}\\
|AB|_{s,\g}&\leq C(s)|A|_{s,\g}|B|_{s,\g}.\label{eq:2.11b}\\
||Ah||_{s,\g}&\leq C(s)(|A|_{s_0,\g}||h||_{s,\g}
+|A|_{s,\g}||h||_{s_0,\g}),\label{eq:2.13b}\end{align}
\end{subequations}
\end{lemma}

Lemma \ref{ bubbole} implies that for any $n\geq0$ one has
\begin{equation}\label{eq:2.12}
|A^{n}|_{s_{0},\g}\leq [C(s_{0})]^{n-1}|A|^{n}_{s_{0},\g}, \quad {\rm and} \quad 
|A^{n}|_{s,\g}\leq n[C(s_{0})]^{n-1}C(s)|A|_{s,\g}, \quad \forall\; s\geq s_{0}.
\end{equation}

The following Lemma shows how to invert linear operators which 
are ''\emph{near}'' to the identity in norm $|\cdot|_{s}$. 

\begin{lemma}\label{inverse} Let $C(s_0)$ be as in Lemma \ref{ bubbole}. Consider an operator of the form $\Phi=\id+\Psi$
where $\Psi=\Psi(\la)$ depends in a Lipschitz way on $\la\in\Lambda\subset\RRR$. 
Assume that $C(s_{0})|\Psi|_{s_{0},\g}\leq1/2$. Then $\Phi$ is invertible and,
for all $s\geq s_{0}\geq (d+1)/2$, 
\begin{equation}\label{eq:2.14}
 \quad |\Phi^{-1}|_{s_{0},\g}\leq 2,\quad
|\Phi^{-1}-\id|_{s,\g}\leq C(s)|\Psi|_{s,\g}
\end{equation}
Moreover, if one  has $\Phi_{i}=\id+\Psi_{i}$, $i=1,2$ such that $C(s_{0})|\Psi_{i}|_{s_0,\g}\leq1/2$, then
\begin{equation}\label{eq:2.15}
|\Phi^{-1}_{2}-\Phi^{-1}_{1}|_{s,\g}\leq C(s)\left(|\Psi_{2}-\Psi_{1}|_{s,\g}+(|\Psi_{1}|_{s,\g}+|\Psi_{2}|_{s,\g})
|\Psi_{2}-\Psi_{1}|_{s_0,\g}\right).
\end{equation}

\end{lemma}

\prova
One has that $(\id +\Psi)^{-1}=\sum_{k\geq0}\frac{(-1)^{k}}{k!}\Psi^{k}$, then by (\ref{eq:2.12}) we get
bounds (\ref{eq:2.14}). Now, we can note that
\begin{equation*}
\begin{aligned}
|\Phi^{-1}_{2}-\Phi^{-1}_{1}|_{s,\g}&= |\Phi_{1}^{-1}(\Psi_{1}-\Psi_{2})\Phi_{2}^{-1}|_{s,\g}
\stackrel{{(\ref{eq:2.11a})}}{\leq}  C(s)|\Phi_{1}^{-1}|_{s_0,\g}|\Psi_{1}-\Psi_{2}|_{s_0,\g}|\Phi_{2}^{-1}|_{s,\g}\\
&+C(s)|\Phi_{1}^{-1}|_{s_0,\g}|\Psi_{1}-\Psi_{2}|_{s,\g}|\Phi_{2}^{-1}|_{s_0,\g}
+C(s)|\Phi_{1}^{-1}|_{s,\g}|\Psi_{1}-\Psi_{2}|_{s_0,\g}|\Phi_{2}^{-1}|_{s_0,\g}\\
&\stackrel{(\ref{eq:2.14})}{\leq}
C(s)(|\Psi_{1}-\Psi_{2}|_{s,\g}+(|\Psi_{1}|_{s,\g}+|\Psi_{2}|_{s,\g})|\Psi_{1}-\Psi_{2}|_{s_0,\g})
\end{aligned}
\end{equation*}	
that is the  (\ref{eq:2.15}).
\EP

\subsubsection{T\"opliz-in-time matrices}

We introduce now a special class of operators, the so-called 
{\em T\"opliz in time} matrices, i.e.
\begin{equation}\label{eq:2.16}
A_{i}^{i'}=A_{(\s,j,p)}^{(\s',j',p')}:=A_{\s,j}^{\s'j'}(p-p'), 
\quad {\rm for} \quad i,i'\in \CC\times\ZZZ_+\times\ZZZ^{d}.
\end{equation}
To simplify the notation in this case, we shall write
$A_{i}^{i'}=A_{k}^{k'}(\ell)$,
 $i=(k,p)=(\s,j,p)\in \CC\times\ZZZ_+\times \ZZZ^{d}$, 
$i'=(k',p')=(\s',j',p')\in \CC\times \ZZZ_+\times\ZZZ^{d}$,
 with  $k,k'\in  \CC\times\ZZZ_+$.

They are relevant because
one can identify the matrix $A$ with a one-parameter family of operators, acting on the space
${\bf H}^{s}_{x}$,
which depend on the time, namely
\begin{equation}\label{pozzo10}
A(\f):=(A_{\s,j}^{\s',j'}(\f))_{\substack{\s,\s'\in C \\ j,j'\in\ZZZ_+}}, \quad 
A_{\s,j}^{\s',j'}(\f):=\sum_{\ell\in\ZZZ^{d}}A_{\s,j}^{\s',j'}(\ell)e^{i\ell\cdot\f}.
\end{equation}

To obtain the stability result on the solutions we will strongly use this property.

\begin{lemma}\label{1.4} If $A$ is a T\"opliz in time matrix as in (\ref{eq:2.16}), and 
$\gots_{0}:=(d+2)/2$, then one has
\begin{equation}\label{aaaaa}
|A(\f)|_{s}\leq C(\gots_{0})|A|_{s+\gots_{0}}, \quad \forall \; \f\in\TTT^{d}.
\end{equation}

\end{lemma}

\prova We can note that, for any $\f\in \TTT^{d}$,
\begin{equation}
\begin{aligned}
|A(\f)|^{2}_{s}&:=\sup_{\s,\s'\in\CCC}\sum_{h\in\ZZZ_+}\langle h\rangle^{2s}\sup_{j-j'=h}|A_{\s,j}^{\s',j'}(\f)|^{2}\leq
C(\gots_{0})\sup_{\s,\s'\in\CC}\sum_{h\in\ZZZ_+}\langle h\rangle^{2s}\sup_{j-j'=h}
\sum_{\ell\in\ZZZ^{d}}|A_{\s,j}^{\s',j'}(\ell)|^{2}\langle\ell\rangle^{2\gots_{0}}\\
&\leq C(\gots_{0})\sup_{\s.\s'\in\CC}\sum_{h\in\ZZZ_+}\sup_{j-j'=h}|A_{\s,j}^{\s;j}(\ell)|^{2}\langle\ell,h\rangle^{2(s+\gots_{0})}
\leq C(\gots_{0})\sup_{\s,\s'\in\CC}\sum_{\substack{h\in\ZZZ_+ \\ \ell\in\ZZZ^{d}}}\sup_{j-j'=h}|A_{\s,j}^{\s',j'}(\ell)|^{2}
\langle\ell,h\rangle^{2(s+\gots_{0})}\\
&\stackrel{(\ref{decay})}{\leq} C(\gots_{0})|A|_{s+\gots_{0}}^{2},
\end{aligned}
\end{equation}
that is the assertion.
\EP

\begin{defi}\label{def:smooth} {\bf (Smoothing operator)}
Given $N\in\NNN$, we the define the \emph{smoothing operator} $\Pi_{N}$ as
\begin{equation}\label{smoothop}
(\Pi_{N}A)_{\s,j,\ell }^{\s',j',\ell'}=\left\{
\begin{aligned}
& A_{\s,j,\ell}^{\s',j',\ell'} \,,\quad |\ell-\ell'|\leq N,\\
& 0 \quad  {\rm otherwise}
\end{aligned}\right.
\end{equation} 
\end{defi}

\begin{lemma}

Let $\Pi_{N}^{\perp}:=\id-\Pi_{N}$,

 if $A=A(\la)$ is a Lipschitz family $\la\in\Lambda$, then
\begin{equation}\label{eq:2.22}
|\Pi_{N}^{\perp}A|_{s,\g}\leq N^{-\be}|A|_{s+\be,\g}, \quad \be\geq0.
\end{equation}
\end{lemma}

\prova
Note that one has,
\begin{equation}\label{eq:2.23}
\begin{aligned}
|\Pi_{N}^{\perp}A|_{s}^{2}&=N^{-2\be}\sup_{\s,\s'\in\CC}\sum_{\substack{h\in\ZZZ_+ \\ |\ell|> N}}
\sup_{j-j'=h}|A_{\s,j}^{\s',j'}(\ell)|^{2}\langle\ell,h\rangle^{2s} N^{2\be}\\
&\leq
N^{-2\be}\sup_{\s,\s'\in\CC}\sum_{\substack{h\in\ZZZ_+ \\ |\ell|>N}}\sup_{j-j'=h}|A_{\s,j}^{\s',j'}(\ell)|^{2}
\langle\ell,h\rangle^{2(s+\be)}
\leq N^{-2\be}|A|_{s+\be}^{2},
\end{aligned}
\end{equation}
The estimate on the Lipschitz norm follows similarly.
\EP

\begin{rmk}{\bf (Multiplication operator)} 
We have already seen that if 
the decay norm is finite the operator has a ''good'' out diagonal decay. Although this  this property is strictly stronger than just being bounded, this class contains many useful operators in particular multiplication ones. 
Indeed, 
let $\TT_{a} : G^{s}_1\to G^{s}_2$, where $G^{s}_{1,2}=X^{s},Y^{s},Z^{s}$, 
be the multiplication operator by a function
$a\in G^{s}$ with $G^{s}=X^{s},Y^{s},Z^{s}$, i.e $(\TT_{a}h)=ah$. Then one can check, in coordinates, that
it is represented by the matrix $T$ such that
\begin{subequations}
\begin{align}
|T|_{s}&\leq ||a||_{s}.\label{eq:2.8b}
\end{align}
\end{subequations}
Moreover, if $a=a(\la)$ is a Lipschitz family of functions,
\begin{equation}\label{eq:2.9}
|T|_{s,\g}\leq ||a||_{s,\g}.
\end{equation}

At the beginning of our algorithm we actually deal with multiplication operators, so that one should try to control the operator by using only the Sobolev norms of functions. Unfortunately, it is not possible
since the class of multiplication operators is not closed under our alghorithm. This is the reason 
we have introduced the decay norms that control decay in more general situations.  
\end{rmk}

\subsubsection{Matrix representation}

In this paragraph we  give a caratherizations of reversible operators in the Fourier space.
We need it to deal with a more general class of operators than the multiplication operators.

\begin{lemma}\label{lemFou} We have that, for $G^{s}=X^{s},Y^{s}, Z^{s}$,

\begin{equation}\label{eq:2.25}
R : G^{s}\to G^{s}\quad  \Leftrightarrow \quad R_{j}^{j'}(\ell)=\ol{R_{j}^{j'}(\ell)}, \; \forall \ell\in\ZZZ^{d},
\;\forall \; j,j'\in \ZZZ_+.
\end{equation}
Moreover,

\begin{equation}\label{eq:2.26}
R : X^{s}\to Z^{s} \quad \Rightarrow R_{j}^{j'}(\ell)=-\ol{R_{j}^{j'}(\ell)},
\; \forall \ell\in\ZZZ^{d}, \; \forall \; j,j'\geq1.
\end{equation}

\end{lemma}

\prova
 One can consider a function $a(\f,x)\in G^{s}$ where $G^{s}=X^{s},Y^{s},Z^{s}$, and develop it in a suitable basis $e_{\ell,j}$, $(\ell,j)\in\ZZZ^{d}\times\ZZZ_{+}$ (to fix the idea we can think $e_{\ell, j}=e^{i\ell \f}\sin jx$, that is the correct basis for $X^{s}$). One has that the coefficients of the function $a$ satisfies $a_{j}(\ell)=\ol{a_{j}(\ell)}$ for $G^{s}=X^{s},Y^{s}$ while $a_{j}(\ell)=-\ol{a_{j}(\ell)}$ if $G^{s}=Z^{s}$. Then (\ref{eq:2.25}) and (\ref{eq:2.26}) follow by  applying the definitions of reversibility or reversibility preserving in (\ref{rever}) and (\ref{rever2}).
 \EP

\begin{lemma}\label{lemrev5}
Consider operators  $A : {\bf G}^{s}\to {\bf G }^{s}$ 
with $G^{s}=X^{s}, Y^{s},Z^{s}$ of the form
$A:=(A_{\s}^{\s'})_{\s,\s'=\pm1}$, then
\begin{equation}\label{rmkrev5}
\left(\begin{matrix}A_{1}^{1} & A_{-1}^{1} \\  A_{-1}^{1} & A_{-1}^{-1} \end{matrix}\right)
\left(\begin{matrix} u \\ \bar{u} \end{matrix}\right)=\left(\begin{matrix} w \\ \bar{w}\end{matrix}\right)
\in {\bf G^{s}}, \quad {\rm for \;\; any} \;\; (u,\bar{u})\in {\bf G}^{s}
\end{equation}
if and only if
\begin{equation}\label{rmkrev6}
A_{\s,j}^{\s',j'}(\ell)=\ol{A_{\s,j}^{\s',j'}(\ell)}, \quad {\rm and} \quad 
\ol{A_{\s,j}^{\s',j'}(-\ell)}=A_{-\s,j}^{-\s',j'}(\ell), \quad \forall\;\; \s,\s'=\pm1, \; \ell\in\ZZZ^{d},\; j,j'\in \ZZZ_+.
\end{equation}
An operator $B :{\bf X^{s}}\to {\bf Z^{s}}$ if and only if
\begin{equation}\label{rmkrev7}
B_{\s,j}^{\s',j'}(\ell)=-\ol{B_{\s,j}^{\s',j'}(\ell)}, \quad {\rm and} \quad 
\ol{B_{\s,j}^{\s',j'}(-\ell)}=B_{-\s,j}^{-\s',j'}(\ell), \quad \forall\;\; \s,\s'=\pm1, \; \ell\in\ZZZ^{d},\; j,j'\geq1.
\end{equation}
\end{lemma}

\prova
 Lemma \ref{lemFou} implies only that 
$A_{\s,j}^{\s',j'}(\ell)=\ol{A_{\s,j}^{\s',j'}(\ell)}$. Since we need that
the complex conjugate of the first component of $A{\bf u}$, with ${\bf u}\in G^{s}$, is
equal to the second one,  the components of $A$ have to satisfy
\begin{equation}\label{rmkrev4}
\ol{A_{\s,j}^{\s',j'}(-\ell)}=A_{-\s,j}^{-\s',j'}(\ell), \quad \forall \s,\s'=\pm1, \ell\in\ZZZ^{d}, j,j'\in \ZZZ_+.
\end{equation}
In this case we say that the operator $A : {\bf G^{s}}\to {\bf G^{s}}$ is \emph{reversibility-preserving}.

Following the same reasoning we have that for \emph{reversible} operators the (\ref{rmkrev7}) hold.
\EP

\subsection{Reduction Algorithm}

We prove Theorem \ref{KAMalgorithm} by means of the following 
Iterative Lemma on the class of  linear operators
\begin{defi}
\begin{equation}\label{eq:4.11}
\oo\cdot\del_{\f}\id+\DD+\RR : \; {\bf X}^{0}\to {\bf Z}^{0},
\end{equation}
where $\oo=\la\bar\oo$, and
\begin{equation}\label{eq:4.12}
\begin{aligned}
\DD=(-i\s m(\la,{\bf u}(\la))D^{2} )_{\s=\pm1}, \qquad \RR=E_{1}D+E_{0} 
\end{aligned}
\end{equation}
with $D:={\rm diag}_{j\in\NNN}\{ j\}$, and where, if we write $k=(\s,j,p)\in\CC\times\NNN\times \ZZZ^{d}$,
\begin{equation}\label{eq:4.12bis}
\begin{aligned}
E_{q}&=\left((E_{q})_{k}^{k'}\right)_{k,k'\in\CC\times\NNN\times \ZZZ^{d}}=
\left((E_{q})_{\s,j}^{\s',j'}(p-p')\right)_{k,k'\in\CC\times\NNN\times \ZZZ^{d}}, \quad q=0,1,\\
&(E_{1})_{\s,j}^{\s,j'}(p-p')\equiv0, \quad \forall\; j,j'\in\NNN, \;\; p,p'\in\ZZZ^{d}.
\end{aligned}
\end{equation}
\end{defi}

 Note that the operator
$\calL_{4}$ has the form (\ref{eq:4.11}) and satisfies the (\ref{eq:4.12}) and (\ref{eq:4.12bis}) as well as the estimates
(\ref{eq:3.2.6a}) and (\ref{eq:3.2.6b}). 
Note that each component
$(E_{q})_{\s}^{\s'}$, $q=0,1$, represent the matrix of the multiplication operator by a function.
This fact is not necessary for our analysis, and it cannot be preserved during the algorithm.

Define
\begin{equation}\label{eq:4.13}
N_{-1}:=1, \quad N_{\nu}:=N_{\nu-1}^{\chi}=N_{0}^{\chi^{\n}}, \;\;\forall\;\n\geq0, \;\;\chi=\frac{3}{2}.
\end{equation}
and
\begin{equation}\label{eq:4.14}
\al=7\tau+3, \qquad \h_{2}:=\h_{1}+\be,
\end{equation}
where $\h_{1}$ is defined in (\ref{eq:3.2.0}) and $\be=7\tau+5$.
 Consider
 $ \calL_4=\calL_0= \oo\cdot\del_{\f}\id+\DD_0+\RR_0$ with $\RR_0= E_{1}^{0}D+ E_{0}^{0}$,
we define
\begin{equation}\label{piccolo1}
\de_{s}^{0}:=|E_{1}^{0}|_{s,\g}+|E_{0}^{0}|_{s,\g}, \quad {\rm for} \quad s\geq0.
\end{equation}
%

\begin{lemma}[{\bf KAM iteration}]\label{teo:KAM}
Let $q>\h_{1}+\gots_{0}+\be$. There exist constant $C_{0}>0$, $N_{0}\in\NNN$ large, such that
if
\begin{equation}\label{eq:4.15}
N_{0}^{C_{0}}\g^{-1}\de_{\gots_{0}+\be}^{0}\leq1,
\end{equation}
then, for any $\nu\geq0$, one has:

\noindent
$({\bf S1.})_{\nu}$ Set
$\Lambda^{\g}_{0}:=\Lambda_{o}$ and for $\nu\geq1$ 
\begin{equation}\label{eq:419bis}
\begin{aligned}
&\Lambda_{\nu}^{\g}:=
\left\{\la\in\Lambda_{\nu-1}^{\g} : 
|\oo\cdot\ell\!+\! \mu_{h}^{\nu-1}(\la)\!-\! \mu_{h'}^{\nu-1}(\la)|\geq
\frac{\g|\s j^{2}-\s'j'^{2}|}{\langle\ell\rangle^{\tau}},
\forall |\ell|\leq N_{\nu-1}, \! h,h'\in\CC\times\NNN
\right\},
\end{aligned}
\end{equation}
For any $\la\in\Lambda_{\nu}^{\g}:=\Lambda_{\nu}^{\g}({\bf u})$, there exists an invertible map $\Phi_{\n-1}$ of the form $\Phi_{-1}=\id$ and for $\nu\geq1$,  $\Phi_{\nu-1}:=\id+\Psi_{\nu-1}: {\bf H}^{s}\to{\bf H}^{s}$, with the following properties.

 The maps $\Phi_{\nu-1}$, $\Phi_{\nu-1}^{-1}$
are reversibility-preserving according to Definition \ref{reversible},
moreover $\Psi_{\nu-1}$ is T\"oplitz in time, $\Psi_{\nu-1}:=\Psi_{\nu-1}(\f)$ (see (\ref{eq:2.16})) and satisfies the bounds:
\begin{equation}\label{eq:4.22}
\begin{aligned}
|\Psi_{\nu-1}|_{s,\g}\leq 
\de_{s+\be}^{0} N_{\nu-1}^{2\tau+1}N_{\nu-2}^{-\al},
\end{aligned}
\end{equation}
Setting, for $\nu\geq 1$, $\calL_{\nu}:=\Phi_{\nu-1}^{-1}\calL_{\nu-1}\Phi_{\nu-1}$, we have:
\begin{equation}\label{eq:4.16}
\begin{aligned}
\calL_{\nu}&=\oo\cdot\del_{\f} \id +\DD_{\nu}+\RR_{\nu}, \qquad
\DD_{\nu}={\rm diag}_{h\in\CC\times\NNN}\{\mu_{h}^{\nu}\},
\\
\quad
\mu_{h}^{\nu}(\la)&=\mu_{\s,j}^{\nu}=\mu_{\s,j}^{0}(\la)+r_{\s,j}^{\nu}(\la), \quad 
\mu_{\s,j}^{0}(0)=-\s i\,m(\la,{\bf u}(\la))j^{2},
\end{aligned}
\end{equation}
and
\begin{equation}\label{eq:4.18}
\RR_{\nu}=E_{1}^{\nu}(\la)D+E_{0}^{\nu}(\la),
\end{equation}
where $\RR_\nu$ is reversible and the matrices
$E_{q}^{\nu}$ satisfy (\ref{eq:4.12bis}) for $q=1,2$.
For $\nu\geq0$ one has $r_{h}^{\nu}\in i \RRR$,  $r_{\s,j}^{\nu}=-r_{-\s,j}^{\nu}$ and the following bounds hold:
\begin{equation}\label{eq:4.20}
|r_{h}^{\nu}|_{\g}:=|r_{h}^{\nu}|_{\Lambda_{\nu}^{\g},\g}\leq\e C.
\end{equation}
Finally, if we define
\begin{equation}\label{piccolo2}
\de_{s}^{\nu}:=|E_{1}^{\nu}|_{s,\g}+|E_{0}^{\nu}|_{s,\g},
\quad \forall s\geq0,
\end{equation}
one has $\forall\;s\in[\gots_{0},q-\h_{1}-\be]$ ($\al$ is defined in \eqref{eq:4.14}) and $\nu\geq0$
\begin{equation}\label{eq:4.21}
\begin{aligned}
\de_{s}^{\nu}&\leq \de_{s+\be}^{0}N_{\nu-1}^{-\al},\\
\de_{s+\be}^{\nu}&\leq \de_{s+\be}^{0}N_{\nu-1}.
\end{aligned}
\end{equation}

\noindent
$({\bf S2})_{\nu}$ For all $j\in\NNN$ there exists  Lipschitz extensions
$\tilde{\mu}_{h}^{\nu}(\cdot) : \Lambda \to \RRR$ of $\mu_{h}^{\nu}(\cdot):\Lambda_{\nu}^{\g}\to \RRR$,
such that for $\nu\geq1$,
\begin{equation}\label{eq:4.23}
|\tilde{\mu}_{h}^{\nu}-\tilde{\mu}_{h}^{\nu-1}|_{\g}\leq
\de_{\gots_{0}}^{\nu-1}, \qquad \forall\; k\in\CC\times\NNN.
\end{equation}

\noindent
$({\bf S3})_{\nu}$ Let ${\bf u}_{1}(\la)$, ${\bf u}_{2}(\la)$ be Lipschitz families of Sobolev functions, defined for $\la\in\Lambda_{o}$ such that (\ref{eq:4.2}), (\ref{eq:4.15}) hold with
$\RR_{0}=\RR_{0}({\bf u}_{i})$ with $i=1,2$.
Then for $\nu\geq0$, for any $\la\in\Lambda_{\nu}^{\g_{1}}\cap\Lambda_{\nu}^{\g_{2}}$,
with $\g_{1},\g_{2}\in[\g/2,2\g]$, one has
\begin{subequations}\label{eq:4.24}
\begin{align}
|E_{1}^{\nu}({\bf u}_{1})-E_{1}^{\nu}({\bf u}_{2})|_{\gots_{0}}+|E_{0}^{\nu}({\bf u}_{1})-E_{0}^{\nu}({\bf u}_{2})|_{\gots_{0}}&\leq
\e N_{\nu-1}^{-\al}||{\bf u}_{1}-{\bf u}_{2}||_{\gots_{0}+\h_{2}},
\label{eq:4.24a}\\
|E_{1}^{\nu}({\bf u}_{1})-E_{1}^{\nu}({\bf u}_{2})|_{\gots_{0}+\be}+|E_{0}^{\nu}({\bf u}_{1})-E_{0}^{\nu}({\bf u}_{2})|_{\gots_{0}+\be}&\leq
\e N_{\nu-1}||{\bf u}_{1}-{\bf u}_{2}||_{\gots_{0}+\h_{2}},\label{eq:4.24b}
\end{align}
\end{subequations}
and moreover, for $\nu\geq1$, for any $s\in[\gots_{0},\gots_{0}+\be]$, for any $k\in\CCC\times\NNN$
and for any $\la\in\Lambda_{\nu}^{\g_{1}}\cap\Lambda_{\nu}^{\g_{2}}$,
\begin{subequations}\label{eq:4.24bis}
\begin{align}
|(r_{h}^{\nu}({\bf u}_{2})-r_{h}^{\nu}({\bf u}_{1}))-(r_{h}^{\nu-1}({\bf u}_{2})-r_{h}^{\nu-1}({\bf u}_{1}))|
&\leq|E_{0}^{\nu-1}({\bf u}_{1})-E_{0}^{\nu-1}({\bf u}_{2})|_{\gots_{0}},\\
|(r_{h}^{\nu}({\bf u}_{2})-r_{h}^{\nu}({\bf u}_{1}))|\leq \e C||{\bf u}_{1}-{\bf u}_{2}||_{\gots_{0}+\h_{2}}.\label{aaa}
\end{align}
\end{subequations}

\noindent
$({\bf S4})_{\nu}$ Let $u_{1},u_{2}$ be as in $({\bf S3})_{\nu}$ and $0<\rho<\g/2$.
For any $\nu\geq0$ one has
\begin{equation}\label{eq:4.25}
\e CN_{\nu-1}^{\tau}\sup_{\la\in\Lambda_{o}}||{\bf u}_{1}-{\bf u}_{2}||_{\gots_{0}+\h_{2}} \leq \rho
\quad \Rightarrow \quad 
\Lambda_{\nu}^{\g}({\bf u}_{1})\subset \Lambda_{\nu}^{\g-\rho}({\bf u}_{2}),
\end{equation}

\end{lemma}

\proof
We start by proving that ${\bf (Si)_{0}}$ hold for $i=0,\ldots,4$.

${\bf (S1)_{0}}$. Clearly the properties (\ref{eq:4.20})-(\ref{eq:4.21}) hold by (\ref{eq:4.11}), (\ref{eq:4.12})
and the form of $\mu_{k}^{0}$ in (\ref{eq:4.16}), recall that $r_{k}^{0}=0$ .
Moreover, $m$ real implies that $\mu_{k}^{0}$ are imaginary. In addition to this,  our hypotheses guarantee that $\RR_{0}=E_{1}^{0}\del_{x}+E_{0}^{0}$ and
$\calL_{0}$ are reversible operators.

${\bf (S2)_{0}}$. We have to extend the eigenvalues $\mu_{k}^{0}$ from
the set $\Lambda_{0}^{\g}$ to the entire $\Lambda $. Namely we extend
the funtion $m(\la)$ to a $\tilde{m}(\la)$ that is lipschitz in $\Lambda$, with the same sup norm
and Lipschitz semi-norm, 
by Kirszbraun theorem.

${\bf (S3)_{0}}$. It holds by (\ref{eq:3.2.6b}) for $\gots_{0}$, $\gots_{0}+\be$ using 
(\ref{eq:4.2}) and (\ref{eq:4.14}).

${\bf (S4)_{0}}$. By definition one has 
$\Lambda_{0}^{\g}({\bf u}_{1})=\Lambda_{o}=\Lambda_{0}^{\g-\rho}({\bf u}_{2})$, then the  
(\ref{eq:4.25}) follows trivially.

\subsubsection{Kam step}

In this Section we  show in detail one step of the KAM iteration. In other words we will
show how to define the transformation $\Phi_{\nu}$ and $\Psi_{\n}$ that trasform the operator
$\calL_{\n}$ in the operator $\calL_{\nu+1}$. For simplicity we shall avoid to write the index,
but we will only write $+$ instead of $\n+1$. 

We consider a transformation of the form $\Phi=\id+\Psi$, with $\Psi:=(\Psi_{\s}^{\s'})_{\s,\s'=\pm1}$, acting on the operator 
$$\calL= \oo\cdot\del_{\f} \id+\DD+\RR
$$
with $\DD$ and $\RR$  as in \eqref{eq:4.16}, \eqref{eq:4.18}.
Then, $\forall\; {\bf h}\in {\bf H}^{s}$, one has
\begin{equation}\label{eq:4.1.22}
\begin{aligned}
\calL\Phi{\bf h}&=\oo\cdot\del_{\f}(\Phi({\bf h}))+\DD\Phi{\bf h}+\RR\Phi{\bf h}\\
&=\Phi\left(\oo\cdot\del_{\f}{\bf h}+\DD{\bf h}\right)+\left(\oo\cdot\del_{\f}\Psi
+[\DD,\Psi]+\Pi_{N}\RR\right){\bf h}+\left(\Pi_{N}^{\perp}\RR+\RR\Psi\right){\bf h},
\end{aligned}
\end{equation}
where $[\DD,\Phi]:=\DD\Phi-\Phi\DD$, and $\Pi_{N}$ is defined in (\ref{smoothop}).
The smoothing operator $\Pi_{N}$ is necessary for technical reasons: it will be used in order to obtain
suitable estimates on the high norms of the transformation $\Phi$.

In the following Lemma we will show how to solve the \emph{homological} equation
\begin{equation}\label{homeq}
\oo\cdot\del_{\f}\Psi+[\DD,\Psi]+\Pi_{N}\RR=[\RR], \quad {\rm where} \quad
[\RR]_{k}^{k'}:=
\left\{
\begin{aligned}
& (E_{0})_{k}^{k}=(E_{0})_{\s,j}^{\s,j}(0), \quad k=k',\\
& 0 \quad \quad k \neq k',
\end{aligned}\right.
\end{equation}
for $k,k'\in \CC\times\NNN\times\ZZZ^{d}$.

\begin{lemma}[{\bf Homological equation}]\label{lemhom}
For any $\la\in\Lambda^{\g}_{\nu+1}$ there exists a unique solution $\Psi=\Psi(\f)$ of the homological equation
(\ref{homeq}), such that
\begin{equation}\label{eq:4.1.33}
|\Psi|_{s,\g}\leq C N^{2\tau+1}\g^{-1}
\de_{s}
\end{equation}
Moreover, for $\g/2\leq \g_{1},\g_{2}\leq 2\g$, and if $u_{1}(\la), u_{2}(\la)$ are Lipschitz functions,
then $\forall \; s\in[\gots_{0},\gots_{0}+\be]$, 
$\la\in \Lambda^{\g_{1}}_{+}(u_{1})\cap\Lambda^{\g_{2}}_{+}(u_{2})$, one has
\begin{equation}\label{eq:4.1.44}
|\Delta_{12}\Psi|_{s}\leq C N^{2\tau+1}\g^{-1}\left((|E_{1}({\bf u}_{2})|_{s}+|E_{0}({\bf u}_{2})|_{s})
||{\bf u}_{1}-{\bf u}_{2}||_{\gots_{0}+\h_{2}}+
|\Delta_{12}E_{1}|_{s}+|\Delta_{12}E_{0}|_{s}\right),
\end{equation}
where we define $\Delta_{12}\Psi=\Psi({\bf u}_{1})-\Psi({\bf u}_{2})$.\\
Finally, one has $\Psi : {\bf X}^{s}\to {\bf X}^{s}$, i.e. the operator $\Psi$ is reversible preserving.
\end{lemma}

\prova
On each component $k=(\s,j,p),k'=(\s',j',p')\in C\times\NNN\times\ZZZ^{d}$, the equation (\ref{homeq}) reads
\begin{equation}\label{eq:4.1.55}
i\oo\cdot(p-p')\Psi_{k}^{k'}+\DD_{k}^{k}\Psi_{k}^{k'}-\Psi_{k}^{k'}\DD_{k'}^{k'}+\RR_{k}^{k'}=[\RR]_{k}^{k'}.
\end{equation}
then, by defining
\begin{equation}\label{29}
\begin{aligned}
d_k^{k'}:=i\oo\cdot(p-p')+\mu_{k}-\mu_{k'}
&\stackrel{(\ref{eq:4.16})}{=}i\oo\cdot (p-p') -im(\s j^{2}-\s' j'^{2})
+r_{\s,j}-r_{\s',j'}
\end{aligned}
\end{equation}
we get 
\begin{equation}\label{eq:4.1.77}
\Psi_{k}^{k'}=
\frac{-\RR_{k}^{k'}}{d_{k}^{k'}}, \;\;\;\;\; k \neq k', \;\; |p-p'|\leq N, \\
\end{equation}
and $\Psi_{k}^{k'}\equiv0$ otherwise.
Clearly the  solution has the form $\Psi_{\s,j,p}^{\s',j',p'}=\Psi_{\s,j}^{\s',j'}(p-p')$ and hence we can define a time-dependent change of variables as $\Psi_{\s,j}^{\s',j'}(\f))  =\sum_{\ell\in\ZZZ^{d}}\Psi_{\s,j}^{\s',j'}(\ell)e^{i\ell\cdot\f}$.

Note that, by (\ref{eq:419bis}) and (\ref{dio}) 
one has for all $k\neq k'\in C\times\NNN\times\ZZZ^{d},$ setting $k=(\s,j,p)$, $k'=(\s',j',p')$ and $\ell=p-p'$
\begin{equation}\label{eq:4.1.99}
|d^{k'}_{k}|\geq \left\{\begin{aligned} 
& \frac{\g(j^{2}+j'^{2})}{\langle\ell\rangle^{\tau}}, \quad \s= -\s' ,\\
&  \frac{\g(j+j')}{\langle\ell\rangle^{\tau}}, \quad {\rm if}\;\;\; \s=\s'\;  j\neq j',\\
& \frac{\g}{\langle\ell\rangle^{\tau}}, \quad {\rm if} \;\;\;\s=\s'\;  j=j' \; p\neq p'
\end{aligned}\right.
\end{equation}
This implies that,  for  $\s\neq\s'$, we have
\begin{equation}\label{eq:4.1.112}
|\Psi_{k}^{k'}|\leq \g^{-1}|\ell|^{\tau}\left(|(E_{1})_{k}^{k'}|+|(E_{0})_{k}^{k'}|
\right)\frac{j}{j^{2}+j'^{2}},
\end{equation}
while, for  $\s=\s'$,
\begin{equation}\label{eq:4.1.112bis}
|\Psi_{k}^{k'}|\leq \left\{
\begin{aligned}
& \g^{-1}\langle\ell\rangle^{\tau}|(E_{0})_{k}^{k'}|\frac{1}{j+j'}, \quad j\neq j',\\
&\g^{-1}\langle\ell\rangle^{\tau}|(E_{0})_{k}^{k'}|, \quad j=j',
 \end{aligned}\right.
\end{equation}
and 
we can estimate the divisors $\de_{k}^{k'}$ from below, hence, 
by the definition of the $s-$norm  in (\ref{decay}) in any case we obtain the estimate
\begin{equation}\label{eq:4.1.113}
|\Psi|_{s}\leq \g^{-1} N^{\tau}\de_{s}.
\end{equation}
If we define the operator $A$ as
\begin{equation}\label{eqrmk1}
A_{k}^{k'}=A_{\s,j}^{\s',j'}(\ell):=\left\{
\begin{aligned}
&\Psi_{\s,j}^{\s,j}(\ell), \quad (\s,j)=(\s',j')\in \CC\times\NNN, \;\; \ell\in\ZZZ^{d},\\
&0, \quad {\rm otherwise},
\end{aligned}\right.
\end{equation}
we have proved  the following Lemma
\begin{lemma}\label{rmk200}
The operator $\Psi-A$ is regularizing, indeed,
\begin{equation}\label{34}
\begin{aligned}
|D(\Psi-A)|_{s}^{2}&:=\sup_{\s,\s'\in \CC}\sum_{\substack{k\in\NNN,\\ \ell\in\ZZZ^{d}}}
\sup_{\substack{j-j'=k \\ j\neq j'}}|\Psi_{\s,j}^{\s',j'}(\ell)|^{2}\langle\ell,k\rangle^{2s},
\stackrel{(\ref{eq:419bis})}{\leq_{s}} 
\g^{-2} N^{2\tau} \de_{s},
\end{aligned}
\end{equation}
where $D$ is defined in (\ref{eq:4.12bis}). 
\end{lemma}
This Lemma will be used in the study of the remainder of the conjugate operator. 
In particular we will use it to prove that the reminder is still in the class of operators described in $(\ref{eq:4.12})$.

Now we need a bound on the lipschitz semi-norm of the transformation. Then, given 
$\la_{1},\la_{2}\in\Lambda^{\g}_{\nu+1}$, one has, for $k=(\s,j,p),k'=(\s',j',p')\in\CC\times\NNN\times\ZZZ^{d}$, and
$\ell:=p-p'$,
\begin{equation}\label{fiorellino2}
|\Psi_{k}^{k'}(\la_{1})-\Psi_{k}^{k'}(\la_{2})|\leq
\frac{|\RR_{k}^{k'}(\la_{1})-\RR_{k}^{k'}(\la_{2})|}{|d_{k}^{k'}(\la_{1})|}+
|\RR_{k}^{k'}(\la_{2})|
\frac{|d_{k}^{k'}(\la_{1})-
d_{k}^{k'}(\la_{2})|}{|d_{k}^{k'}(\la_{1})||d_{k}^{k'}(\la_{2})|},
\end{equation}
Now, recall that $\oo=\la\bar{\oo}$,  by using that 
$
\g|m|^{\rm lip}=\g|m-1|^{\rm lip}\stackrel{(\ref{eq:3.2.4a})}{\leq}\e C,
$
and by (\ref{eq:4.20}), we obtain
\begin{equation}\label{eq:4.1.115}
\begin{aligned}
|d_{k}^{k'}(\la_{1})-d_{k}^{k'}(\la_{2})|&\stackrel{(\ref{29}), (\ref{eq:4.16})}{\lessdot}
|\la_{1}-\la_{2}|\cdot\left(|\ell|+\e\g^{-1}|\s j^{2}-\s'j'^{2}|+\e\g^{-1}
\right).
\end{aligned}
\end{equation}
Then, for $\s,\s'=\pm1$, $j\neq j'$ and $\e\g^{-1}\leq1$,
\begin{equation}\label{eq:4.1.116}
\frac{|d_{k}^{k'}(\la_{1})-
d_{k}^{k'}(\la_{2})|}{|d_{k}^{k'}(\la_{1})||d_{k}^{k'}(\la_{2})|}
\stackrel{(\ref{eq:4.1.115}),(\ref{eq:419bis})}{\lessdot}
|\la_{1}-\la_{2}|\frac{N^{2\tau+1}\g^{-2}}{|\s j^{2}-\s' j'^{2}|}
\end{equation}
for $|\ell|\leq N$. By using that, for any $|\ell|\leq N$, $j,j'\geq1$, $j\neq j'$, $\s,\s'=\pm1$
$$
|\RR_{k}^{k'}|/|\s j^{2}-\s' j'^{2}|\leq |(E_{1})_{k}^{k'}|+|(E_{0})_{k}^{k'}|
$$
for any $\la\in\Lambda^{\g}_{\nu+1}$, the (\ref{eq:4.1.116}), the fact that 
 $|d_{k}^{k'}|\geq\g/\langle\ell\rangle^{-\tau}$ for $\s=\s'$ and $j=j'$, one has
 and finally the (\ref{eq:4.1.113}), we get
\begin{equation}\label{eq:4.1.118}
|\Psi|_{s,\g}:=|\Psi|^{{\rm sup}}_{s}+\g \sup_{\la_{1}\neq \la_{2}}
\frac{|\Psi(\la_{1})-\Psi(\la_{2})|_{s}}{|\la_{1}-\la_{2}|}\leq
\g^{-1} C N^{2\tau+1}\left(
|E_{1}|_{s,\g}+|E_{0}|_{s,\g}\right),
\end{equation}
and hence the (\ref{eq:4.1.33}) is proved.

Let us check the (\ref{eq:4.1.44}). For $\la\in\Lambda^{\g_1}_{\nu+1}\cap\Lambda^{\g_2}_{\nu+1}$, 
if $k=(\s,j,p)\neq(\s',j',p')=k'$, one has
\begin{equation}\label{eq:4.1.119}
\begin{aligned}
 |\Delta_{12}\Psi_{k}^{k'}|\leq 
\frac{|\Delta_{12}\RR_{k}^{k'}|}{|d_{k}^{k'}({\bf u}_1)|}
+&|\RR_{k}^{k'}({\bf u}_2)|
\frac{|\Delta_{12}d_{k}^{k'}|}{|d_{k}^{k'}({\bf u}_1)||d_{k}^{k'}({\bf u}_{2})|}
\stackrel{(\ref{eq:3.2.4b}), (\ref{aaa})}{\leq}
N^{2\tau}\g^{-1}\left(|(\Delta_{12}E_{1})_{k}^{k'}|+|(\Delta_{12}E_{0})_{k}^{k'}|\right.\\
&\left.+\left(|(E_{1})_{k}^{k'}({\bf u}_{2})|+
|(E_{0})_{k}^{k'}({\bf u_{2}})|\right)||{\bf u}_{1}-{\bf u}_{2}||_{\gots_{0}+\eta_{2}}
\right)
\end{aligned}
\end{equation}
where we used $\e\g^{-1}\leq1$ $\g_{1}^{-1},\g_{2}^{-1}
\leq\g^{-1}$, hence (\ref{eq:4.1.119})
implies the (\ref{eq:4.1.44}).

Since $\ol{\mu_{\s,j}}=-\mu_{\s,j}$ and the operator $\RR$ is reversible (see (\ref{rmkrev7})), 
by (\ref{eq:4.1.77}), we have that
\begin{equation}\label{eq:4.1.122}
\ol{\Psi_{k}^{k'}}=\ol{\Psi_{\s,j}^{\s',j'}(\ell)}=
\frac{- \ol{\RR_{\s,j}^{\s',j'}(\ell)}}{-i\oo\cdot\ell+\ol{\mu_{\s',j'}}- \ol{\mu_{\s, j}}}
=\frac{\RR_{\s,j}^{\s',j'}(\ell)}{-(i\oo\cdot\ell+\mu_{\s',j'}-\mu_{\s,j} )}=
\Psi_{\s,j}^{\s',j'}(\ell)=\Psi_{k}^{k'},
\end{equation}
so that, by Lemma \ref{lemFou}, for any $\s,\s'=\pm1$, the operators $\Psi_{\s}^{\s'}$ are 
reversibility preserving. In the same way,  
again thanks to the reversibility of $\RR$, one can check 
$\ol{\Psi_{\s,j}^{\s',j'}(-\ell)}=\Psi_{-\s,j}^{-\s',j}(\ell)$
which implies $\Psi: {\bf X}^{s}\to{\bf X}^{s}$, i.e. $\Psi$ is reversibility
preserving.
\EP

By Lemma \ref{inverse}, for $\de_{\gots_{0}}$ small enough, we have by (\ref{eq:4.1.33}) for $s=\gots_{0}$
\begin{equation}\label{40}
C(\gots_{0})|\Psi|_{\gots_{0}}\leq\frac{1}{2},
\end{equation}
then, the operator $\Phi=\calI+\Psi$ si invertible. In this case we can conjugate the operator $\calL$ to an operator
$\calL_{+}$  as shown in the next Lemma.

\begin{lemma}[{\bf The new operator $\calL_{+}$}]\label{lemma2}
Consider the operator $\Phi=\calI+\Psi$ with $\Psi$ definend in Lemma \ref{lemhom}. Then, one has
\begin{equation}\label{41}
\calL_{+}:=\Phi^{-1}\calL\Phi:=\oo\cdot\del_{\f}\id+\DD_{+}+\RR_{+},
\end{equation}
where the diagonal operator $\DD_{+}$ has the form
\begin{equation}\label{42}
\begin{aligned}
\DD_{+}&:=diag_{h\in C\times \NNN}\{\mu_{h}^{+}\}, \\
\mu_{h}^{+}&:=\mu_{h}+(E_{0})_{h}^{h}(0)=\mu_{h}^{0}+
r_{h}+(E_{0})_{h}^{h}=:\mu_{h}^{0}+r_{h}^{+},
\end{aligned}
\end{equation}
with $h:=(\s,j)\in \CC\times\NNN$ and the remainder,
\begin{equation}\label{43}
\RR_{+}:=E_{1}^{+}D+E_{0}^{+}
\end{equation}
where  $E_{i}^{+}$ are linear bounded operators of the form (\ref{eq:4.12bis}) for $i=0,1$.

Morever, the eigenvalues $\mu_{h}^{+}$ satisfiy
\begin{equation}\label{45}
|\mu_{h}^{+}-\mu_{h}|^{{\rm lip}}=|r_{h}^{+}-r_{h}|^{{\rm lip}}=|(E_{0})_{h}^{h}(0)|^{{\rm lip}}\leq
|E_{0}|^{{\rm lip}}_{\gots_{0}}, \quad h\in C\times\NNN,
\end{equation}
while the remainder $\RR_{+}$ satisfies
\begin{equation}\label{46}
\begin{aligned}
&\de_{s}^{+}:=|E_{1}^{+}|_{s,\g}+|E_{0}^{+}|_{s,\g}\leq_{s} N^{-\be}\de_{s+\be}+N^{2\tau+1}\g^{-1}
\de_{s}\de_{\gots_{0}},\\ &\de_{s+\be}^{+}\leq_{s+\be} \de_{s+\be}+N^{2\tau+1}\g^{-1}
\de_{s+\be}\de_{\gots_{0}}.
\end{aligned}
\end{equation}

Finally,  for $\g/2\leq \g_{1},\g_{2}\leq 2\g$, and if $u_{1}(\la), u_{2}(\la)$ are Lipschitz functions,
then $\forall \; s\in[\gots_{0},\gots_{0}+\be]$, 
$\la\in \Lambda^{\g_{1}}_{+}(u_{1})\cap\Lambda^{\g_{2}}_{+}(u_{2})$,
setting 
$|\Delta_{12}E_{1}|_{s}+|\Delta_{12}E_{0}|_{s}=\Delta_{s}$,
we have:
\begin{equation}\label{eq:4.2.19}
\begin{aligned}
\Delta^+_{s}&\leq
|\Pi_{N_{ }}^{\perp}\Delta_{12}E_{0}^{ }|_{s}+|\Pi_{N_{ }}^{\perp}\Delta_{12}E_{1}^{ }|_{s}\\
&+N_{ }^{2\tau+1}\g^{-1}
(\de^{ }_{s}({\bf u}_{1})+\de^{ }_{s}({\bf u}_{2}))(\de_{\gots_{0}}^{ }({\bf u}_{1})+
\de_{\gots_{0}}^{ }({\bf u}_{2}))||{\bf u}_{1}-{\bf u}_{2}||_{s+\h_{2}}\\
&+
N_{ }^{2\tau+1}\g^{-1}(\de_{s}^{ }({\bf u}_{1})+
\de_{s}^{ }({\bf u}_{2}))\Delta_{\gots_{0}}^{ } +N_{ }^{2\tau+1}\g^{-1}(\de_{\gots_{0}}^{ }({\bf u}_{1})
+
\de_{\gots_{0}}^{ }({\bf u}_{2}))\Delta_{s}^{ },
\end{aligned}
\end{equation}

\end{lemma}

\prova The expression \eqref{42} follows by \eqref{homeq}, the bound (\ref{45}) follows by (\ref{decay2}). 

The bound (\ref{46}) is more complicated. First of all we note that, by (\ref{eq:4.1.22}) and (\ref{homeq}), we have
\begin{equation}\label{48}
\begin{aligned}
\RR_{+}:=\Phi^{-1}\left(\Pi_{N}^{\perp}\RR+\RR \Psi-\Psi[\RR]
\right):=E_{1}^{+}D+E_{0}^{+}, 
\end{aligned}
\end{equation}
where
\begin{equation}\label{49}
\begin{aligned}
E_{1}^{+}&:=\Phi^{-1}\left(\Pi_{N}^{\perp}E_{1}+ E_{1}A\right),\\
E_{0}^{+}&:=\Phi^{-1}\left(\Pi_{N}^{\perp}E_{0}+E_{0}\Psi-\Psi[\RR]+E_{1}D(\Psi-A)\right),
\end{aligned}
\end{equation}
where $A$ is defined in (\ref{eqrmk1}).

We can estimate the first of the (\ref{49}) by
\begin{equation}\label{51}
\begin{aligned}
|E_{1}^{+}|_{s,\g}&\stackrel{(\ref{eq:2.11a}),(\ref{eq:2.14})}{\leq_{s}}
2|\Pi_{N}^{\perp}E_{1}|_{s,\g}+(1+|\Psi|_{s,\g})\left(|\Pi_{N}^{\perp}E_{1}|_{\gots_{0},\g}+
|E_{1}|_{\gots_{0},\g}|A|_{\gots_{0},\g}
\right)+2(|E_{1}|_{s,\g}|A|_{\gots_{0},\g}+|E_{1}|_{\gots_{0},\g}|A|_{s,\g})\\
&\stackrel{(\ref{eq:4.1.33}),(\ref{eq:2.22})}{\leq_{s}}
N^{-\be}|E_{1}|_{s+\be,\g}+N^{2\tau+1}\g^{-1}\de_{\gots_{0}}\de_{s}.
\end{aligned}
\end{equation}

The  bound  on $E_0^+$ is obtained in the same way by recalling that, by Lemma \ref{rmk200},
\begin{equation}\label{52}
|D(\Psi-A)|_{s,\g}\leq \g^{-1}N^{2\tau+1}\de_{s}.
\end{equation}
The second bound in (\ref{46}) follows exactly in the same way.

Now, consider $\Delta_{12}E_{1}^{ }+\Delta_{12}E_{0}^{ }$, that is defined for
$\la\in\Lambda_{ }^{\g_{1}}({\bf u}_{1})\cap\Lambda_{ }^{\g_{2}}({\bf u}_{2})$. 
Define also $E_{1,i}^{ }:=E_{1}^{ }({\bf u}_{i})$ and $E_{0,i}^{ }:=E_{0}^{ }({\bf u}_{i})$, for $i=1,2$. We prove the bounds only for $E_{0}^{ +1}$, which is the hardest case, the bounds on $E_{1}^{ +}$ follow in the same way.
By Lemma \ref{ bubbole} and the definition of $E_{0}^{ +}$ (see \eqref{49})
one has
\begin{equation}\label{eq:4.2.18}
\begin{aligned}
|\Delta_{12}E_{0}^{ +}|_{s}
&\stackrel{(\ref{eq:4.1.33}),(\ref{eq:4.1.44})}{\leq_{s}}
|\Pi_{N_{ }}^{\perp}\Delta_{12}E_{0}^{ }|_{s}+N_{ }^{2\tau+1}\g^{-1}
(\de^{ }_{\gots_{0}}({\bf u}_{1})+\de^{ }_{\gots_{0}}({\bf u}_{2}))(\de_{s}^{ }({\bf u}_{1})+
\de_{s}^{ }({\bf u}_{2}))||{\bf u}_{1}-{\bf u}_{2}||_{s+\h_{2}}\\
&+N_{ }^{2\tau+1}\g^{-1}(\de_{s}^{ }({\bf u}_{1})+
\de_{s}^{ }({\bf u}_{2}))|\Delta_{12}E_{0}^{ }|_{\gots_{0}}
+N_{ }^{2\tau+1}\g^{-1}(\de_{\gots_{0}}^{ }({\bf u}_{1})+
\de_{\gots_{0}}^{ }({\bf u}_{2}))|\Delta_{12}E_{0}^{ }|_{s},
\end{aligned}
\end{equation}
We prove equivalent bounds for $E_1^{ +}$; then, we obtain \eqref{eq:4.2.19}
using the bounds given in Lemmata \ref{lemhom} and \ref{inverse}.
to estimate the norms of the transformation $\Phi$ .

\EP

In the next Section we will show that it is possible to iterate the procedure described above infinitely many times.

\subsubsection{The iterative Scheme}\label{iteration}

Here we complete the proof of the Lemma \ref{teo:KAM} by induction on $\nu\geq0$.
Hence, assume that $({\bf Si})_{\nu}$ hold. Then we prove $({\bf Si})_{\nu+1}$
for $i=1,2,3,4$. We will use the
estimates obtained in the previous Section.

\noindent
${\bf (S1)}_{\nu+1}$ The eigenvalues $\mu_{h}^{\nu}$ of $\DD_{\nu}$ are defined on $\Lambda^{\g}_{\nu}$,
then also the set $\Lambda_{\nu+1}^{\g}$ is well defined. Then, by Lemma \ref{lemhom},
for any $\la\in\Lambda_{\nu+1}^{\g}$ there exists a unique solution $\Psi_{\nu}$ of the equation
(\ref{homeq}), such that, by inductive hypothesis ${\bf (S1)_{\nu}}$,
\begin{equation}\label{58}
|\Psi_{\nu}|_{s,\g}\stackrel{(\ref{eq:4.1.33})}{\lessdot}
\g^{-1}N_{\nu}^{2\tau+1}\de_{s}^{\nu}\stackrel{(\ref{eq:4.21})}{\leq}
\g^{-1}N_{\nu}^{2\tau+1}N^{-\al}_{\nu-1}\de_{s+\be}^{0},
\end{equation}
hence the (\ref{eq:4.22}) holds at the step $\nu+1$. Moreover, by (\ref{58}) and hypothesis (\ref{eq:4.15}), one has
for $s=\gots_{0}$
\begin{equation}\label{59}
C(\gots_{0})|\Psi_{\nu}|_{\gots_{0},\g} \leq C(\gots_{0})\g^{-1}N_{\nu}^{2\tau+1}N_{\nu-1}^{-\al}
\de_{\gots_{0}+\be}^{0}\leq\frac{1}{2},
\end{equation}
for $N_{0}$ large enough and using, for $\nu=1$ the smallness condition (\ref{eq:4.15}). 
In this case, by Lemma \ref{inverse},
we have that the transformation $\Phi_{\nu}:=\calI+\Psi_{\nu}$ is invertible with
\begin{equation}\label{60}
|\Phi^{-1}_{\nu}|_{\gots_{0},\g}\leq2, 
\quad |\Phi^{-1}_{\nu}|_{s,\g}\leq1+C(s)|\Psi_{\nu}|_{s,\g}.
\end{equation}
Now, by Lemma \ref{lemma2}, we have $\calL_{\nu+1}:=\Phi^{-1}_{\nu}\calL_{\nu}\Phi_{\nu}=\oo\cdot\del_{\f}\id+\DD_{\nu+1}+\RR_{\nu+1}$, where
\begin{equation}\label{61}
\begin{aligned}
\DD_{\nu+1}&:=diag_{h\in C\times\NNN }\{\mu^{\nu+1}_{h}\},\quad
\mu_{h}^{\nu+1}:=\mu_{h}^{\nu}+(E_{0}^{\nu})_{h}^{h}(0)=\mu_{h}^{0}+r_{h}^{\nu+1},\\
\RR_{\nu+1}&=\Phi^{-1}_{\nu}\left(\Pi_{N_{\nu}}^{\perp}\RR_{\nu}+\RR_{\nu}\Psi_{\nu}-\Psi_{\nu}[\RR_{\nu}]\right)
=E_{1}^{\n+1}D+E_{0}^{\nu+1},
\end{aligned}
\end{equation}
where $E_{i}^{\nu+1}\rightsquigarrow E_{i}^+$, see (\ref{49}).
Let us check the (\ref{eq:4.21}) on the reminder $\RR_{\nu+1}$.  By (\ref{46}) in Lemma \ref{lemma2}, we have
\begin{equation}\label{64}
\begin{aligned}
\de_{s}^{\n+1}&\leq_{s} N_{\nu}^{-\be}\de_{s+\be}^{\nu}+\g^{-1}N_{\nu}^{2\tau+1}\de_{\gots_{0}}^{\nu}\de_{s}^{\n}\\
&\stackrel{(\ref{eq:4.21})}{\leq_{s}} N_{\nu}^{-\be}N_{\nu-1}\de_{s+\be}^{0}+
\g^{-1}N_{\n}^{2\tau+1}N_{\nu-1}^{-2\al}\de_{\gots_{0}+\be}^{0}\de_{s+\be}^{0}
\stackrel{(\ref{eq:4.14}),(\ref{27}),(\ref{eq:4.15})}{\leq_{s}}
\de_{s+\be}^{0}N_{\nu}^{-\al},
\end{aligned}
\end{equation}
that is the first of the (\ref{eq:4.21}) for $\n\rightsquigarrow \n+1$. In the last inequality we used that
$\chi=3\backslash2$, $\be>\al+1$ and $\chi(2\tau+1+\al)<2\al$, and this gives us a reason for the choices of $\be$
and $\al$ in  (\ref{eq:4.14}). Now,  by using the (\ref{46}) we have
\begin{equation}\label{65}
\begin{aligned}
\de_{s+\be}^{\nu+1}\leq_{s+\be} \de_{s+\be}^{\nu} + \g^{-1}N_{\nu}^{2\tau+1}\de_{\gots_{0}}\de_{s+\be}^{\nu}
\leq_{s+\be}\de_{s+\be}^{0}N_{\nu},
\end{aligned}
\end{equation}
for $N_{0}=N_{0}(s,\be)$ large enough. This completes the proof of the (\ref{eq:4.21}).

By using (\ref{45}) in Lemma \ref{lemma2}, we have, $\forall\; h\in C\times\NNN$,
\begin{equation}\label{67}
|\mu_{h}^{\n+1}-\mu_{h}^{\n}|_{\g}=
|r_{h}^{\n+1}-r_{h}^{\n}|_{\g} \leq\de_{\gots_{0}}^{\n}
\stackrel{(\ref{eq:4.21})}{\leq}\de_{\gots_{0}+\be}^{0}N_{\n-1}^{-\al},
\end{equation}
hence, we get the (\ref{eq:4.20}) by
$|r_{h}^{\n+1}|_{\g} \leq \sum_{i=0}^{\n}|r_{h}^{\n+1}-r_{h}^{\n}|_{\g}
\stackrel{(\ref{67})}{\leq}\de_{\gots_{0}+\be}^{0}K.$

Finally, we have to check that $\ol{\mu_{\s,j}^{\nu+1}}=-\mu_{\s,j}^{\nu+1}=\mu_{-\s,j}^{\nu+1}$. 
It follows by the inductive hypotheses since, by \eqref{rmkrev7}, one has $$\ol{(E^{\nu}_{0})_{\s,j}^{\s,j}(0)}=-{(E^{\nu}_{0})_{\s,j}^{\s,j}(0)}={(E^{\nu}_{0})_{-\s,j}^{-\s,j}(0)}$$

\noindent
${\bf (S2)_{\n+1}}$ Thanks to (\ref{67}), we can extend, by Kirszbraun theorem, the function
$\mu_{h}^{\n+1}-\mu_{h}^{\n}$ to a Lipschitz funtion on  $\Lambda$. Defining $\tilde{\mu}_{k}^{\nu+1}$ in this way, this extension has the same Lipschitz norm , so that the 
bound (\ref{eq:4.23}) hold. 

\noindent
${\bf (S3)_{\nu+1}}.$ Let 
$\la\in\Lambda^{\g_{1}}_{\nu}({\bf u}_{1})\cap\Lambda_{\n}^{\g_2}({\bf u}_{2})$,
then by Lemma \ref{lemhom} we can construct
operators $\Psi_{\nu}^{i}:=\Psi_{\nu}({\bf u}_{i})$ and $\Phi_{\nu}^{i}=\Phi_{\nu}({\bf u}_{i})$
for $i=1,2$. Using the (\ref{eq:4.1.44}),
we have that
\begin{equation}\label{eq:4.2.14}
\begin{aligned}
|\Delta_{12}\Psi_{\nu}|_{\gots_{0}}
&\stackrel{(\ref{eq:4.21}),(\ref{eq:4.24})}{\lessdot}
N_{\nu}^{2\tau+1}N_{\nu-1}^{-\al} \g^{-1}\left(\de^{0}_{\gots_{0}+\be}
+\e\right)||{\bf u}_{1}-{\bf u}_{2}||_{\gots_{0}+\h_{2}}\\
&\stackrel{(\ref{eq:4.15})}{\lessdot}
N_{\nu}^{2\tau+1}N_{\nu-1}^{-\al}||{\bf u}_{1}-{\bf u}_{2}||_{\gots_{0}+\h_{2}}\stackrel{(\ref{eq:4.14})}{\leq}
||{\bf u}_{1}-{\bf u}_{2}||_{\gots_{0}+\h_{2}},
\end{aligned}
\end{equation}
where we used  the fact that $\e \g^{-1}$ is small. Moreover, one can note that
\begin{equation}\label{eq:4.2.15}
|\Delta_{12}\Phi^{-1}_{\nu}|_{s}\stackrel{(\ref{eq:2.15}), (\ref{eq:4.2.14})}{\leq_{s}}
\left(|\Psi^{1}_{\nu}|_{s}+|\Psi_{\nu}^{2}|_{s}\right)||{\bf u}_{1}-{\bf u}_{2}||_{\gots_{0}+\h_{2}}
+|\Delta_{12}\Psi_{\nu}|_{s},
\end{equation}
then, by using the inductive hypothesis (\ref{eq:4.22}), the (\ref{eq:4.15}) and the (\ref{eq:4.2.15}) for $s=\gots_{0}$,
one obtain
\begin{equation}\label{eq:4.2.16}
|\Delta_{12}\Phi^{-1}_{\nu}|_{\gots_{0}}\lessdot ||{\bf u}_{1}-{\bf u}_{2}||_{\gots_{0}+\h_{2}}.
\end{equation}

The (\ref{eq:4.2.19}) with $s=\gots_{0}$, togeter with (\ref{eq:4.15}), (\ref{eq:4.21}) and (\ref{eq:4.24})
implies
\begin{equation}\label{eq:4.2.20}
\begin{aligned}
|\Delta_{12}E_{1}^{\nu+1}|_{\gots_0}+|\Delta_{12}E_{0}^{\nu+1}|_{\gots_{0}}&\leq_{\gots_{0}}
\left(\e N_{\nu-1}N_{\nu}^{-\be}+N_{\nu}^{2\tau+1}N_{\nu-1}^{-2\al}\e^{2}\g^{-1}
\right)||{\bf u}_{1}-{\bf u}_{2}||_{\gots_{0}+\h_{2}}\\
&\leq \e N_{\nu}^{-\al}||{\bf u}_{1}-{\bf u}_{2}||_{\gots_{0}+\h_{2}},
\end{aligned}
\end{equation}
for $N_{0}$ large enough and $\e\g^{-1}$ small. Moreover, consider the (\ref{eq:4.2.19})
with $s=\gots_{0}+\be$, then by (\ref{eq:4.15}),  (\ref{eq:4.24}) and (\ref{eq:4.21}), we obtain
for $N_{0}$ large enough
\begin{equation}\label{eq:4.2.21}
\begin{aligned}
|\Delta_{12}E_{1}^{\nu+1}|_{\gots_{0}+\be}+|\Delta_{12}E_{0}^{\nu+1}|_{\gots_{0}+\be}
&\leq_{\gots_{0}+\be}
(\de_{\gots_{0}+\be}^{\nu}({\bf u}_{1})+
\de_{\gots_{0}+\be}^{\nu}({\bf u}_{2}))||{\bf u}_{1}-{\bf u}_{2}||_{s+\h_{2}}+
|\Delta_{12}E_{1}^{\nu}|_{\gots_{0}+\be}+|\Delta_{12}E_{0}^{\nu}|_{\gots_{0}+\be}\\
&\leq C(\gots_{0}+\be)\e N_{\nu-1}||{\bf u}_{1}-{\bf u}_{2}||_{\gots_{0}+\h_{2}}\leq
\e N_{\nu}||{\bf u}_{1}-{\bf u}_{2}||_{\gots_{0}+\h_{2}}.
\end{aligned}
\end{equation}
Finally note that the (\ref{eq:4.24bis}) is implied by (\ref{45}) that has been proved in Lemma \ref{lemma2}.

\noindent
${\bf (S4)}_{\nu+1}$. Let $\la\in\Lambda_{\nu+1}^{\g}$, then
by (\ref{eq:419bis}) and the inductive hypothesis ${\bf (S4)_{\nu}}$ one has
that $\Lambda_{\nu+1}^{\g}({\bf u}_{1})\subseteq\Lambda_{\nu}^{\g}({\bf u}_{1})\subseteq\Lambda_{\nu}^{\g-\rho}({\bf u}_{2})\subseteq\Lambda_{\nu}^{\g/2}({\bf u}_{2})$.
Hence the eigenvalues $\mu_{h}^{\nu}(\la,{\bf u}_{2}(\la))$ are well defined by the ${\bf (S1)_{\nu}}$.
Now, since $\la\in\Lambda_{\nu}^{\g}({\bf u}_{1})\cap\Lambda_{\nu}^{\g/2}({\bf u}_{2})$,
we have for $h=(\s,j) \in\CC\times\NNN$
and setting $h'=(\s',j')\in\CC\times\NNN$
\begin{equation}\label{eq:4.2.22}
\begin{aligned}
|(\mu_{h}^{\nu}&-\mu_{h'}^{\nu})(\la,{\bf u}_{2}(\la))-(\mu_{h}^{\nu}-\mu_{h'}^{\nu})(\la,{\bf u}_{1}(\la))|
\stackrel{(\ref{eq:3.2.4})}{\leq} 
|(\mu_{h}^{0}-\mu_{h'}^{0})(\la,{\bf u}_{2}(\la))-(\mu_{h}^{0}-\mu_{h'}^{0})(\la,{\bf u}_{1}(\la))|\\
&+2\sup_{h\in\CC\times\NNN}|r_{h}^{\nu}(\la,{\bf u}_{2}(\la))-r_{h}^{\nu}(\la,{\bf u}_{1}(\la))|
\stackrel{(\ref{eq:4.24bis})}{\leq} 
\e C|\s j^{2}-\s'j'^{2}| ||{\bf u}_{2}-{\bf u}_{1}||_{\gots_{0}+\h_{2}},
\end{aligned}
\end{equation}
The (\ref{eq:4.2.22}) implies that for any $|\ell|\leq N_{\nu}$ and $j\neq j'$,
\begin{equation}\label{eq:4.2.23}
\begin{aligned}
|i\oo\cdot\ell+\mu_{h}^{\nu}({\bf u}_{2})-\mu_{h'}^{\nu}({\bf u}_{2})|
&\stackrel{(\ref{eq:419bis}),(\ref{eq:4.2.22})}{\geq}
\g|\s j^{2}-\s' j'^{2}|\langle\ell\rangle^{-\tau}-C|\s j^{2}-\s'j'^{2}| ||{\bf u}_{2}-{\bf u}_{1}||_{\gots_{0}+\h_{2}}\\
&\stackrel{{\bf (S4)_{\nu}}}{\geq}(\g-\rho)|\s j^{2}-\s' j'^{2}|\langle\ell\rangle^{-\tau},
\end{aligned}
\end{equation}
where we used that, for any $\la\in\Lambda_{0}$, one has $C\e N^{\tau}_{\nu}||{\bf u}_{1}-{\bf u}_{2}||_{\gots_{0}+\h_{2}}\leq \rho$.
Now, the (\ref{eq:4.2.23}), 
 imply 
that if $\la\in\Lambda_{\nu+1}^{\g}({\bf u}_{1})$ then 
$\la\in\Lambda_{\nu+1}^{\g-\rho}({\bf u}_{2})$, that is the ${\bf (S4)_{\nu+1}}$.

\subsubsection{Proof of Theorem \ref{KAMalgorithm}}

We want apply Lemma \ref{teo:KAM} to the linear operator
$\calL_{0}=\calL_{4}$ defined in (\ref{eq:3.5.9})  where $\RR_{0}:=E_{1}^{0}D+E_{0}^{0}$ defined
in (\ref{eq:4.12bis}),
and we have defined for $s\in[\gots_{0}, q-\h_{1}-\be]$,
$\de_{s}^{0}:=|E_{1}^{0}|_{s,\g}+|E_{0}^{0}|_{s,\g},$
then
\begin{equation}\label{eq:4.1.2}
\de_{\gots_{0}+\be}^{0}\stackrel{(\ref{eq:3.2.7})}{\leq}\e C(\gots_{0}+\be)(1+||{\bf u}||_{\be+\gots_{0}+\h_{1},\g})
\stackrel{(\ref{eq:4.2})}{\leq}2\e C(\gots_{0}+\be)
, \qquad \Rightarrow \qquad
N_{0}^{C_0}\de_{\gots_{0}+\be}^{0}\g^{-1}\leq 1,
\end{equation}
if $\e\g^{-1}\leq \epsilon_{0}$ is small enough, that is the (\ref{eq:4.15}). We first prove that there exists
a final transformation $\Phi_{\infty}$. 
For any $\la\in\cap_{\nu\geq0}\Lambda_{\nu}^{\g}$ we define
\begin{equation}\label{eq:4.1.3}
\tilde{\Phi}_{\nu}:=\Phi_{0}\circ\Phi_{1}\circ\ldots\circ\Phi_{\nu}.
\end{equation}
One can note that $\tilde{\Phi}_{\n+1}=\tilde{\Phi}_{\nu}\circ\Phi_{\nu+1}=\tilde{\Phi}_{\nu}+
\tilde{\Phi}_{\nu}\Psi_{\nu+1}$. Then, one has
\begin{equation}\label{eq:4.1.4}
|\tilde{\Phi}_{\nu+1}|_{\gots_{0},\g}\stackrel{(\ref{eq:2.11b})}{\leq}
|\tilde{\Phi}_{\nu}|_{\gots_{0},\g}+C|\tilde{\Phi}_{\nu}|_{\gots_{0},\g}|\Psi_{\nu+1}|_{\gots_{0},\g}
\stackrel{(\ref{eq:4.22})}{\leq}
|\tilde{\Phi}_{\nu}|_{\gots_{0},\g}(1+\e^{(\gots_0)}_{\nu}),
\end{equation}
where we have defined for $s\geq \gots_0$,
\begin{equation}\label{eq:4.1.5}
\e_{\nu}^{(s)}:=K \g^{-1}N_{\nu+1}^{2\tau+1}N_{\nu}^{-\al}\de_{s}^{0}, 
\end{equation}
for some constant $K>0$. Now, by iterating (\ref{eq:4.1.4}) and using the (\ref{eq:4.15}), (\ref{eq:4.22}),
we obtain
\begin{equation}\label{eq:4.1.6}
|\tilde\Phi_{\nu+1}|_{\gots_{0},\g}\leq|\tilde\Phi_{0}|_{\gots_{0},\g}\prod_{\nu\geq0}(1+\e^{(\gots_{0})}_{\nu})
\leq2
\end{equation}

The estimate on the high norm follows by
\begin{equation}\label{eq:4.1.7}
\begin{aligned}
|\tilde{\Phi}_{\nu+1}|_{s,\g}&\stackrel{(\ref{eq:2.11a}),(\ref{eq:4.1.6})}{\leq}
|\tilde\Phi_{\nu}|_{s,\g}(1+C(\gots_0)|\Psi_{\nu+1}|_{\gots_{0},\g})+ C(s)|\tilde\Phi_{\nu}|_{\gots_{0},\g}|\Psi_{\nu+1}|_{s,\g} \\
&\stackrel{(\ref{eq:4.22}),(\ref{27})}{\leq}|\tilde{\Phi}_{\nu}|_{s,\g}(1+\e^{(\gots_{0})}_{\nu})+\e^{(s)}_{\nu}\leq C\left(
\sum_{j=0}^{\infty}\e_{j}^{(s)}+|\tilde\Phi_{0}|_{s,\g}\right)
\stackrel{(\ref{eq:4.22})}{\leq}C(s)\left(1+\de_{s+\be}^{0}\g^{-1}\right)
\end{aligned}
\end{equation}
where we used the inequality $\prod_{j\geq0}(1+\e_{j}^{(\gots_{0})})\leq 2$.
Thanks to (\ref{eq:4.1.7}) we can prove that the sequence $\tilde{\Phi}_{\n}$ is a Chauchy sequence in norm $|\cdot|_{s}$. Indeed,
\begin{equation}\label{eq:4.1.8}
\begin{aligned}
|\tilde{\Phi}_{\nu+m}-\tilde{\Phi}_{\nu}|_{s,\g}&\leq\sum_{j=\nu}^{\nu+m-1}
|\tilde{\Phi}_{j+1}-\tilde{\Phi}_{j}|_{s,\g}\stackrel{(\ref{eq:2.11a})}{\leq}C(s)\sum_{j=\nu}^{\nu+m-1}
(|\tilde\Phi_{j}|_{s,\g}|\Psi_{j+1}|_{\gots_{0},\g}+|\tilde\Phi_{j}|_{\gots_{0},\g}|\Psi_{j+1}|_{s,\g}\\
&\stackrel{(\ref{eq:4.22}),(\ref{eq:4.1.6}),(\ref{eq:4.1.7}),(\ref{eq:4.15})}{\leq}
C(s)\sum_{j\geq\nu}\de_{s+\be}^{0}\g^{-1}N_{j}^{-1}
\leq C(s)
\de_{s+\be}^{0}\g^{-1}N_{\nu}^{-1}.
\end{aligned}
\end{equation}
As consequence one has that $\tilde{\Phi}_{\nu}\stackrel{|\cdot|_{s,\g}}{\to}\Phi_{\infty}$, Moreover,
 (\ref{eq:4.1.8}) used with $m=\infty$ and $\nu=0$ and $|\tilde{\Phi}_{0}-\id|_{s,\g}=|\Psi_{0}|_{s,\g}\leq
 \g^{-1}\de_{s+\be}^{0}$ imply
\begin{equation}\label{eq:4.1.9}
\begin{aligned}
|\Phi_{\infty}-\id|_{s,\g}&\leq C(s)\g^{-1}\de_{s+\be}^{0},\qquad
|\Phi_{\infty}^{-1}-\id|_{s,\g}&
\stackrel{(\ref{eq:2.14})}{\leq}C(s)\g^{-1}\de_{s+\be}^{0}.
\end{aligned}
\end{equation}
Hence the (\ref{eq:4.8}) is verified.

Let us now define for $k=(\s,j)\in \CC\times\NNN$,
\begin{equation}\label{eq:4.1.1}
\mu_{k}^{\infty}:=\mu_{\s,j}^{\infty}(\la)=\lim_{\nu\to+\infty}\tilde{\mu}^{\nu}_{\s,j}(\la)
=\tilde{\mu}^{0}_{\s,j}(\la)+
\lim_{\nu\to+\infty}\tilde{r}_{\s,j}^{\n}.
\end{equation}
We can note that, for any $\nu,j\in \NNN$, the following important estimates on the eigenvalues hold:
\begin{equation}\label{eq:4.1.10}
\begin{aligned}
|\mu_{k}^{\infty}-\tilde{\mu}_{k}^{\nu}|_{\Lambda,\g}&\leq
\sum_{m=\nu}^{\infty}|\tilde{\mu}^{m+1}_{k}-\tilde{\mu}_{k}^{m}|_{\Lambda,\g}
\stackrel{(\ref{eq:4.23}),(\ref{eq:4.21})}{\leq}
C \de_{\gots_{0}+\be}^{0}N_{\nu-1}^{-\al},
\end{aligned}
\end{equation}
and moreover,
\begin{equation}\label{eq:4.1.11}
|\mu_{k}^{\infty}-\tilde{\mu}^{0}_{k}|_{\Lambda,\g} \leq 
C\de_{\gots_{0}+\be}^{0}.
\end{equation}

As seen in Lemma \ref{teo:KAM}, the corrections $r_{\s,j}^{\nu}=(E_{0}^{\nu})_{k}^{k}=
(E_{0}^{\nu})_{\s,j}^{\s,j}(0)$.

The following Lemma gives us a connection between the Cantor sets defined in Lemma
\ref{teo:KAM} and Theorem \ref{KAMalgorithm}.

\begin{lemma}\label{lem:4.4}
One has that
\begin{equation}\label{eq:4.1.12}
\Lambda^{2\g}_{\infty}\subset\cap_{\nu\geq0}\Lambda_{\nu}^{\g}.
\end{equation}
\end{lemma}

\prova
Consider $\la\in\Lambda^{2\g}_{\infty}$. We show by induction that for any $\nu>0$ then 
$\la\in\Lambda_{\nu}^{\g}$,
since by definition we have $\Lambda^{2\g}_{\infty}\subset\Lambda_{0}^{\g}:=\Lambda_{o}$.
Assume that $\Lambda^{2\g}_{\infty}\subset \Lambda_{\nu-1}^{\g}$. Hence $\mu^{\nu}_{h}$ are well defined and coincide with their extension.
Then, for any fixed  $k=(\s,j,p), k'=(\s',j',p')\in\CC\times\NNN\times\ZZZ^{d}$,
 we have
\begin{equation}\label{eq:4.1.13}
\begin{aligned}
|\oo\cdot\ell+\mu_{\s,j}^{\nu}-\mu_{\s',j'}^{\nu}|
&\stackrel{(\ref{martina10}),(\ref{eq:4.1.10})}{\geq}\frac{2\g|\s j^{2}-\s' j'^{2}|}{\langle\ell\rangle^{\tau}}
-2C\de_{\gots_{0}+\be}^{0}N_{\nu-1}^{-\al}.
\end{aligned}
\end{equation}
Now, by the smallness hypothesis (\ref{eq:4.15}), 
we can estimate for  $|p-p'|=|\ell|\leq N_{\nu}$,
\begin{equation}\label{eq:4.1.15}
|\oo\cdot\ell+\mu_{\s, j}^{\nu}-\mu_{\s',j'}^{\nu}|\geq\frac{\g|\s j^{2}-\s' j'^{2}|}{\langle\ell\rangle^{\tau}},
\end{equation}
that implies $\la\in\Lambda_{\nu}^{\g}$.
\EP

Now, for any $\la\in \Lambda^{2\g}_{\infty}\subset\cap_{\n\geq0}\Lambda_{\n}^{\g}$ (see (\ref{eq:4.1.12})),
one has
\begin{equation}\label{100}
|\DD_{\n}-\DD_{\infty}|_{s,\g}=
\sup_{k\in \CC\times\NNN\times\ZZZ^{d}}| \mu_{\s,j}^{\n}-\mu_{\s',j'}^{\infty}|_{\g}
\stackrel{(\ref{eq:4.1.10}), (\ref{eq:4.1.11})}{\leq}K \de_{\gots_{0}+\be}^{0}N_{\n-1}^{-\al}, \quad
\de_{s}^{\n}\stackrel{(\ref{eq:4.21})}{\leq}\de_{s+\be}^{\be}N_{\n-1}^{-\al},
\end{equation}
that implies
\begin{equation}\label{101}
\calL_{\n}\stackrel{(\ref{eq:4.16})}{=}\DD_{\n}+\RR_{\n}\stackrel{|\cdot|_{s,\g}}{\to}
\DD_{\infty}=:\calL_{\infty}, \quad \DD_{\infty}:=diag_{k\in C\times\NNN\times\ZZZ^{n}}\mu_{k}^{\infty}.
\end{equation}

By applying iteratively the (\ref{eq:4.16}) we obtain $\calL_{\nu}=\tilde{\Phi}_{\n-1}^{-1}\calL_{0}\tilde{\Phi}_{\n-1}$
where $\tilde{\Phi}_{\n-1}$ is defined in (\ref{eq:4.1.3}) and, by (\ref{eq:4.1.8}), $\tilde{\Phi}_{\n-1}\to\Phi_{\infty}$
in norm $|\cdot|_{s,\g}$.
Passing to the limit we get
\begin{equation}\label{102}
\calL_{\infty}=\Phi_{\infty}^{-1}\circ\calL_{0}\circ\Phi_{\infty},
\end{equation}
that is the (\ref{eq:4.6}), 
while the (\ref{eq:4.5}) follows by (\ref{eq:4.1.2}), (\ref{eq:4.1.10}) and (\ref{eq:4.1.11}).
Finally, (\ref{eq:2.11a}), (\ref{eq:2.13b}), Lemma \ref{1.4} and (\ref{eq:4.8}) implies the bounds (\ref{eq:4.9}). 
This concludes the proof.
\EP

\setcounter{equation}{0}
\section{Conclusion of the diagonalization algorithm  and  inversion of $\calL({\bf u})$}\label{sec:5}

In the previous Section we have conjugated the operator $\calL_{4}$ (see (\ref{eq:3.5.9}))
to a diagonal operator $\calL_{\infty}$. 
In conclusion, we have that
\begin{equation}\label{eq:4.4.1}
\calL=W_{1}\calL_{\infty}W_{2}^{-1}, \quad W_{i}=\VV_{i}\Phi_{\infty}, \quad 
 \VV_{1}:=\TT_{1}\TT_{2}\TT_{3}\rho\TT_{4}, \quad
\VV_{2}=\TT_{1}\TT_{2}\TT_{3} \TT_{4}.
\end{equation}

We have the following result
\begin{lemma}\label{lemma5.8}
Let $\gots_{0}\leq s\leq q-\be-\h_{1}-2$, with $\h_{1}$ define in (\ref{eq:3.2.0}) and $\be$  in Theorem (\ref{KAMalgorithm}). Then, for $\e\g^{-1}$ small enough, and
\begin{equation}\label{eq:4.4.2}
||{\bf u}||_{\gots_{0}+\be+\h_{1}+2,\g}\leq1,
\end{equation}
one has for any $\la\in\Lambda^{2\g}_{\infty}$,
\begin{equation}\label{eq:4.4.3}
\begin{aligned}
||W_{i}{\bf h}||_{s,\g}+||W_{i}^{-1}{\bf h}||_{s,\g}&
\leq C(s)\left(||{\bf h}||_{s+2,\g}+
||{\bf u}||_{s+\be+\h_{1}+4,\g}||{\bf h}||_{\gots_{0},\g}\right),
\end{aligned}
\end{equation}
for $i=0,1$. Moreover, $W_{i}$ and $W_{i}^{-1}$ are reversible preserving.
\end{lemma}

\prova
Each $W_i$  is composition of two operators, the $\VV_i$ satisfy the (\ref{eq:3.2.3}) while $\Phi_\infty$ satisfies (\ref{eq:4.8}). We use  (\ref{eq:2.13b}) in order to pass to the operator norm. Then Lemma \ref{lem5} and \eqref{A2} with $p=s-\gots_0$, $q=2$ implies the bounds (\ref{eq:4.4.3}). Moreover the transofrmations $W_{i}$ and $W_{i}^{-1}$
 are reversibility preserving because each transformations $\VV_{i}$,$\VV_{i}^{-1}$ and
 $\Phi_{\infty}$, $\Phi_{\infty}^{-1}$ is reversibility preserving.
 \EP 

\noindent
\subsection{Proof of Proposition \ref{teo2}}
We fix $\h=\h_{1}+\be+2$ and $q>\gots_{0}+\h$.
Let $\mu^{\infty}_{h}$ be the functions defined in (\ref{eq:4.1.1}).
Then by 
Theorem \ref{KAMalgorithm} and Lemma \ref{lemma5.8}
for $\la\in \Lambda_{\infty}^{2\g}$ we have
the (\ref{1.2.2}). Hence item (i) is proved.

Item (ii) follows by applying the dynamical system point of view. We have already proved that
\begin{equation}\label{eq:6.2}
\calL=\TT_{1}\TT_{2}\TT_{3}\rho\TT_{4}\Phi_{\infty}\calL_{\infty}\Phi^{-1}_{\infty}\TT_{4}^{-1}\TT_{3}^{-1}
\TT_{2}^{-1}\TT_{1}^{-1}.
\end{equation}
 By Lemma \ref{lem:3.9} all the changes of variables in \eqref{eq:6.2} can be seen as transformations of the phase space
${\bf H}^{s}_{x}$ depending in a quasi-periodic way on time plus quasi periodic reparametrization of time ($\TT_3$).
With this point of view, consider a dynamical system of the form
\begin{equation}\label{eq:6.3}
\del_{t}{\bf u}=L(\oo t){\bf u}.
\end{equation}
Under a transformation of the form ${\bf u}=A(\oo t){\bf v}$, one has that the system (\ref{eq:6.3})
become
\begin{equation}\label{eq:6.4}
\del_{t}{\bf v}=L_{+}(\oo t){\bf v}, \qquad
 L_{+}(\oo t)=A(\oo t)^{-1} L(\oo t)A(\oo t)-A(\oo t)^{-1}\del_{t}A(\oo t)
\end{equation}
The transformation $A(\oo t)$ acts on the functions ${\bf u}(\f,x)$ as
\begin{equation}\label{eq:6.5}
( {A}{\bf u})(\f,x):=(A(\f){\bf u}(\f,\cdot))(x):=A(\f){\bf u}(\f,x), \qquad ({A}^{-1}{\bf u})(\f,x)=A^{-1}(\f){\bf u}(\f,x).
\end{equation}
Then the operator on the quasi-periodic functions
\begin{equation}\label{eq:6.6}
\calL:=\oo\cdot\del_{\f}-L(\f),
\end{equation}
associated to the system (\ref{eq:6.3}), is transformed by $ {A}$ into
\begin{equation}\label{eq:6.7} 
 {A}^{-1}\calL  {A}=\oo\cdot\del_{\f}-L_{+}(\f),
\end{equation}
that represent the system in (\ref{eq:6.4}) acting on quasi-periodic functions.
The same considerations hold for transformations of the type
\begin{equation}\label{eq:6.8}
\begin{aligned}
&\tau:=\psi(t):=t+\al(\oo t), \quad t=\psi^{-1}(\tau):=\tau+\tilde{\al}(\oo \tau),\\ 
& (B{\bf u})(t):={\bf u}(t+\al(\oo t)), \qquad (B^{-1}{\bf v})(\tau)=v(\tau+\tilde{\al}(\oo \tau)).
\end{aligned}
\end{equation}
with $\al(\f)$, $\f\in\TTT^{d}$ is $2\pi-$periodic in all the $d$ variables. The operator $B$
is nothing but the operator on the functions induced by the diffeomorfism of the torus
$t\to t+\al(\oo t)$.
The transformation ${\bf u}=B{\bf v}$ transform the system (\ref{eq:6.3}) into 
\begin{equation}\label{eq:6.9}
\del_{t}{\bf v}=L_{+}(\oo t){\bf v}, \qquad
 L_{+}(\oo \tau):=\left(\frac{L(\oo t)}{1+(\oo\cdot\del_{\f}\al)(\oo t)}\right)_{| t=\tilde{\psi}(\tau)}
\end{equation}
If we consider the operator $B$ acting on the quasi-periodic functions as 
$(B{\bf u})(\f,x)={\bf u}(\f+\oo\al(\f),x)$ and $(B^{-1}{\bf u})(\f,x):={\bf u}(\f+\oo\tilde{\al}(\f),x)$, 
we have that
\begin{equation}\label{eq:6.10}
B^{-1}\calL B=\rho(\f)\calL_{+}=\rho(\f)\left(\oo\cdot\del_{\f}-L_{+}(\f)\right)=
\rho(\f)\left(\oo\cdot\del_{\f}-\frac{1}{\rho(\f)}L(\f+\oo\tilde{\al}(\f))
\right),
\end{equation}
and $\rho(\f):=B^{-1}(1+\oo\cdot\del_{\f}\al)$, that means that $\calL_{+}$ is the linear system
(\ref{eq:6.9}) acting on quasi-periodic functions.

By these arguments, we have simply that
a curve ${\bf u}(t)$ in the phase space of functions of $x$, i.e. ${\bf H}^{s}_{x}$, 
solves the linear dynamical system (\ref{1.2.1}) if and only if
the curve
\begin{equation}\label{eq:6.11}
{\bf v}(t):= \Phi_{\infty}^{-1}\TT_{4}^{-1} {\TT_{3}}^{-1} {\TT_{2}}^{-1}\TT_{1}^{-1}(\oo t){\bf h}(t)
\end{equation}
solves the system (\ref{1.2.6}). This completely justify Remark \ref{ciccio}.  In Lemma \ref{lem:3.9} and the (\ref{eq:4.9}) 
 we have checked that these transformations are well defined.
\EP

\noindent
The result of Proposition \ref{teo2} holds for $\la$ in a suitable Cantor set. 

\subsection{Proof of Lemma \ref{inverseofl}}
 
As explained in the Introduction, we now study the invertibility of
\begin{equation}\label{eq:4.4.5}
\calL_{\infty}:=diag_{k\in\CC\times\NNN\times\ZZZ^{d}}\{i\oo\cdot\ell+\mu_{\s,j}^{\infty}\}, \quad
\mu_{\s,j}^{\infty}\la=-i\s m(\la)j^{2}+r_{\s,j}^{\infty}(\la).
\end{equation}
in order to obtain a better understanding of the set $\calG_{\infty}$ of the Nash-Moser Proposition
\ref{teo4}.
.

\begin{lemma}\label{inverselinfty}
For  ${\bf g}\in{\bf Z}^{s}$,
 consider the equation
\begin{equation}\label{eq:4.4.7}
\calL_{\infty}({\bf u}){\bf h}={\bf g}.
\end{equation}
If $\la\in \Lambda^{2\g}_{\infty}({\bf u})\cap P^{2\g}_{\infty}({\bf u})$  (defined respectively in \eqref{martina10} and \eqref{eq:primedimerda}),  then there exists a unique solution 
$\calL_{\infty}^{-1}{\bf g}:={\bf h}=(h,\bar{h})\in{\bf X}^{s}$. Moreover, for all Lipschitz family ${\bf g}:=
{\bf g}(\la)\in {\bf Z}^{s}$ one has
\begin{equation}\label{eq:4.4.8}
||\calL_{\infty}^{-1}{\bf g}||_{s,\g}\leq C \g^{-1}||{\bf g}||_{s+2\tau+1,\g}.
\end{equation}
\end{lemma}

\prova By solving the (\ref{eq:4.4.7})  one obtain the solution ${\bf h}:=(h_{+},h_{-})$ of the form
\begin{equation}\label{eq:4.4.9}
\begin{aligned}
h_{+}(\f,x)&:=\sum_{\ell\in\ZZZ^{d}, j\geq1} \frac{g_{j}(\ell)}{i\oo\cdot\ell+\mu_{1,j}^{\infty}}
e^{i\ell\cdot\f}\sin jx, \\
h_{-}(\f,x)&:=
\sum_{\ell\in\ZZZ^{d}, j\geq1 } 
\frac{\ol{g_{j}(-\ell)}}{i\oo\cdot\ell+\mu_{-1,j}^{\infty}}
e^{i\ell\cdot\f}\sin jx=\sum_{\ell\in\ZZZ^{d}, j\geq1 } 
\frac{\ol{g_{j}(\ell)}}{-i\oo\cdot\ell+\mu_{-1,j}^{\infty}}e^{-i\ell\cdot\f}\sin jx.
\end{aligned}
\end{equation}
Now, by the hypothesis of reversibility, we have already seen that 
$\mu_{1,j}^{\infty}=-\ol{\mu_{1,j}^{\infty}}$ and $\mu_{-1,j}^{\infty}=-\mu_{1,j}^{\infty}$, then
one has that $\ol{h_{-}}=h_{+}:=h$. 
Moreover, one has
\begin{equation}\label{eq:4.4.10}
\ol{h_{j}(\ell)}=\ol{\frac{g_{j}(\ell)}{i\oo\cdot\ell+\mu_{1,j}^{\infty}}}
{=}\frac{-g_{j}(\ell)}{-(i\oo\cdot\ell+\mu_{1,j}^{\infty})}
=h_{j}(\ell)
\end{equation}
then the Lemma (\ref{lemFou}) implies that ${\bf h}\in{\bf X}^{s}$.

Now, since $\la\in \Lambda^{2\g}_{\infty}({\bf u})\cap P^{2\g}_{\infty}({\bf u})$
then, by (\ref{eq:primedimerda}), we can estimate the (\ref{eq:4.4.9})
\begin{equation}\label{eq:4.4.11}
||h||_{s}\leq C \g^{-1}||{\bf g}||_{s+\tau}.
\end{equation}
The  Lipschitz bound on $h$ follow  exactly as in formul\ae (\ref{fiorellino2})-(\ref{eq:4.1.116}) 
and we obtain
\begin{equation}\label{eq:4.4.17}
||{\bf h}||_{s,\g}=||{\bf h}||^{{\rm sup}}_{s}+\g||{\bf h}||^{{\rm lip}}_{s}\lessdot
\g^{-1}||{\bf g}||_{s+2\tau+1,\g},
\end{equation}
that is the (\ref{eq:4.4.8}).
\EP

\noindent
Since $W_{i}^{\pm1}$ are reversibility preserving, we show in the next Lemma how to solve the
equation $\calL {\bf h}={\bf g}$ for ${\bf g}\in {\bf Z}^{s}$:

\smallskip
\noindent {\it PROOF OF LEMMA \ref{inverseofl}}.
By (\ref{eq:4.4.1}) one has that the equation
$\calL{\bf h}={\bf g}$ si equivalent to $\calL_{\infty}W_{2}^{-1}{\bf h}=W_{1}^{-1}{\bf g}$.
By Lemma \ref{inverselinfty} this second equation has a unique solution
 $W_{2}^{-1}{\bf h}\in {\bf X}^{s}$. Note that this is true because $W_{1}^{-1}$ is reversibility-preserving,
so that $W_{1}^{-1}{\bf g}\in {\bf Z}^{s}$ if ${\bf g}\in{\bf Z}^{s}$. Hence the solution with zero average
of 
$\calL{\bf h}={\bf g}$ is of the form
\begin{equation}\label{eq:4.4.24}
{\bf h}:=
W_{2}\calL_{\infty}^{-1}W_{1}^{-1}{\bf g},
\end{equation}
 Moreover, one has that
${\bf h}\in {\bf X}^{s}$, because $W_{2}$ is reversibility-preserving, and because, 
by Lemma \ref{inverselinfty} one has that $\calL_{\infty}^{-1} : {\bf  Z}^{0} \to {\bf X}^{0} $.

Now we have  
\begin{equation}\label{eq:4.4.25}
\begin{aligned}
||{\bf h}||_{s,\g}&\stackrel{(\ref{eq:4.4.3})}{\leq} C(s)\left(
||\calL_{\infty}^{-1}W_{1}^{-1}{\bf g}||_{s+2,\g}+
||{\bf u}||_{s+\be+\h_{1}+4,\g}||\calL_{\infty}^{-1}
W_{1}^{-1}{\bf g}||_{\gots_{0},\g}
\right)\\
&\leq C(s)\g^{-1} \left( ||{\bf g}||_{s+2\tau+5,\g}+
||{\bf u}||_{s+4\tau+\be+10+\h_{1},\g}||{\bf g}||_{\gots_{0},\g}
\right),
\end{aligned}
\end{equation}
where, in the second inequality we used (\ref{eq:4.4.8}) on $\calL_{\infty}^{-1}$, again the (\ref{eq:4.4.3}) for $W_{1}^{-1}$ and (\ref{eq:4.4.19}). Finally we used the (\ref{A2}) with $a=\gots_{0}+2\tau+\h_{1}+\be+7, b=\gots_{0}$ and  
$p=s-\gots_{0}, q=2\tau+3$.
The (\ref{eq:4.4.25}) implies  the (\ref{eq:4.4.23}) with $\zeta$ 
defined in (\ref{eq:4.4.18}) where we already fixed $\h:=\h_{1}+\be+2$ in the proof of Proposition \ref{teo2}.
\EP


\setcounter{equation}{0}
\section{Measure estimates and conclusions}
\label{sec6}

The aim of this Section is to use the information obtained in Sections $\ref{sec:3}$ and $\ref{sec:4}$,
in order to apply Theorem \ref{NM} to our problem and prove Theorem \ref{teo1}.
First of all we prove the approximate reducibility  Lemma \ref{teo3}.

\smallskip
\noindent {\it PROOF OF LEMMA \ref{teo3}}. We first apply the change of variables defined in \eqref{tra} to $\calL(\bf v)$ in order to reduce to $\calL_4(\bf v)$. We know that Lemma \ref{teo:KAM} holds for $\calL_4(\bf u)$, now we fix $\nu$  such that $N_{\nu-1}\leq N\leq N_{\nu}$ and apply $({\bf S3}_{\nu})-({\bf S4}_{\nu})$ with ${\bf u}_1={\bf u}$, ${\bf u}_2={\bf v}$. This implies our claim since, by Lemma \ref{lem:4.4}, we have  $\Lambda_\infty^{2\g}({\bf u})\subseteq \Lambda_{\nu}^\g({\bf u})\subseteq \Lambda_{\nu}^{\g-\rho}({\bf v})$. Finally for all $\la\in\Lambda_{\nu+1}^{\g-\rho}({\bf v})$ we can perform $\nu+1$ steps in Lemma \ref{teo:KAM}. Fixing $\kappa= 2\alpha/3$ we obtain the bounds on the changes of variables and remainders,  using formul\ae\phantom{.} \eqref{60} and \eqref{67}.

\subsection{Proof of Proposition \ref{measure}.}
Recall that we have set
$$
\g_{n}:=\g\left(1+\frac{1}{2^{n}}\right),
$$
 $({\bf u}_{n})_{\geq0}$ is the sequence of approximate solutions introduced in Theorem 
\ref{NM}. which  is well defined in $\calG_{n}$ and satisfies the hypothesis
of Proposition \ref{teo2}.   $\calG_{n}$ in turn is defined in $(N1)_n$ and Definition \ref{invertibility}. 
For notational convenience we extend
 the eigenvalues $\mu_{\s,j}^{\infty}({\bf u}_{n})$ introduced in  Proposition \ref{teo2}) , which are  defined only for $j\in\NNN$,
to a function defined for $j\in\ZZZ_{+} $ in the following way:
\begin{equation}\label{eq141}
\Omega_{\s,j}({\bf u}_{n}):=
\mu_{\s,j}^{\infty}({\bf u}_{n}), \quad (\s,j)\in\CC\times\NNN, \qquad
\Omega_{\s,j}({\bf u}_{n})\equiv0,\quad \;\;\;\;\s\in\CC, \;\; j=0.
\end{equation}
We now define inductively a sequence of nested sets $G_{n}\subseteq\calG_{n}$ for $n\geq0$. 
Set $G_{0}=\Lambda$ and 
\begin{equation}\label{eq142bis}
G_{n+1}:=
\left\{\begin{aligned}
 \la\in G_{n} \; :\; &|i\oo\cdot\ell+\Omega_{\s,j}({\bf u}_{n})-\Omega_{\s',j'}({\bf u}_{n})|
\geq\frac{2\g_{n}|\s j^{2}-\s' j'^{2}|}{\langle\ell\rangle^{\tau}}, \\
&\;\forall\ell\in\ZZZ^{n},\;\;  \s,\s'\in\CC, \;\; j,j'\in\ZZZ_{+}
\end{aligned}\right\},
\end{equation}
The following Lemma implies \eqref{eq136b}.
\begin{lemma}\label{megalemma} Under the Hypotheses of Proposition \ref{measure},
for any $n\geq0$, one has
\begin{equation}\label{eq137}
G_{n+1}\subseteq \calG_{n+1}.
\end{equation}
\end{lemma}

\prova
For any $n\geq0$
and if $\la\in G_{n+1}$,  one has 
by Lemmata
\ref{inverselinfty} and \ref{inverseofl}, (recalling that $\g\leq \g_{n}\leq 2\g$ and $2\tau+5<\zeta$)
\begin{equation}\label{eq139}
\begin{aligned}
||\calL^{-1}({\bf u}_{n}){\bf g}||_{s,\g}&\leq C(s)\g^{-1}\left(
||{\bf g}||_{s+\zeta,\g}+||{\bf u}_{n}||_{s+\zeta,\g}||{\bf g}||_{\gots_{0},\g}\right),\\
||\calL^{-1}({\bf u}_{n})||_{\gots_{0},\g}&\leq C(\gots_{0})\g^{-1}
N_{n}^{\zeta}||{\bf g}||_{\gots_{0},\g},
\end{aligned}
\end{equation}
for $\gots_{0}\leq s\leq q-\mu$, for any ${\bf g}(\la)$ Lipschitz family. 
The (\ref{eq139}) are nothing but the (\ref{eq104}) in Definition \ref{invertibility} with $\mu=\zeta$ .
It represents the loss of regularity that you have when you perform the
regularization procedure in Section \ref{sec:3} and during the diagonalization algorithm in
Section \ref{sec:4}. This justifies our choice of $\mu$ in Proposition \ref{measure}.
\EP

\noindent
By Lemma \ref{megalemma}, in order to obtain the bound (\ref{eq136}), it is enough to prove that
\begin{equation}\label{eq140}
|\Lambda\backslash \cap_{n\geq0}G_{n}|\to 0, \;\;\; {\rm as}\;\;\;\; \g\to0.
\end{equation}
We will prove by induction that, for any $n\geq0$, one has
\begin{equation}\label{eq144}
|G_{0}\backslash G_{1}|\leq C_{\star}\g, \quad |G_{n}\backslash G_{n+1}|\leq C_{\star}\g N_{n}^{-1},\;\;\; n\geq1.
\end{equation}
First of all, write
\begin{equation}\label{eq145}
\begin{aligned}
&G_{n}\backslash G_{n+1}:=\bigcup_{\substack{\s,\s'\in\CC, j,j'\in\ZZZ_{+} \\ \ell\in\ZZZ^{n}}} 
R_{\ell j j'}^{\s,\s'}({\bf u}_{n})\\
R_{\ell j j'}^{\s,\s'}({\bf u}_{n})&:=\left\{\la \in G_{n}\; : \; 
|i\la\bar{\oo}\cdot\ell+\Omega_{\s,j}({\bf u}_{n})-
\Omega_{\s',j'}({\bf u}_{n})|<2\g_{n}|\s j^{2}-\s' j'^{2}|
\langle\ell\rangle^{-\tau}
\right\}.
\end{aligned}
\end{equation}
Assume in the following that, if $\s=\s'$, then $j\neq j'$, since one has 
$R_{\ell j j}^{\s,\s}({\bf u_{n}})=\emptyset$.
Important properties of the sets $R_{\ell j j'}^{\s,\s'}$ are the following. The proofs are quite standard
and follow very closely
 Lemmata 5.2 and 5.3 in \cite{BBM1}.  For completeness we give a proof in the Appendix $C$.
\begin{lemma}\label{lemma6.13}
For any $n\geq0$, $|\ell|\leq N_{n}$, one has, for $\e\g^{-1}$ small enough,
\begin{equation}\label{eq146}
R_{\ell j j'}^{\s,\s'}({\bf u}_{n})\subseteq R_{\ell j j'}^{\s,\s'}({\bf u}_{n-1}).
\end{equation}
Moreover, 
\begin{equation}\label{eq147}
{\rm if} \;\;\;\;\; R_{\ell j j'}^{\s,\s'}\neq \emptyset, \;\;\;\;\;\; {\rm then}\;\;\;\;\; 
|\s j^{2}-\s' j'^{2}|\leq 8|\bar{\oo}\cdot\ell|.
\end{equation}
\end{lemma}
\begin{lemma}\label{capperi}
For all $n\geq0$, one has
\begin{equation}\label{eq149}
|R_{\ell j j'}^{\s,\s'}({\bf u}_{n})|\leq C \g \langle\ell\rangle^{-\tau}.
\end{equation}
\end{lemma}

\noindent

We now prove (\ref{eq140}) by assuming Lemmata \ref{lemma6.13} and \ref{capperi}. 
By (\ref{eq145}) one has $R_{\ell j j'}^{\s,\s'}({\bf u}_{n})\subset G_{n}$, and at the same time
for all $|\ell|\leq N_{n}$ one has $R_{\ell j j'}^{\s,\s'}({\bf u}_{n})\subseteq R^{\s,\s'}_{\ell j j'}({\bf u}_{n-1})$ by
(\ref{eq146}). Hence, if $|\ell|\leq N_{n}$, one has $R_{\ell j j'}^{\s,\s'}({\bf u}_{n})=\emptyset$ since
$R_{\ell j j'}^{\s, \s'}({\bf u}_{n-1})\cap G_{n}=\emptyset$ by definition (\ref{eq142bis}).
This implies that
\begin{equation}\label{eq148}
G_{n}\backslash G_{n+1}\subseteq \bigcup_{\substack{ \s,\s'\in\CC, j,j'\in\ZZZ_{+} \\ |\ell|> N_{n}}}
R_{\ell j j'}^{\s, \s'}({\bf u}_{n})
\end{equation}

Now, consider the sets $R_{\ell j j'}^{\s,\s'}(0)$.
By (\ref{eq147}), we know that if $R_{\ell j j'}^{\s,\s'}(0)\neq \emptyset$ then we must have
$j+j'\leq 16 |\bar{\oo}||\ell|$. Indeed, if $\s=\s'$, then
\begin{equation}\label{eq154}
|j^{2}-j'^{2}|=|j-j'|(j+j')\geq \frac{1}{2}(j+j'), \quad \forall j,j'\in\ZZZ_{+}, \;\;j\neq j',
\end{equation}
while, if $\s\neq\s'$, one has $(j+j')/2\leq (j^{2}+j'^{2})\leq 8 |\bar{\oo}||\ell|$ see (\ref{eq147}). 
Then, for $\tau> d+2$, we obtain the first of (\ref{eq144}), by
\begin{equation}\label{eq155}
|G_{0}\backslash G_{1}|\leq \sum_{\substack{\s,\s'\in\CC,\\ j,j'\in\ZZZ_{+} \\ \ell\in\ZZZ^{d}}}
|R_{\ell j j'}^{\s,\s'}(0)|\leq
 \sum_{\substack{\s,\s'\in\CC, \\ (j+j')\leq 16|\bar\oo||\ell| \\ \ell\in\ZZZ^{d}}}
|R_{\ell j j'}^{\s,\s'}(0)|\stackrel{(\ref{eq149})}{\leq} C\g \sum_{\ell\in\ZZZ^{d}}\langle\ell\rangle^{-(\tau-1)}
\leq C\g.
\end{equation}
Finally, we have for any $n\geq1$,
\begin{equation}\label{eq156}
\begin{aligned}
|G_{n}\backslash G_{n+1}|&\stackrel{(\ref{eq148})}{\leq}\!\!
\sum_{\substack{\s,\s'\in\CC, \\ (j+j')\leq 16|\bar\oo||\ell| \\ |\ell|>N_{n}}}\!\!\!
|R_{\ell j j'}^{\s,\s'}({\bf u}_{n})|\stackrel{(\ref{eq149})}{\leq}\!\!\!
\sum_{\substack{ |\ell|>N_{n}}}
 C\g\langle\ell\rangle^{-(\tau-1)}\leq C\g N_{n}^{-\tau+d+1}\leq C \g N_{n}^{-1},
\end{aligned}
\end{equation}
that implies the (\ref{eq144}).
Now we have
\begin{equation}\label{eq157}
|\Lambda\backslash \cap_{n\geq0}G_{n}|\leq \sum_{n\geq0}|G_{n}\backslash G_{n+1}|\leq
C\g +C \g\sum_{n\geq1}N_{n}^{-1}\leq C \g\to 0, \quad {\rm as}\;\;\;\;\; \g\to 0.
\end{equation}
By  (\ref{eq137}), we have that $\cap_{n\geq0}G_{n} \subseteq \calG_{\infty}$. Then, by (\ref{eq157}), we obtain then (\ref{eq136}). 

\EP

\subsection{Proof of Theorem \ref{teo1} }

\noindent
Fix $\g:=\e^{a}$, $a\in(0,1)$. Then the smallness condition $\e\g^{-1}=\e^{1-a}<\epsilon_{0}$ of Theorem
\ref{NM} is satisfied. Then we can apply it with $\mu=\zeta$ in (\ref{eq:4.4.18}) (see
Lemma \ref{measure}).
 Hence by (\ref{eq112}) we have that the function
${\bf u}_{\infty}$ in ${\bf X}^{\gots_{0}+\zeta}$ is a solution of the perturbed NLS with
$\oo=\la\bar{\oo}$. Moreover, one has
\begin{equation}\label{eq:6.1}
|\Lambda\backslash\calG_{\infty}|\stackrel{(\ref{eq136})}{\to}0,
\end{equation}
as $\e$
 tends to zero. To complete the proof of the theorem, it remains to prove 
 the linear stability of the solution.

Since the eigenvalues $\mu_{\s,j}^{\infty}$ are purely imaginary, we know that
the Sobolev norm of the solution ${\bf v}(t)$ of (\ref{1.2.6}) is constant in time. We show
that the Sobolev norm of ${\bf h}(t)=W_{2}^{-1}{\bf v}(t)$, solution of (\ref{1.2.1}) does not grow on time. To do this
we first note that, by (\ref{eq:3.93f}) and (\ref{eq:4.9}), one has
\begin{equation}\label{eq:6.12}
\begin{aligned}
&||\TT_{i}^{\pm1}(\oo t){\bf g}||_{H_{x}^{s}}+||(\TT_{4}\Phi_{\infty})^{\pm1}(\oo t){\bf g}||_{H^{s}_{x}}
\leq C(s)||{\bf g}||_{H_{x}^{s}}, \quad \forall t\in\RRR, \;\;\forall {\bf g}={\bf g}(x)\in {\bf H}_{x}^{s},\\
&||(\TT_{i}^{\pm1}(\oo t)-\id){\bf g}||_{H_{x}^{s}}+
||((\TT_{4}\Phi_{\infty})^{\pm1}(\oo t)-\id){\bf g}||_{H^{s}_{x}}\leq\e\g^{-1}C(s)||{\bf g}||_{H_{x}^{s+1}},
\quad \forall t\in \RRR,\;\; \forall {\bf g}\in{\bf H}_{x}^{s}.
\end{aligned}
\end{equation}
with $i=1,2$.
In both cases, the constant $C(s)$ depends on $||{\bf u}||_{s+\gots_{0}+\be+\h_{1}}$.
Now, we will show that there exists a constant $K>0$ such that the following bounds hold:
\begin{subequations}\label{eq:6.14}
\begin{align}
||{\bf h}(t)||_{H_{x}^{s}}&\leq K||{\bf h}(0)||_{H_{x}^{s}},\label{eq:6.14a}\\
||{\bf h}(0)||_{H_{x}^{s}}-\e^{b}K||{\bf h}(0)||_{H_{x}^{s+1}}\leq&||{\bf h}(t)||_{H_{x}^{s}}\leq 
||{\bf h}(0)||_{H_{x}^{s}}+\e^{b}K||{\bf h}(0)||_{H_{x}^{s+1}}, \quad b\in(0,1).\label{eq:6.14b}
\end{align}
\end{subequations}
The (\ref{eq:6.14}) imply the linear stability of the solution.

\noindent
Recalling that $\TT_{3}f(t):=f(t+ \al(\oo t))=f(\tau)$ and $\TT_{3}^{-1}f(\tau)=f(\tau+\widehat{\al}(\oo \tau))=f(t)$, fixing $\tau_{0}=\al(0)$, 
one has,
\begin{equation}\label{eq:6.15}
\begin{aligned}
||{\bf h}(t)||_{H_{x}^{s}}&\stackrel{(\ref{eq:6.11})}{=}
||\TT_{1}\TT_{2} \TT_{3} \TT_{4}\Phi_{\infty}{\bf v}(t)||_{H_{x}^{s}}\stackrel{(\ref{eq:6.12})}{\leq}C(s)
|| \TT_{3} \TT_{4}\Phi_{\infty}{\bf v}(t)||_{H_{x}^{s}}=||\TT_{4}\Phi_{\infty}{\bf v}(\tau)||_{H_{x}^{s}}\\
&\stackrel{(\ref{eq:6.12})}{\leq}C(s)||{\bf v}(\tau)||_{H_{x}^{s}}\stackrel{(\ref{1.2.7})}{=}
C(s)||{\bf v}(\tau_0)||_{H_{x}^{s}}
\stackrel{(\ref{eq:6.11})}{=}C(s)
||\Phi_{\infty}^{-1}\TT_{4}^{-1} {\TT_{3}}^{-1} {\TT_{2}}^{-1}\TT_{1}^{-1}{\bf h}(\tau_0)||_{H_{x}^{s}}\\
&\stackrel{(\ref{eq:6.12})}{\leq}C(s)|| {\TT_{3}}^{-1}\TT_{2}^{-1}\TT_{1}^{-1}{\bf h}(\tau_0)||_{H_{x}^{s}}
=C(s)|| {\TT_{2}}^{-1}\TT_{1}^{-1}{\bf h}(0)||_{H_{x}^{s}}\stackrel{(\ref{eq:6.12})}{\leq}C(s)||{\bf h}(0)||_{H_{x}^{s}},
\end{aligned}
\end{equation}
%
Then (\ref{eq:6.14a}) is proved.
Following the same procedure, we obtain
\begin{equation}\label{eq:6.16}
\begin{aligned}
||{\bf h}(t)||_{H_{x}^{s}}&\stackrel{(\ref{eq:6.11})}{=}
||\TT_{1} \TT_{2}\TT_{3}\TT_{4} \Phi_{\infty}{\bf v}(t)||_{H_{x}^{s}}\leq
|| \TT_{3}\TT_{4}\Phi_{\infty}{\bf v}(t)||_{H_{x}^{s}}+
||(\TT_{1}\TT_{2}-\id)\TT_{3}\TT_{4}\Phi_{\infty}{\bf v}(t)||_{H_{x}^{s}}\\
&\stackrel{(\ref{eq:6.12})}{\leq}||{\bf v}(\tau)||_{H_{x}^{s}}+
||(\TT_{4}\Phi_{\infty}-\id){\bf v}(\tau)||_{H_{x}^{s}}+\e\g^{-1}C(s)||\TT_{4}\Phi_{\infty}{\bf v}(\tau)||_{H_{x}^{s+1}}\\
&\stackrel{(\ref{1.2.7}),(\ref{eq:6.12})}{\leq}
||{\bf v}(\tau_{0})||_{H^{s}_{x}}+\e\g^{-1}C(s)||{\bf v}(\tau_0)||_{H_{x}^{s+1}},\\
&\stackrel{(\ref{eq:6.11}),(\ref{eq:6.12})}{\leq}
||{\bf h}(0)||_{H^{s}_{x}}+\e\g^{-1}C(s)||{\bf h}(0)||_{H_{x}^{s+1}},
\end{aligned}
\end{equation}
where we used $\tau_{0}=\al(0)$ and in the last inequality we have performed the same
triangular inequalities used
in the first two lines only with the $\TT_i^{-1}$. Then, using that $\g=\e^{a}$, with $a\in(0,1)$, we get
the second of (\ref{eq:6.14b}) with $b=1-a$. The first is obtained in the same way.
This concludes the proof of Theorem \ref{teo1}.
\EP

\setcounter{equation}{0}
\appendix 
\section{General Tame and Lipschitz estimates}

Here we want to illustrate some standard estimates for composition of functions and changes of variables that we use in the paper. We start with classical embedding, algebra, interpolation and tame estimate in Sobolev spaces $H^{s}:=H^{s}(\TTT^{d},\CCC)$ and $W^{s,\infty}:=W^{s,\infty}$, $d\geq1$.

\begin{lemma}\label{A} Let $s_{0}>d/2$. Then
\begin{enumerate}
\item[(i)] {\bf Embedding.} $||u||_{L^{\infty}}\leq C(s_{0})||u||_{s_{0}}$, $\forall \; u\in H^{s_{0}}$.
\item[(ii)] {\bf Algebra.} $||uv||_{s_{0}}\leq C(s_{0})||u||_{s_{0}}||v||_{s_{0}}$, $\forall\; u,v\in H^{s_{0}}$.
\item[(iii)] {\bf Iterpolation.} For $0\leq s_{1}\leq s\leq s_{2}$, $s=\la s_{1}+(1-\la)s_{2}$,
\begin{equation}\label{A1}
||u||_{s}\leq ||u||^{\la}_{s_{1}}||u||_{s_{2}}^{1-\la}, \quad \forall\; u\in H^{s_{2}}.
\end{equation}
Let $a,b\geq0$ and $p,q>0$. For all $u\in H^{a+p+q}$ and $v\in H^{b+p+q}$ one has
\begin{equation}\label{A2}
||u||_{a+p}||v||_{b+q}\leq ||u||_{a+p+q}||v||_{b}+||u||_{a}||v||_{b+p+q}.
\end{equation}
Similarly, for the $|u|^{\infty}_{s }:=\sum_{|\al|\leq s}||D^{\al}u||_{L^{\infty}}$ norm, one has
\begin{equation}\label{A3}
|u|^{\infty}_{s }\leq C(s_{1},s_{2})(|u|^{\infty}_{s_{1}})^{\la}(|u|_{s_{2}}^{\infty})^{1-\la}, 
\quad \forall\; u\in W^{s_{2},\infty},
\end{equation}
and $\forall\; u\in W^{a+p+q,\infty}, \; v\in W^{b+p+q,\infty}$,
\begin{equation}\label{A4}
|u|^{\infty}_{a+p }|v|^{\infty}_{b+q }\leq C(a,b,p,q)(
|u|^{\infty}_{a+p+q }|v|^{\infty}_{b }+|u|^{\infty}_{a }|v|^{\infty}_{b+p+q }).
\end{equation}
\item[(iv)] {\bf Asymmetric tame product.} For $s\geq s_{0}$ one has
\begin{equation}\label{A5}
||uv||_{s}\leq C(s_{0})||u||_{s}||v||_{s_{0}}+C(s)||u||_{s_{0}}||v||_{s}, \quad \forall\; u,v\in H^{s}.
\end{equation}
\item[(v)] {\bf Asymmetric tame product in $W^{s,\infty}$.} For $s\geq0$, $s\in \NNN$ one has
\begin{equation}\label{A6}
|uv|^{\infty}_{s }\leq\frac{3}{2}||u||_{L^{\infty}}|v|^{\infty}_{s }+C(s)|u|^{\infty}_{s }||v||_{L^{\infty}},
\forall\; u,v\in W^{s,\infty}.
\end{equation}
\item[(vi)] {\bf Mixed norms asymmetric tame product.} For $s\geq0$, $s\in\NNN$ one has
\begin{equation}\label{A7}
||uv||_{s}\leq\frac{3}{2}||u||_{L^{\infty}}||v||_{s}+C(s)|u|_{s,\infty}||v||_{0},
\forall\; u\in W^{s,\infty}, v\in H^{s}.
\end{equation}
If $u:=u(\la)$ and $v:=v(\la)$ depend in a lipschitz way on $\la\in\Lambda\subset\RRR$,
all the previous statements hold if one replace the norms
$||\cdot||_{s}$, $|\cdot|^{\infty}_{s }$ with $||\cdot||_{s,\g}$, $|\cdot|^{\infty}_{s,\g }$.
\end{enumerate}
\end{lemma}

Now we recall classical tame estimates for composition of functions.

\begin{lemma}\label{lemA2}{\bf Composition of functions} Let $f : \TTT^{d}\times B_{1}\to\CCC$, where
$B_{1}:=\left\{y\in\RRR^{m} : |y|<1\right\}$. it induces the composition operator on $H^{s}$
\begin{equation}\label{A11}
\tilde{f}(u)(x):=f(x,u(x),Du(x),\ldots,D^{p}u(x))
\end{equation}
where $D^{k}$ denotes the partial derivatives $\del_{x}^{\al}u(x)$ of order $|\al|=k$.

Assume $f\in C^{r}(\TTT^{d}\times B_{1})$. Then

\noindent
$(i)$ For all $u\in H^{r+p}$ such that $|u|_{p,\infty}<1$, the composition operator $(\ref{A11})$
is well defined and
\begin{equation}\label{A12}
||\tilde{f}(u)||_{r}\leq C||f||_{C^{r}}(||u||_{r+p}+1),
\end{equation}
where the constant $C$ depends on $r,p,d$. If $f\in C^{r+2}$, then, for all $|u|^{\infty}_{s },
|h|^{\infty}_{p }<1/2$, one has
\begin{equation}\label{A13}
\begin{aligned}
||\tilde{f}(u+h)-\tilde{f}(u)||_{r}&\leq C||f||_{C^{r+1}}(||h||_{r+p}+|h|^{\infty}_{p }||u||_{r+p}),\\
||\tilde{f}(u+h)-\tilde{f}(u)-\tilde{f}'(u)[h]||_{r}&\leq
C||f||_{C^{r+2}}|h|^{\infty}_{p}(||h||_{r+p}+|h|^{\infty}_{p }||u||_{r+p}).
\end{aligned}
\end{equation}

\noindent
$(ii)$ the previous statement also hold replacing $||\cdot||_{r}$ with the norm $|\cdot|_{\infty}$.
\end{lemma}

\prova For the proof see \cite{Ba2} and \cite{Moser-Pisa-66}.
\EP

\begin{lemma}\label{lemA3}{\bf Lipschitz esitmate on parameters} Let $d\in\NNN$, 
$d/2<s_{0}\leq s$, $p\geq0$, $\g>0$. Let $F : \Lambda\times H^{s}\to\CCC$, for 
$\Lambda\subset\RRR$, be a $C^{1}-$map in $u$
 satisfying the tame esitimates: $\forall\; ||u||_{s_{0}+p}\leq1$, $h\in H^{s+p}$,
 \begin{subequations}\label{A14}
 \begin{align}
 ||F(\la_{1},u)-F(\la_{2},u)||_{s}&\leq C(s) |\la_{1}-\la_{2}|(1+||u||_{s+p}), \quad 
 \la_{1},\la_{2}\in\Lambda\label{A14aa}\\
 \sup_{\la\in\Lambda}||F(\la,u)||_{s}&\leq C(s)(1+||u||_{s+p}),\label{A14a}\\
 \sup_{\la\in\Lambda}||\del_{u}F(\la,u)[h]||_{s}&\leq C(s)
 (||h||_{s+p}+||u||_{s+p}||h||_{s_{0}+p})\label{A14b}.
 \end{align}
\end{subequations}
Let $u(\la)$ be a Lipschitz family of functions with
$||u||_{s_{0}+p,\g} \leq1$. Then one has
\begin{equation}\label{A15}
||F(\cdot,u)||_{s,\g}\leq C(s)(1+||u||_{s+p,\g}).
\end{equation}
The same statement holds when the norms $||\cdot||_{s}$ are replaced by $|\cdot|^{\infty}_{s }$.
\end{lemma}
\prova We first note that, by (\ref{A14a}), one has $sup_{\la}||F(\la,u(\la))||_{s}
\leq C(s)(1+||u||_{s+p,\g})$. Then, denoting $h=u(\la_{2})-u(\la_{1})$, we have
\begin{equation}\label{A16}
\begin{aligned}
\!\!||F(\la_{2},u(\la_{2}))-&F(\la_{1},u(\la_{1}))||_{s}\leq
||F(\la_{2},u(\la_{2}))-F(\la_{1},u(\la_{2}))||_{s}+||F(\la_{1},u(\la_{2}))-F(\la_{1},u(\la_{1}))||_{s}\\
&\leq |\la_{2}-\la_{1}|C(1+||u(\la_{2})||_{s+p})+
\int_{0}^{1}||\del_{u}F(u(\la_{1})+t(u(\la_{2})-u(\la_{1})))[h]||_{s}d t\\
&\stackrel{(\ref{A14b})}{\leq} 
 C(s)\left[||h||_{s+p}+||h||_{s_{0}+p}\int_{0}^{1}
\left((1-t)||u(\la_{1})||_{s+p}+t||u(\la_{2})||_{s+p}
\right)d t
\right]\\
&+|\la_{2}-\la_{1}|C(1+||u(\la_{2})||_{s+p}),
\end{aligned}
\end{equation}
so that
\begin{equation}\label{A17}
\begin{aligned}
\g\sup_{\substack{\la_{1},\la_{2}\in\Lambda \\ \la_{1}\neq\la_{2}}}&
\frac{||F(u(\la_{1},\la_{1}))-F(\la_{2},u(\la_{2}))||_{s}}{|\la_{1}-\la_{2}|}
\leq C\g (1+\sup_{\la_{2}\in\Lambda}||u(\la_{2})||_{s+p})\\
&+C(s)\left[
||u||_{s+p,\g}+||u||_{s_{0}+p,\g} \frac{1}{2}\sup_{\la_{1},\la_{2}}(||u(\la_{1})||_{s+p,\g}+||u(\la_{2})||_{s+p,\g})
\right]\\
&\leq C(s)\left[||u||_{s+p,\g}^{2}+||u||_{s_{0}+p,\g}||u||_{s+p,\g}\right]+C(s)(1+||u||_{s+p,\g}),
\end{aligned}
\end{equation}
since $||u||_{s_{0}+p,\g}\leq1$, then the lemma follows, because
\begin{equation}\label{A17bis}
\sup_{\la\in\Lambda}||F(\la,u(\la))||_{s}+
\g\sup_{\substack{\la_{1},\la_{2}\in\Lambda \\ \la_{1}\neq\la_{2}}}
\frac{||F(u(\la_{1},\la_{1}))-F(\la_{2},u(\la_{2}))||_{s}}{|\la_{1}-\la_{2}|}
\leq C(s)(1+||u||_{s+p,\g}).
\end{equation}
\EP

In the following we will show some estimates on changes of variables. The lemma is classical,
one can see for instance \cite{Ba2}.

\begin{lemma}\label{change}{\bf (Change of variable)} Let $p : \RRR^{d}\to \RRR^{d}$
be a $2\pi-$periodic function in $W^{s,\infty}$, $s\geq1$, with $|p|^{\infty}_{1 }\leq1/2$.
Let $f(x)=x+p(x)$. Then one has
$(i)$ $f$ is invertible, its inverse is $f^{-1}(y)=g(y)=y+q(y)$ where $q$ is $2\pi-$periodic,
$q\in W^{s,\infty}(\TTT^{d};\RRR^{d})$ and $|q|^{\infty}_{s }\leq C|p|^{\infty}_{s }$. More precisely,
\begin{equation}\label{A18}
|q|_{L^{\infty}}=|p|_{L^{\infty}}, \; |dq|_{L^{\infty}}\leq 2|dp|_{L^{\infty}},
\; |dq|^{\infty}_{s-1 }\leq C|dp|^{\infty}_{s-1 },
\end{equation}
where the constant $C$ depends on $d,s$. 

Moreover, assume that $p=p_{\la}$ depends in a Lipschitz way by a parameter 
$\la\in\Lambda\subset\RRR$, an suppose, as above, that $|d_{x}p_{\la}|_{L^{\infty}}\leq1/2$
for all $\la$. Then $q=q_{\la}$ is also Lipschitz in $\la$, and
\begin{equation}\label{A19}
|q|^{\infty}_{s,\g}\leq C\left(|p|^{\infty}_{s,\g}+
\left[\sup_{\la\in\Lambda}|p_{\la}|^{\infty}_{s+1 }\right]
|p|_{L^{\infty},\g}
\right)\leq C|p|^{\infty}_{s+1,\g},
\end{equation}
the constant $C$ depends on $d,s$ (it is independent on $\g$).

$(ii)$ If $u\in H^{s}(\TTT^{d};\CCC)$, then $u\circ f(x)=u(x+p(x))\in H^{s}$, and, with the same $C$ as
in $(i)$ one has
\begin{subequations}\label{A20}
\begin{align}
||u\circ f||_{s}&\leq C(||u||_{s}+|dp|^{\infty}_{s-1 }||u||_{1}),\label{A20a}\\
||u\circ f-u||_{s}&\leq C(|p|_{L^{\infty}}||u||_{s+1}+|p|^{\infty}_{s }||u||_{2}),\label{A20b}\\
||u\circ f||_{s,\g}&\leq 
C(|u|_{s+1,\g}+|p|^{\infty}_{s,\g}||u||_{2,\g})\label{A20c}.
\end{align} 
\end{subequations}
The (\ref{A20a}), (\ref{A20b}) and (\ref{A20c}) hold also for $u\circ g$.

$(iii)$ Part $(ii)$ also holds with $||\cdot||_{s}$ replaced by $|\cdot|^{\infty}_{s }$, and
$||\cdot||_{s,\g}$ replaced by $|\cdot|^{\infty}_{s,\g}$, namely
\begin{subequations}
\begin{align}
|u\circ f|^{\infty}_{s }&\leq C(|u|^{\infty}_{s }+|dp|^{\infty}_{s-1 }|u|^{\infty}_{1 }),\label{A21a}\\
|u\circ f|^{\infty}_{s,\g}&\leq C(|u|^{\infty}_{s+1,\g}+
|dp|^{\infty}_{s-1,\g}|u|^{\infty}_{2,\g})\label{A21b}.
\end{align}
\end{subequations}
\end{lemma}

\begin{lemma} \label{lem5} {\bf (Composition).} Assume that
for any $||u||_{s_{0}+\mu_{i},\g}\leq1$ the operator $\calQ_{i}(u)$ satisfies
\begin{equation}\label{A35}
||\calQ_{i}h||_{s,\g}\leq C(s)(||h||_{s+\tau_{i},\g}+||u||_{s+\mu_{i},\g}||h||_{s_{0}+\tau_{i}\g}),
\quad i=1,2.
\end{equation}
Let $\tau:=\max\{\tau_{1},\tau_{2}\}$, and $\mu:=\max\{\mu_{1},\mu_{2}\}$. Then,
for any
\begin{equation}\label{A36}
||u||_{s_{0}+\tau+\mu,\g}\leq1,
\end{equation}
one has that the composition operator $\calQ:=\calQ_{1}\circ \calQ_{2}$ satisfies
\begin{equation}\label{A37}
||\calQ h||_{s,\g}\leq C(s)(||h||_{s+\tau_{1}+\tau_{2},\g}+||u||_{s+\tau+\mu,\g}||h||_{s_{0}+\tau_{1}+\tau_{2},\g}).
\end{equation}
\end{lemma}

\prova It is sufficient to apply the estimates (\ref{A35}) to $\calQ_{1}$ first, then to $\calQ_{2}$ and 
using the condition (\ref{A36}).
\EP

%
\setcounter{equation}{0}
\section{Proof of Lemmata \ref{capperi} and \ref{lemma6.13} }\label{C}

\noindent
\emph{ Proof of Lemma \ref{capperi}}.  
\noindent
  Define the function $\psi: \Lambda\to \CCC$,
\begin{equation}\label{eq150}
\begin{aligned}
\psi(\la)&:=i\la\bar{\oo}\cdot\ell+\Omega_{\s,j}(\la)-\Omega_{\s',j'}(\la)
\stackrel{(\ref{eq:4.4})}{=}
i\la\bar{\oo}\cdot\ell-i m(\la)(\s j^{2}-\s' j'^{2})+r_{\s,j}^{\infty}(\la)-r_{\s',j'}^{\infty}(\la),
\end{aligned}
\end{equation}
where with abuse of notation we set $r^{\infty}_{\s,0}\equiv0$.
Note that, by $(N1)_{n}$ of Theorem \ref{NM}, we have $||{\bf u}_{n}||_{\gots_{0}+\mu_{2},\g} \leq1$ on $G_{n}$. Then (\ref{eq:4.5}) holds and we have
\begin{equation}\label{eq151}
|\Omega_{\s,j}-
\Omega_{\s',j'}|^{{\rm lip}}\leq|m|^{{\rm lip}}|\s j^{2}-\s'j'^{2}|+|r_{\s,j}^{\infty}|^{{\rm lip}}
+|r_{\s',j'}^{\infty}|^{{\rm lip}}\leq C\e\g^{-1}|\s j^{2}-\s' j'^{2}|\stackrel{(\ref{eq147})}{\leq}C \e\g^{-1}|\bar{\oo}\cdot\ell|.
\end{equation}
We can estimate, for any $\la_{1},\la_{2}\in\Lambda$,
\begin{equation}\label{eq152}
\frac{|\psi(\la_{1})-\psi(\la_{2})|}{|\la_{1}-\la_{2}|}
\stackrel{(\ref{eq147}), (\ref{eq151})}{\geq}
\left(\frac{1}{8}-C\e\g^{-1}\right)|\bar{\oo}\cdot\ell|\geq\frac{|\s j^{2}-\s'j'^{2}|}{9},
\end{equation}
if $\e\g^{-1}$ is small enough. Then, using standard measure estimates on sub-levels of Lipschitz functions, we conclude
\begin{equation}\label{eq153}
|R_{\ell j j'}^{\s,\s'}|\leq 4\g_{n}|\s j^{2}-\s' j'^{2}|\langle\ell\rangle^{-\tau}\frac{9}{|\s j^{2}-\s' j'^{2}|}\leq
C \g\langle\ell\rangle^{-\tau}.
\end{equation}
\EP

\noindent
\emph{Proof of Lemma \ref{lemma6.13}}.
\noindent
We first prove the (\ref{eq147}); note that if $(\s,j)=(\s',j')$ then it is trivially true.
If $R_{\ell j j'}^{\s,\s'}({\bf u}_{n})\neq\emptyset$, then, by definition (\ref{eq145}), there exists a $\la\in\Lambda$
such that
\begin{equation}\label{eq158}
|\Omega_{\s,j}({\bf u}_{n})-\Omega_{\s',j'}({\bf u}_{n})|<2\g_{n}|\s j^{2}-\s'j'^{2}|\langle\ell\rangle^{-\tau}+2|\bar{\oo}\cdot\ell|.
\end{equation}
On the other hand, for $\e$ small and since $(\s,j)\neq(\s',j')$,
\begin{equation}\label{eq159}
|\Omega_{\s,j}({\bf u}_{n})-\Omega_{\s',j'}({\bf u}_{n})|
\stackrel{(\ref{eq:3.2.4a}),(\ref{eq:4.5})}{\geq}
\frac{1}{2}|\s j^{2}-\s'j'^{2}|-C\e\geq\frac{1}{3}|\s j^{2}-\s' j'^{2}|.
\end{equation}
By the (\ref{eq158}), (\ref{eq159}) and $\g_{n}\leq 2\g$ follows
\begin{equation}\label{eq160}
2|\bar{\oo}\cdot\ell|\geq\left(\frac{1}{3}-\frac{4\g}{\langle\ell\rangle^{\tau}}\right)|\s j^{2}-\s'j'^{2}|
\geq \frac{1}{4}|\s j^{2}-\s' j'^{2}|,
\end{equation}
since $\g\leq \g_{0}$, by choosing $\g_{0}$ small enough. It is sufficient $\g_{0}<1/48$. Then, the (\ref{eq147}) hold.

\noindent
In order to prove the (\ref{eq146}) we need to understand the variation of the eigenvalues $\Omega_{\s,j}({\bf u})$ with respect to the function ${\bf u}$. 
We have to study the difference
\begin{equation}\label{eq161}
\begin{aligned}
\Omega_{\s,j}({\bf u}_{n})-&\Omega_{\s,j}({\bf u}_{n-1})=
-i(m({\bf u}_{n})-m({\bf u}_{n-1}))(\s j^{2}-\s'j'^{2})+(r_{\s,j}^{\infty}({\bf u}_{n})-r_{\s,j}^{\infty}({\bf u}_{n-1}))
\end{aligned}
\end{equation}
Indeed if we assume that
\begin{equation}\label{eq167}
|(\Omega_{\s,j}-\Omega_{\s',j'})({\bf u}_{n})-(\Omega_{\s,j}-\Omega_{\s',j'})({\bf u}_{n-1})|
\leq C \e |\s j^{2}-\s'j'^{2}|N_{n}^{-\al},
\end{equation}
then, for $j\neq j'$, $|\ell|\leq N_{n}$, and $\la\in G_{n}$, we have
\begin{equation}\label{eq168}
\begin{aligned}
|i\la\bar{\oo}\cdot\ell+\Omega_{\s,j}({\bf u}_{n})-\Omega_{\s',j'}({\bf u}_{n})|
&\stackrel{(\ref{eq167})}{\geq}2\g_{n-1}|\s j^{2}-\s' j'^{2}|\langle\ell\rangle^{-\tau}
-C\e|\s j^{2}-\s'j'^{2}|N_{n}^{-\al}\\
&\geq2\g_{n}|\s j^{2}-\s'j'^{2}|\langle\ell\rangle^{-\tau},
\end{aligned}
\end{equation}
because $C\e\g^{-1}N_{n}^{\tau-\al}2^{n+1}\leq1$ if $\e\g^{-1}$ small enough. In order to complete the proof of
the (\ref{eq146}) we only need to verify the  (\ref{eq167}).

By Lemma \ref{teo:KAM}, using  the $({\bf S4})_{n+1}$ with
$\g=\g_{n-1}$ and $\g-\rho=\g_{n}$,  and with ${\bf u}_{1}={\bf u}_{n-1}$, ${\bf u}_{2}={\bf u}_{n}$,
we have
\begin{equation}\label{eq162}
\Lambda_{n+1}^{\g_{n-1}}({\bf u}_{n-1})\subseteq\Lambda_{n+1}^{\g_{n}}({\bf u}_{n}),
\end{equation}
since, for $\e\g^{-1}$ small enough,
\begin{equation}\label{eq163}
\e C N_{n}^{\tau}\sup_{\la\in G_{n}}||{\bf u}_{n}-{\bf u}_{n-1}||_{\gots_{0}+\mu}\stackrel{(\ref{eq109})}{\leq}
\e^{2}\g^{-1}C C_{\star} N_{n}^{\tau-\ka_{3}}\leq \g_{n-1}-\g_{n}=:\rho=\g2^{-n}.
\end{equation}
where $\ka_{3}$ is defined in (\ref{eq109}) with
$\nu=2$, $\mu$ defined in (\ref{eq:4.4.18}) with $\h=\h_{1}+\be$, $\mu>\tau$ 
(see Lemmata \ref{measure}, \ref{megalemma} and (\ref{eq:4.14}), (\ref{eq:3.2.0})).
 We also  note that,
\begin{equation}\label{eq164}
G_{n}\stackrel{(\ref{eq142bis}), (\ref{martina10})}{\subseteq}\Lambda_{\infty}^{2\g_{n-1}}({\bf u}_{n-1})
\stackrel{(\ref{eq:4.1.12})}{\subseteq}\cap_{\nu\geq0}\Lambda_{\nu}^{\g_{n-1}}({\bf u}_{n-1})
\subseteq\Lambda_{n+1}^{\g_{n-1}}({\bf u}_{n-1})\stackrel{(\ref{eq162})}{\subseteq}
\Lambda_{n+1}^{\g_{n}}({\bf u}_{n}).
\end{equation}
This means that 
$\la\in G_{n}\subset \Lambda_{n+1}^{\g_{n-1}}({\bf u}_{n-1})\cap\Lambda_{n+1}^{\g_{n}}({\bf u}_{n})$,
and hence, we can apply the ${\bf (S3)}_{\nu}$, with $\nu=n+1$, in Lemma \ref{teo:KAM} to get
\begin{equation}\label{eq165}
\begin{aligned}
|r_{\s,j}^{\infty}({\bf u}_{n})-&r_{\s,j}^{\infty}({\bf u}_{n-1})|\leq
|r_{\s,j}^{n+1}({\bf u}_{n})-r_{\s,j}^{n+1}({\bf u}_{n-1})|+
|r_{\s,j}^{\infty}({\bf u}_{n})-r_{\s,j}^{n+1}({\bf u}_{n})|
+
|r_{\s,j}^{\infty}({\bf u}_{n-1})-r_{\s,j}^{n+1}({\bf u}_{n-1})|\\
&\stackrel{(\ref{eq:4.1.10}),(\ref{eq:3.2.6a}),(\ref{aaa})}{\leq}
\e C ||{\bf u}_{n}-{\bf u}_{n-1}||_{\gots_{0}+\h_{2}}
+\e\left(1+||{\bf u}_{n-1}||_{\gots_{0}+\h_{1}+\be}+||{\bf u}_{n}||_{\gots_{0}+\h_{1}+\be}\right)N_{n}^{-\al}\\
&\stackrel{(\ref{eq109})}\leq C\e^{2}\g^{-1} N_{n}^{-\mu_{3}}
+\e\left(1+||{\bf u}_{n-1}||_{\gots_{0}+\h_{1}+\be}+||{\bf u}_{n}||_{\gots_{0}+\h_{1}+\be}\right)N_{n}^{-\al}.
\end{aligned}
\end{equation}
Now, first of all $\mu_{3}>\al$ by (\ref{eq109}), (\ref{eq:4.14}), moreover
$\h_{1}+\be<\h_{5}$ then by ${\bf (S1)}_{n}$, $({\bf S1})_{n-1}$, one has 
$||{\bf u}_{n-1}||_{\gots_{0}+\h_{5}}+||{\bf u}_{n}||_{\gots_{0}+\h_{5}}\leq 2$, we obtain
\begin{equation}\label{eq166}
|r_{\s,j}^{\infty}({\bf u}_{n})-r_{\s,j}^{\infty}({\bf u}_{n-1})|\stackrel{(\ref{eq165})}{\leq}
\e N_{n}^{-\al}.
\end{equation}
Then, by (\ref{eq161}), (\ref{eq:3.2.4b}) and (\ref{eq166}) one has that the (\ref{eq167}) hold and the proof of Lemma (\ref{lemma6.13}) is complete.
\EP

This research was supported by the European Research Council under
FP7
and partially by the PRIN2009
grant ``Critical point theory and perturbative methods for
nonlinear differential equations".

\end{document}